\documentclass[12pt]{article}
\usepackage{amsmath,amssymb,amsthm,cite}
\usepackage[curve,frame,line,matrix,arrow,graph,poly]{xy}
\usepackage{color}

\advance\topmargin by -\headsep \textheight=9.5in
\font \eightrm=cmr8

\newcommand{\delete}[1]{}
\newcommand{\invlim}{\varprojlim}

\DeclareMathOperator{\Ad}{Ad}       
\renewcommand{\a}{\alpha}           
\newcommand{\A}{\mathcal{A}}        
\DeclareMathOperator{\ad}{ad}       
\newcommand{\alg}{\mathrm{alg}}     
\renewcommand{\b}{\beta}            
\newcommand{\barox}{\mathrel{\overline\otimes}} 
\newcommand{\C}{\mathbb{C}}         
\DeclareMathOperator{\CBHD }{CBHD}  
\newcommand{\D}{\mathbb{D}}         
\DeclareMathOperator{\Der}{Der}     
\newcommand{\Dl}{\Delta}            
\DeclareMathOperator{\Diff}{Diff}   
\newcommand{\del}{\partial}         
\newcommand{\ddto}[1]{\frac{d}{d#1}\biggr|_{#1=0}} 
\newcommand{\diff}{\mathfrak{diff}} 
\newcommand{\dl}{\delta}            
\newcommand{\dst}[2]{\langle#1,#2\rangle} 
\DeclareMathOperator{\End}{End}     
\newcommand{\eps}{\varepsilon}      
\newcommand{\F}{\mathcal{F}}        
\def\fd#1#2{\frac{d#1}{d#2}}        
\newcommand{\G}{\mathcal{G}}        
\newcommand{\Ga}{\Gamma}            
\newcommand{\g}{\mathfrak{g}}       
\newcommand{\ga}{\gamma}            
\DeclareMathOperator{\gr}{gr}       
\DeclareMathOperator{\Hom}{Hom}     
\newcommand{\hookto}{\hookrightarrow} 
\newcommand{\id}{{\mathrm{id}}}     
\renewcommand{\L}{\mathcal{L}}      
\newcommand{\la}{\lambda}           
\DeclareMathOperator{\Map}{Map}     
\newcommand{\nn}{\nonumber}         
\newcommand{\Om}{\Omega}            
\newcommand{\op}{\oplus}            
\newcommand{\ox}{\otimes}           
\newcommand{\pd}{\partial}          
\newcommand{\pds}[2]{\frac{\partial#1}{\partial#2}} 
\newcommand{\Q}{\mathbb{Q}}         
\DeclareMathOperator{\Qder}{Qder}   
\newcommand{\R}{\mathbb{R}}         
\newcommand{\Rr}{\mathcal{R}}       
\newcommand{\Sf}{\mathbb{S}}        
\DeclareMathOperator{\Sh}{Sh}       
\newcommand{\sepword}[1]{\quad\mbox{#1}\quad} 
\newcommand{\set}[1]{\{\,#1\,\}}    
\newcommand{\sg}{\sigma}            
\newcommand{\sll}{\mathfrak{sl}}    
\DeclareMathOperator{\T}{T}         
\newcommand{\Te}{\mathcal{T}}       
\newcommand{\Th}{\Theta}            
\newcommand{\thalf}{\tfrac{1}{2}}   
\newcommand{\tthird}{\tfrac{1}{3}}  
\newcommand{\U}{\mathcal{U}}        
\newcommand{\vf}{\varphi}           
\newcommand{\X}{\mathfrak{X}}       
\newcommand{\x}{\times}             
\renewcommand{\.}{\cdot}            
\renewcommand{\:}{:\,}              
\let\hash=\#
\renewcommand{\#}{\mathbin{\hash}}  
\def\<#1,#2>{\langle#1,#2\rangle}   
\def\wick:#1:{\mathopen:#1\mathclose:} 

\makeatletter
\def\section{\@startsection{section}{1}{\z@}{-3.5ex plus -1ex minus
               -.2ex}{2.3ex plus .2ex}{\large\bf}}
\def\subsection{\@startsection{subsection}{2}{\z@}{-3.25ex plus -1ex
               minus -.2ex}{1.5ex plus .2ex}{\normalsize\bf}}
\makeatother

\numberwithin{equation}{section}    

\theoremstyle{plain}
\newtheorem{thm}{Theorem}[section]  
\newtheorem{prop}[thm]{Proposition} 
\newtheorem{lema}[thm]{Lemma}       
\newtheorem{corl}[thm]{Corollary}   

\theoremstyle{definition}
\newtheorem{defn}{Definition}       

\theoremstyle{remark}

\begin{document}

\title{Hopf algebras in dynamical systems theory}

\author{Jos\'e F. Cari\~nena*, \
Kurusch Ebrahimi-Fard\dag,\ \\[6pt]
H\'ector Figueroa\ddag\
\ and \
Jos\'e M. Gracia-Bond\'{\i}a\S$\|$\
\\[20pt]
*Departamento de F\'{\i}sica Te\'orica, Universidad de Zaragoza,
\\
Zaragoza 50006, Spain
\\[6pt]
\dag Max-Planck-Institut f\"ur Mathematik, Vivatsgasse 7
\\
D-53111 Bonn, Germany
\\[6pt]
\ddag Departamento de Matem\'aticas, Universidad de Costa Rica,
\\
San Pedro 2060, Costa Rica
\\[6pt]
\S Departamento de F\'{\i}sica Te\'orica I,
Universidad Complutense,
\\
Madrid 28040, Spain
\\[6pt]
$\|$Departamento de F\'{\i}sica, Universidad de Costa Rica,
\\
San Pedro 2060, Costa Rica}

\date{December 28, 2006 }

\maketitle

\begin{abstract}
The theory of exact and of approximate solutions for non-autonomous
linear differential equations forms a wide field with strong ties to
physics and applied problems. This paper is meant as a stepping stone
for an exploration of this long-established theme, through the tinted
glasses of a (Hopf and Rota--Baxter) algebraic point of view. By
reviewing, reformulating and strengthening known results, we give
evidence for the claim that the use of Hopf algebra allows for a
refined analysis of differential equations. We revisit the renowned
Campbell--Baker--Hausdorff--Dynkin formula by the modern approach
involving Lie idempotents. Approximate solutions to differential
equations involve, on the one hand, series of iterated integrals
solving the corresponding integral equations; on the other hand,
exponential solutions. Equating those solutions yields identities
among products of iterated Riemann integrals. Now, the Riemann
integral satisfies the integration-by-parts rule with the Leibniz rule
for derivations as its partner; and skewderivations generalize
derivations. Thus we seek an algebraic theory of integration, with the
Rota--Baxter relation replacing the classical rule. The methods to
deal with noncommutativity are especially highlighted. We find new
identities, allowing for an extensive embedding of Dyson--Chen series
of time- or path-ordered products (of generalized integration
operators); of the corresponding Magnus expansion; and of their
relations, into the unified algebraic setting of Rota--Baxter maps and
their inverse skewderivations. This picture clarifies the approximate
solutions to generalized integral equations corresponding to
non-autonomous linear (skew)differential equations.
\end{abstract}

\medskip


\smallskip

\noindent
\textsl{2006 PACS}: 02.20.Qs, 02.30.Hq, 02.40.Gh, 45.30.+s

\smallskip


\bigskip

\rightline{\it To Ennackal Chandy George Sudarshan,}
\rightline{\it with admiration, on his 75th birthday}

\tableofcontents
\newpage

\section{Aim, plan of the article and preliminaries}

This paper studies non-autonomous differential equations of the
general type
\begin{align}
\dot g(t)g^{-1}(t) &= \xi(t),  \qquad g(t_0) = 1_G \sepword{or}
\label{eq:gato-encerrado}
\\
g^{-1}(t)\dot g(t) &= \eta(t), \qquad g(t_0) = 1_G,
\label{eq:menos-gatos}
\end{align}
where the unknown $g:\R_t\to G$ is a curve on a (maybe
infinite-dimensional) local Lie group~$G$, with~$1_G$ the neutral
element; and $\xi(t),\eta(t)$ are given curves on the tangent Lie
algebra~$\g$ of~$G$. Before proceeding note that, if $g(t)$
solves~\eqref{eq:gato-encerrado}, then $g^{-1}(t)$
solves~\eqref{eq:menos-gatos} for $\eta(t)=-\xi(t)$.

Such equations are pervasive in mathematics, physics and engineering.
To begin with, $G$ can have a faithful finite-dimensional
representation. For instance, consider (affine) linear differential
equations on~$\R^n$,
\begin{equation}
\dot x = A(t)x + b(t) \sepword{with} x(t_0) = x_0.
\label{eq:manes-de-afin}
\end{equation}
They are exactly solved by
$$
x(t) = \G(t,t_0)x_0 + \int_{t_0}^t\G(t,t')b(t')\,dt';
$$
where the \textit{Green's function}~$\G(t,t_0)$ is the matrix
satisfying
$$
\frac{d\G(t,t_0)}{dt}\,\G^{-1}(t,t_0) = A(t), \quad \G(t_0,t_0)=1_n
$$
of the kind~\eqref{eq:gato-encerrado}. More generally, dynamical
systems admitting a superposition principle can be reduced to the
form~\eqref{eq:gato-encerrado}, with~$G$ a finite dimensional Lie
group. This assertion is part of the classical \textit{Lie--Scheffers}
theory~\cite{NihilNovumSubSole}, reviewed in Section~3 as part of and
motivation for the whole enterprise. Even more generally, any
non-autonomous dynamical system, given in local coordinates by
\begin{equation}
\frac{dx^i}{dt} = Y^i(t;x(t)), \qquad x = (x^1,\ldots,x^n)
\label{eq:alpha-and-omega}
\end{equation}
corresponds to a `time-dependent vector field'~$Y$ with~$Y(t)$
belonging to $\X(M)$, the Lie algebra of all vector fields on a
manifold~$M$. Then the solution of~\eqref{eq:alpha-and-omega} is given
by the solution of an equation like~\eqref{eq:menos-gatos}; this
remark will be formalized in Section~4. The crucial difference is
dimensionality of the (pseudo-)group. In practice, almost always we
must content ourselves with approximate solutions ---Lie--Scheffers
systems are not solvable by quadratures in general--- and actually
those are our main concern.

We reformulate~\eqref{eq:gato-encerrado} within the framework of Hopf
algebra and Rota--Baxter operator theory ---the latter has become
popular recently in relation with the Connes--Kreimer paradigm for
renormalization theory in perturbative quantum field theory. The
convenience of such algebraic approach stems already from that,
unless~$G$ is a matrix group, $\dot g$ and~$g^{-1}$ cannot be
multiplied, and then equations~\eqref{eq:gato-encerrado}
and~\eqref{eq:menos-gatos} have no meaning, strictu sensu. Hopf
algebras generalize both Lie groups and Lie algebras, so the problem
does not present itself in a Hopf algebra formalism. Another advantage
is that Hopf algebra and Rota--Baxter theory allow for efficient and
meaningful comparisons among the different techniques for
solving~\eqref{eq:gato-encerrado}, proposed over the years. The main
aim of this paper is to show these and other benefits of our algebraic
viewpoint. They have been patent for a while to people working on
control theory, but mostly ignored by the wider community of
mathematicians and mathematical physicists.

We presume the readers acquainted with standard tools of
differential analysis, like for instance in~\cite{OldRedBook}:
primarily the notion of tangent map and the exponential
map~$\exp:\g\to G$. For the benefit of the readers, Lie group and
Lie algebra actions are reviewed in Appendix~A. The basics of Hopf
algebra are a prerequisite. Unless otherwise specified, we
consider Hopf algebras over the complex numbers. We briefly
introduce our notations for them, which are like in~\cite{Quaoar};
the pedagogical paper~\cite{FresasSalvajes} is recommended as
well. Given an associative algebra with unit~$H\ni1_H=:u(1_\C)$,
then $H\ox H$ is associative with bilinear multiplication given by
$(a\ox b)(c\ox d)=ac\ox bd$ on decomposable tensors, and
unit~$1_H\ox1_H$. Write just~1 for the unit element of~$H$
henceforth. One says $H$ is a bialgebra if algebra morphisms
$\eta:H\to\C$ (augmentation) and $\Dl:H\to H\ox H$ (coproduct) are
defined, such that the maps $(\eta\ox\id)\Dl$ and
$(\id\ox\eta)\Dl$ from~$H$ to~$H$ coincide with the identity
map~$\id$ and $(\Dl\ox\id)\Dl$ and $(\id\ox\Dl)\Dl$ from~$H$
to~$H\ox H\ox H$ also coincide (we omit the sign $\circ$ for
composition of linear maps). One says $H$ is a Hopf algebra if it
furthermore possesses an antiautomorphism~$S$, the
\textit{antipode}, such that $m(S\ox\id)\Dl=m(S\ox\id)\Dl=u\eta$;
where $m:H\ox H\to H$ denotes the algebra map. Familiarity with
enveloping algebras and the Poincar\'e--Birkhoff--Witt and
Cartier--Milnor--Moore theorems in particular will be helpful.
Both results can be summarized in the statement that a connected
cocommutative Hopf algebra is the enveloping algebra of a Lie
algebra, as an algebra, and cofree, as a coalgebra. At any rate,
we discuss a strong version of the Poincar\'e--Birkhoff--Witt
theorem in Section~4, and the Cartier--Milnor--Moore theorem at
the end of Section~5. The necessary notions of Rota--Baxter
operator theory will be introduced and explained in due course.

Whereas the Hopf algebraic description springs up naturally from the
intrinsic geometrical approach, we have found it expedient to smooth
this transition with the help of Lie--Rinehart algebras: these
constitute the ``noncommutative geometry'' version of Lie algebroid
technology.

The plan of the work is as follows. In this section we explain our
main aims and fix some notations of frequent use. Next we recommend a
look at Appendix~A, indispensable for everything that follows; most
readers will just need to scan it for the notations. In Section~2 we
address for the first time Lie--Scheffers systems; they are intimately
linked to equations~\eqref{eq:gato-encerrado}
and~\eqref{eq:menos-gatos}.

After this, two paths are possible: either reading Appendices~C and~D
for motivation, or not. Most of the stuff in them could be regarded as
preceding Sections~3 and following; but it gets in the way of our
algebraic business, and this is why it has been confined to the end.
Section~3 plunges the reader at once into an application of Hopf
algebra to differential geometry. This is due to Rinehart and
Huebschmann, and deserves to be better known, as it clarifies several
questions; one should compare the treatment of differential operators
given here with that of~\cite[Chapter~3]{Mackenzie}. Readers less
familiar with Hopf algebra might wish to read this in parallel with
Section~4.

Sections~4 to 6 are largely expository. In Section~4 we leisurely
\textit{construct} the Hopf algebra structure governing our approach
to equations~\eqref{eq:gato-encerrado} and~\eqref{eq:menos-gatos} from
the geometric notions. Section~5 recalls some structure results for
Hopf algebras.

After that, our master plan is to transplant the usual paradigmatic
strategies for dealing with~\eqref{eq:gato-encerrado}
and~\eqref{eq:menos-gatos} to the Hopf algebraic soil, which on the
one hand will prove to be their native one, and on the other naturally
leads to far-reaching generalizations. At the outset, in Section~6 we
consider the Campbell--Baker--Hausdorff--Dynkin (CBHD) development,
which we proceed to derive in Hopf algebraic terms. In turn, that
development is the natural father of the Magnus expansion
method~\cite{Magnus}. We eventually derive the Magnus series with the
help of Rota--Baxter theory. For the purpose, the Riemann integral is
treated in this paper as a particular Rota--Baxter operator of
weight~zero.

We first show in Section~7 that skewderivations and Rota--Baxter
operators (of the same weight) are natural inverses. The ordinary
Spitzer formula is revisited in Section~8, together with a nonlinear
CBHD recursion due to one of us. The latter is instrumental in
obtaining the \textit{noncommutative Spitzer formula}. Also, inspired
by the work of Lam, we give a noncommutative generalization of the
Bohnenblust--Spitzer formula. In the next two sections, the Magnus
expansion is arrived at as a limiting case of that formula. All along,
we try to distinguish carefully which statements are valid for general
Rota--Baxter operators, which for Rota--Baxter operators of vanishing
weight, and which just for the Riemann integral. The main alternative
integration method, the Dyson--Chen
`expansional'~\cite{Dyson,Chen,Araki}, flows from the Magnus series,
and vice versa, by our Hopf algebraic means in Section~11. In turn, it
reveals itself useful to understand the quirks of the Magnus
expansion, and to solve the weight-zero CBHD recursion. Section~12
explores by means of pre-Lie algebras with Rota--Baxter maps the
solution of that nonlinear recursion in the general case. Section~13
is the conclusion, whereupon perspectives for research are discussed.

As said, Appendix~A reviews the basics of Lie group and Lie algebra
actions on manifolds. As also hinted at, Appendices~B and~C run a
parallel, complementary strand to the main body of the paper. They
contain more advanced material on the topic of dynamical systems with
symmetry; their treatment here naturally calls for the Darboux
derivative of Lie algebroid theory. The main point is to show how one
is led to the arena of Lie algebra and \textit{geometrical
integration}, for the solution of differential equations we are
concerned with. This provides a rationale for our choice of the Magnus
series, and its generalizations, as the primary approximation method
in the body of the paper ---see the discussion at the beginning of
Section~10.

Appendix~D gives the Hopf algebraic vision of a theorem of~Lie
and~Engel.

\smallskip

Several notational conventions are fixed next. Let~$M$ be a (second
countable, smooth, without boundary) manifold of finite dimension~$n$.
The space~$\F(M)$ of (real or) complex smooth functions on~$M$ is
endowed with the standard commutative and associative algebra
structure. Let $\tau_M:TM\to M$ denote the tangent bundle to~$M$.
Vector fields on~$M$ can be defined either as sections for~$\tau_M$,
that is, maps $X:M\to TM$ such that $\tau_M\circ X=\id_M$, or as
derivations of~$\F(M)$. When we wish to distinguish between those
roles, we denote by~$\L_X$ the differential operator corresponding to
the vector field~$X$. Because the commutator of two derivations is
again a derivation, the space~$\X(M)$ of vector fields has a Lie
algebra structure; we choose to define the bracket there as the
opposite of the usual one: in local coordinates,
$$
[X, Y]^i = Y^j \,\del_j X^i - X^j \,\del_j Y^i;
$$
so $\X(M)\equiv\diff(M)$, the Lie algebra of the infinite-dimensional
Lie group~$\Diff(M)$~\cite{Omori}. Also $\X(M)$ becomes a faithful
$\F(M)$-module when one defines $hX(x)=h(x)X(x)$ for~$h\in\F(M)$.
Together, $\F(M)$ and~$\X(M)$ constitute a Lie--Rinehart algebra in
the sense of~\cite{Bigotes1}.

Given a smooth map $f:N\to M$, the pull-back $f^*:\F(M)\to\F(N)$ is
defined as $f^*h=h\circ f$. A vector field $X\in\X(N)$ is said to be
$f$-related with the vector field $Y\in\X(M)$ if $Tf\circ X=Y\circ f$;
we then say that the vector field $X$ is $f$-projectable onto the
vector field $Y$, and write $X\sim_f Y$. We have $X\sim_f Y$ iff the
maps $\L_X\circ f^*$ and~$f^*\circ\L_Y$ from~$\F(M)$ to~$\F(N)$
coincide. If~$X_1,X_2\in\X(N)$ are $f$-related with~$Y_1,Y_2$
respectively, then $X_1+X_2$ and~$[X_1,X_2]$ are also $f$-related,
respectively with $Y_1+Y_2$ and~$[Y_1,Y_2]$. A given~$X\in\X(N)$ will
not be $f$-projectable in general. However, if~$f$ is a
diffeomorphism, then every vector field $X\in\X(N)$ is projectable
onto a unique vector field on~$M$, to wit, $Y=Tf\circ X\circ f^{-1}$,
and we say~$X$ is the pull-back of~$Y$. A vector field $X\in\X(M)$ is
\textit{invariant} under a diffeomorphism $f$ of~$M$ iff $X\sim_f X$.
On~$\R_t$ (or on an open interval~$I\subset\R_t$) there is a canonical
vector field~$d/dt$. A curve $\gamma:\R_t\to M$ is said to be an
integral curve for a vector field $X\in\X(M)$ if~$d/dt$ and $X$ are
$\gamma$-related: $\dot\gamma:=T\gamma\circ d/dt=X\circ \gamma$.
Well-known theorems assert that the integral curves of a vector field
define a local $\R_t$-action or flow~\cite{Carniceria}.

By a \textit{vector field along}~$f$ we understand a map $Y:N\to
TM$ such that $\tau_M\circ Y=f$:
\[
\xymatrix{&TM\ar[d]^{\tau_M} \\ N\ar[r]_f\ar[ur]^Y &M}
\]
It is clear that the concept is just a particular case of a more
general one: section along the map $f$ over a general bundle $\pi:E\to
M$. Vector fields along~$f$ can also be regarded as $f$-derivations,
in an obvious sense. The right hand side of the non-autonomous
dynamical system~\eqref{eq:alpha-and-omega} is just the vector field
along the map $\pi_2:\R_t\x M\to M$ expressed in local coordinates by
$$
Y = Y^i(t;x(t))\,\pd_i\circ\pi_2.
$$
Also, clearly any curve $\gamma:\R_t\to M$ defines a vector
field~$\dot\gamma$ along~$\gamma$. We envisage here the concept of
integral curves of vector fields~$Y$ along maps~$f$. These are curves
$\gamma:\R_t\to N$ such that the image under $Tf\circ T\gamma$ of the
vector field~$d/dt$ coincides with the vector field
along~$f\circ\gamma$ given by~$Y\circ\gamma$ ---depending on~$f$,
there might be vector fields along it without integral curves.
Under this definition $t\mapsto(t,\ga(t))$ is always the integral
curve of~$(t,\ga(t))\mapsto\dot\gamma(t)$ along~$\pi_2:\R_t\x M\to M$.

\section{The Lie--Scheffers theorem}

\begin{defn}
The system~\eqref{eq:alpha-and-omega} of differential equations admits
a \textit{superposition principle} ---or possesses a set of
fundamental solutions--- if a superposition function
$\Psi:\R^{n(m+1)}\to\R^n$ exists, written
\begin{equation}
x = \Psi(x_{(1)},\ldots,x_{(m)};k_1,\ldots,k_n),
\end{equation}
such that the general solution of~\eqref{eq:alpha-and-omega} can be
expressed (at least for small~$t$) as the functional
\begin{equation}
x(t) = \Psi(x_{(1)}(t),\dots,x_{(m)}(t);k_1,\dots,k_n),
\label{eq:manes-de-LS}
\end{equation}
where~$\set{x_{(a)}:a=1,\dots,m}$ is a set of particular solutions and
$k_1,\dots,k_n$ denote~$n$ arbitrary parameters. The latter must be
\textit{essential} in the sense that they can be solved from the
solution functional:
\begin{equation}
k = \Xi\bigl(x_{(1)}(t),\dots,x_{(m)}(t);x(t)\bigr), \sepword{with} k
:= (k_1,\dots,k_n).
\label{eq:desmanes-de-LS}
\end{equation}
\end{defn}

The Lie--Scheffers theorem~\cite{NihilNovumSubSole} asserts:
\begin{thm}
{}For~\eqref{eq:alpha-and-omega} to admit a superposition principle it
is necessary and sufficient that the time-dependent vector field~$Y$
be of the form
\begin{equation}
Y(t;x) = Z_1(t)X_1(x) + \cdots + Z_r(t)X_r(x),
\label{eq:cascabel-al-gato}
\end{equation}
where, as indicated in the notation, the $r$~scalar functions~$Z_a$
depend only of~$t$ and the $r$~vector fields~$X_a$ depend only on the
variables~$x$; and these fields close to a real Lie algebra~$\g$. That
is, the~$X_a$ are linearly independent and there exist suitable
structure constants~$f_{ab}^c$ such that
\begin{equation*}
[X_a, X_b] = \sum_{c=1}^r f_{ab}^c\,X_c.
\end{equation*}
Moreover, the dimension~$r$ of~$\g$ is not greater than $nm$. Systems
fulfilling the conditions of the theorem are called here
Lie--Scheffers systems associated to~$\g$.
\end{thm}

Modern reviews of this subject
include~\cite{macrote,TresMosqueteros,Ciguenya,PBPepinArturo,
TresMosqueterosbis}. We ponder a few pertinent examples of
Lie--Scheffers systems next.

\smallskip

Linear systems~\eqref{eq:manes-de-afin} are Lie--Scheffers systems. If
$n+1$ particular solutions $x_{(1)},\ldots,x_{(n+1)}$
of~\eqref{eq:manes-de-afin} are known, such that
$x_{(2)}(t)-x_{(1)}(t),\ldots,x_{(n+1)}(t)- x_{(1)}(t)$ are
independent, and $H(t)$ is the regular matrix with these vectors as
columns, then the vector of parameters~\eqref{eq:desmanes-de-LS} is
given by
$$
k = H^{-1}(t)\bigl(x(t) - x_{(1)}(t)\bigr).
$$
This follows from the fact that the transformation
$$
x'(t) = H^{-1}\bigl(x(t) - x_{(1)}(t)\bigr)
$$
reduces the system to $dx'/dt=0$.

A famous example for $n=1$ is provided by the Riccati equation:
\begin{equation}
\dot x = a_0(t) + a_1(t)x + a_2(t)x^2.
\label{eq:manes-de-Riccati}
\end{equation}
One can understand by Hopf algebraic methods why, up to
diffeomorphisms, Riccati's is the only nonlinear Lie--Scheffers
differential equation on the real line; this was indicated
in~\cite{Bruna} and it is spelled in Appendix~D.
Also~\eqref{eq:manes-de-Riccati} is the simplest Lie--Scheffers system
not solvable by quadratures; and other Lie--Scheffers equations on the
line are reductions of it, in an appropriate sense ---see Appendix~B.
The superposition principle for the Riccati equation is given by
\begin{equation*}
k = \frac{(x - x_{(2)})(x_{(1)} - x_{(3)})}{(x - x_{(1)})(x_{(2)} -
x_{(3)})}.
\end{equation*}

The one-dimensional example
\begin{equation}
\dot x = b(t)\chi(x)
\label{eq:tugurito}
\end{equation}
is instructive. We assume $\chi$ does not change sign in the interval
of interest. Let
$$
\phi(x) := \int^x\frac{dx'}{\chi(x')}.
$$
Then
$$
x(t) = \phi^{-1}\Bigl(k' + \int b(t)\,dt\Bigr)
$$
is the general solution. We have therefore a superposition rule of the
form
$$
x(t) = \phi^{-1}\bigl(\phi(x_{(1)}(t)) + k\bigr),
$$
with $m=n=r=1$: only one particular solution is required. Notice that
the local diffeomorphism~$\phi$ projects the vector field
corresponding to the right hand side in~\eqref{eq:tugurito} to the
vector field~$b(t)\pd_x$: in the language of Lie, the flow associated
to this problem is locally similar to a translation.

\smallskip

Now we can tackle at last the question of giving intrinsic geometrical
meaning to~\eqref{eq:gato-encerrado}. It turns out to correspond to a
Lie--Scheffers system on a Lie group. For an arbitrary curve $t\mapsto
g(t)$ on the $r$-dimensional Lie group~$G$, we have the vector field
along the curve given by $\dot g(\cdot)$ as discussed at the end of
Section~1. Then we \textit{define} the left hand side
of~\eqref{eq:gato-encerrado} as
\begin{equation}
\dot g(t)\,g^{-1}(t) := T_{g(t)}R_{g^{-1}(t)}\dot g(t).
\label{eq:madre-del-cordero}
\end{equation}
By construction, for each value of the parameter~$t$, this vector lies
in $T_1G\equiv\g$, the tangent Lie algebra of~$G$. We obtain in this
way a curve on~$\g$. Note that, if $g(t)=\exp(t\eta)$, then simply
$\dot g(t)\,g^{-1}(t)=\eta$. Now, if $\ga^a$ is a basis for~$\g$, then
\begin{equation}
\dot g(t)\,g^{-1}(t) = \sum_{a=1}^r Z_a(t)\ga^a =: \xi(t) \in \g,
\label{eq:eqingr}
\end{equation}
for some functions~$Z_a(t)$. We realize that $g(t)$ is an integral
curve of the right invariant vector field along $\R_t\x G\to G$:
$$
\xi_G(t,g) := \sum_{a=1}^r Z_a(t)\ga^a_G(g) = \sum_{a=1}^r
Z_a(t)X^R_{\ga^a}(g), \sepword{with} \ga_G^a(1_G) = \ga^a.
$$
Here $\ga^a_G$ is the fundamental vector field or infinitesimal
generator of the left group translations generated by~$\ga^a$,
which is a \textit{right} invariant vector field ---see Appendix~A
for the notations. In the language of~\eqref{eq:alpha-and-omega},
the differential system is
\begin{equation}
\dot g(t) = \sum_{a=1}^r Z_a(t)\,\ga^a_G(g(t)).
\label{eq:fund}
\end{equation}
The theorem says that for every choice of the functions $Z_a(t)$
we have a (right invariant) Lie--Scheffers system on the Lie
group~$G$, and any such system is of this form. The reader should
be aware, nevertheless, that for a system of the
type~\eqref{eq:fund} there might be more than one superposition
rule~\cite{TresMosqueterosbis}. The reason one needs only one
particular solution is precisely the right invariance of the last
equation: if $g(t)$ is the solution such that $g(t_0)=1$, then
consider $\bar g(t):=g(t)\,g_0$ for each $g_0\in G$. We have
\begin{eqnarray*}
&&\dot{\bar g}(t) = T_{g(t)}R_{g_0}(\dot g(t)) =
T_{g(t)}R_{g_0}\left[T_1R_{g(t)}\left(\sum_{a=1}^r
Z_a(t)\ga^a\right)\right]
\cr
&&= T_1R_{g(t)g_0}\left(\sum_{a=1}^r Z_a(t)\ga^a\right) = \sum_{a=1}^r
Z_a(t)\ga^a_G(\bar g(t)).
\end{eqnarray*}
Equation $X^R(g) = T_1R_gX^R(1_G)$ has been used. Therefore $\bar
g(t)$ is the solution of the same equation~\eqref{eq:fund} with $\bar
g(t_0)=g_0$: the solution curves of the system are obtained from just
one of them by right-translations. In other words, the superposition
functional~\eqref{eq:manes-de-LS} can be symbolically expressed by
$\Psi(g_{(1)},k)=g_{(1)}k$, with $k\in G$; for which always~$m=1$.

Lie--Scheffers systems live on manifolds which are not groups in
general; however, they are always associated with the action of a
finite-dimensional Lie group on the manifold on which $Y(t;x)$ is
defined; and this symmetry of the differential equation can be
powerfully exploited through the action of the group of curves on
the group manifold on a set of systems of the same type. This
variant of Lie's reduction method is of wide applicability; it is
explained in Appendix~B.

Let us finally note than in control theory, say on~$M\equiv\R^n$,
business is often with equations of a form not unrelated
to~\eqref{eq:fund}:
\begin{equation*}
\dot x(t) = X_1(x(t)) + \sum_{a=2}^r Z_a(t)\,X_a(x(t));
\end{equation*}
the functions $Z_2,\ldots,Z_r$ being the controls. In the most
important cases the $X_1,\ldots,X_r$ vector fields close to a
finite-dimensional Lie algebra, or $X_2,\ldots,X_r$ close to a
finite-dimensional Lie algebra.

\section{Differential operators on Lie--Rinehart algebras}

Let $R$ be a commutative, unital ring and $\A$ a commutative algebra
over $R$ be given. A derivation $\dl$ of~$\A$ is a $R$-linear map
from~$\A$ to itself, such that $\dl(ab)=\dl a\,b+a\,\dl b$. Since $\A$
is commutative, the linear space $\Der(\A)$ of such maps becomes an
$\A$-module when we define $(a\dl)b=a\,\dl b$. Moreover, with the
usual bracket given by the commutator $\Der(\A)$ is a Lie algebra. In
this paper $R=\C$ nearly always.

A left (right) \textit{action} of a Lie algebra~$\g$ on~$\A$ is a Lie
algebra homomorphism (antihomomorphism) $\a:\g\to\Der(\A)$.

\begin{defn}
Assume that we are given a commutative algebra~$\A$ and a Lie algebra
$\g$ which is also a faithful $\A$-module. The pair $(\A,\g)$ is a
\textit{Lie--Rinehart algebra} if there exists an $\A$-module morphism
$\a:\g\to\Der(\A)$, called the \textit{anchor}, satisfying the
compatibility condition
\begin{equation}
a[X, Y] = [X, a\,Y] - \a(X)a\,Y,
\label{eq:compatibility}
\end{equation}
for $a\in\A, X,Y\in \g$. If we write $m_a$ for the multiplication
operator $m_a(X)=a\,X$ and~$\ad_X$, as usual, for the adjoint operator
$\ad_X(Y)=[X, Y]$, the compatibility condition is rewritten as
$$
[\ad_X, m_a] = m_{\a(X)a}.
$$
Often, in the definition of Lie--Rinehart algebra, the apparently
stronger condition that the anchor be also a left action of~$\g$
on~$\A$ is required. But, as it turns out, these two definitions are
equivalent.
\end{defn}

The concept essentially coincides with Kastler and Stora's Lie--Cartan
pairs~\cite{GoodPair}. As already indicated
$\bigl(\F(M),\X(M),\id\bigr)$ is a Lie--Rinehart algebra. A more
general example of Lie--Rinehart algebra may be $\bigl(\F(M),
\Ga(M,E),\a\bigr)$, where $\Ga(M,E)$ is the $\F(M)$-module of sections
of a vector bundle~$E$ over~$M$, on which a Lie bracket and an
anchor~$\a$ (hence a vector bundle map $E\to TM$, denoted in the same
way) are supposed given. If the fibres have dimension bigger than~one,
then a linear map satisfying~\eqref{eq:compatibility} is not only
automatically a Lie algebra morphism, but also a $\F(M)$-module
morphism. These geometrical examples are called Lie algebroids. A Lie
algebroid is called \textit{transitive} when it is onto fibrewise;
totally intransitive when $\a=0$. For examples of this, consider a
principal bundle $P(M,G,\pi)$ over~$M$; if~$VP$ is the vertical bundle
over~$P$, we have the exact sequences of vector bundles
$$
0 \to VP   \hookrightarrow TP   \to TM \to 0 \sepword{and}
0 \to VP/G \hookrightarrow TP/G \to TM \to 0;
$$
the second being essentially the Atiyah sequence; and then
$(C^\infty(M),\Ga(M,TP/G),T\pi/G)$ is a transitive Lie algebroid;
while $\bigl(C^\infty(M),\Ga(M,VP/G)\bigr)$ is totally intransitive.

Whenever we have a Lie--Rinehart algebra, we can algebraically define
a differential calculus. For instance a $n$-form is a skewsymmetric
$n$-linear map from~$\g$ to~$\A$. If we define~$d$ on 1-forms by
$$
d\b(X,Y) = \a(X)\b(Y) - \a(Y)\b(X) - \b\bigl([X, Y]\bigr),
$$
then certainly $d$ can be extended so~$d^2=0$. Also, let $V$ be an
$\A$-module. A $V$-connection in the sense of~\cite{GoodPair,Bigotes2}
is a linear assignment to each element $X\in \g$ of a linear
map~$\rho(X):V\to V$ such that, for $v\in V$,
$$
\bigl(a\rho(X)\bigr)v = a\bigl(\rho(X)\bigr)v; \qquad \rho(X)(av)
= a\,\rho(X)v + \a(X)a\,v.
$$
If~$V$ is moreover a $\g$-module, the connection is flat (as the
curvature defined in the obvious way vanishes).

A morphism $(\A,\g)\to(\A',\g')$ of Lie--Rinehart algebras is a pair
$(\phi,\psi)$ of an algebra morphism $\phi:\A\to\A'$ and an
$\A$-module morphism $\psi:\g\to \g'$, where the action of~$\A$ on~$\g'$
is given by $aX':=\phi(a)X'$, intertwining the anchors:
$$
\phi\bigl(\a(X)\,a\bigr) = \a'\bigl(\psi(X)\bigl)\,\phi(a).
$$
An important example by Grabowski is as follows~\cite{JanuszBifronte}.
A linear operator $D:\g\to\g$ is called a \textit{quasi-derivation}
for~$\A$, and we write $D\in\Qder_\A(\g)$, if for each~$a\in\A$ there
exists $\widehat D(a)\in\A$ ---necessarily unique--- such that
$[D,m_a]=m_{\widehat D(a)}$, where the bracket is the usual
commutator. It is easily seen that $\widehat D\in\Der(\A)$. Then
$(\id_\A,\ad):(\A,\g)\to\bigl(\A,\Qder_\A(\g)\bigr)$, where
$\ad:X\mapsto\ad_X$, is a morphism of Lie--Rinehart algebras.

Our first example of Hopf algebra comes now across: the (universal)
enveloping algebra~$\U(\g)$ of the Lie algebra~$\g$. The enveloping
algebra is Hopf because there is the diagonal algebra
homomorphism
$$
\Dl: \U(\g) \to \U(\g\op \g) \simeq \U(\g) \ox \U(\g) \sepword{by} X
\mapsto X \op X \mapsto X \ox 1 + 1 \ox X, \; \Dl1 = 1 \ox 1,
$$
for every~$X\in \g$.

When $\g=\X(M)$, the enveloping algebra should not be confused with
the algebra (with the usual composition product) of differential
operators~$\D(M)$. It is true that the first-order elements of both
are the vector fields, thus coincide. However, we are going to show
that if first order differential operators are to be considered
primitive elements, then $\D(M)$ cannot be given a natural Hopf
algebra structure; we thank P.~Aschieri for making us aware of the
following argument, since published~\cite{Corfu}. Consider the linear
map $\L:\U(\X(M))\to\D(M)$ obtained by extending the Lie derivative.
This map is not onto because zeroth order differential operators are
functions, whereas the zeroth order elements of $\U(\X(M))$ are just
scalars. Actually, $\D(M)$ is a $\F(M)$-module, while $\U(\X(M))$ is
not. For this very reason $\D(M)$ cannot be a Hopf algebra. Consider
two commuting linearly independent vector fields $X,Y$ nonvanishing on
a common domain (e.g. locally let~$X$ be the partial
derivative~$\del_i$ and~$Y$ a different one~$\del_j$), and the vector
fields $aX,aY$, where $a$ is an arbitrary function, nonvanishing on
the same domain. The composition $a\L_X\L_Y = a\L_Y\L_X$ is an element
in~$\D(M)$. Suppose \textit{arguendo} that there exist on~$\D(M)$ a
coproduct~$\dl$ compatible with composition of operators, and such
that vector fields are primitives. We would have
\begin{align*}
\dl(a\L_X\L_Y) &= a\L_X\L_Y \ox 1 + a\L_X \ox \L_Y + \L_Y \ox a\L_X + 1
\ox a\L_X\L_Y \sepword{and}
\\
\dl(a\L_Y\L_X) &= a\L_Y\L_X \ox 1 + a\L_Y \ox \L_X + \L_X \ox a\L_Y +
1 \ox a\L_Y\L_X.
\end{align*}
Now, the right hand sides are not equal. As a corollary we have
that the map $\L$ is not injective: the fact that in~$\U(\X(M))$
one has a good coproduct implies that $aX\. Y$ is there different
>from $aY\. X$, with $\.$ the product in the enveloping algebra.
Thus the kernel of $\L$ contains $aX\. Y-aY\. X$. Notice that the
argument fails if~$M$ is one-dimensional. Notice as well that
$\dl$ makes $\D(M)$ into a good coalgebra over~$\F(M)$. However,
the product is then not $\F(M)$-linear.

In spite of the above, Hopf algebra renders us a first great service
in helping to manufacture~$\D(M)$ out of~$\U(\X(M))$. This involves a
purely algebraic construction~\cite{Bigotes1,Harbinger} suggested by
the previous discussion and better presented in the context of
Lie--Rinehart algebras.

Assume for simplicity that $\A$ is unital. The \textit{universal
object} of~$(\A,\g)$ is by definition a triple $(\U(\A,\g),\imath_\A,
\imath_\g)$, where $\U(\A,\g)$ is an associative algebra, therefore a
Lie algebra with the usual commutator, together with morphisms
$\imath_\A:\A\to\U(\A,\g)$ and $\imath_\g:\g\to\U(\A,\g)$,
respectively of algebras and Lie algebras, such that
$$
\imath_A(a)\imath_\g(X) = \imath_\g(aX); \qquad
[\imath_\g(X),\imath_A(a)] = \imath_A(\a(X)a);
$$
and $(\U(\A,\g),\imath_\A,\imath_\g)$ is universal among these
triples; that is, for a similar triple $(B,\phi_\A,\phi_\g)$, there is
a unique algebra morphism $\Phi_B:\U(\A,\g)\to B$ such that
$\Phi_B\imath_\A=\phi_\A$ and $\Phi_B\imath_\g=\phi_\g$.

To construct~$\U(\A,\g)$ we employ~$\U(\g)$ in the following way. The
condition that the anchor maps into derivations means precisely
that~$\A$ is a Hopf $\U(\g)$-module~\cite{Quaoar}. One may keep using
the same notation~$\a$ for the new action, as for the generators;
$\a(1)=1$. Consider now the smash product or crossed product algebra
$\A\rtimes\U(\g)$. This is the vector space $\A\ox\U(\g)$ with the
product defined on simple tensors by
$$
(a \ox u)(b \ox v) :=  a\a(u_{(1)})b \ox u_{(2)}v,
$$
where we use the standard Sweedler notation $\Dl u=u_{(1)}\ox
u_{(2)}$. There are obvious morphisms $\imath'_\A:\A\to\A\rtimes
\U(\g)$ and $\imath'_\g:\g\to\A\rtimes\U(\g)$. Now, let~$J$ be the
right ideal generated in~$\A\rtimes U(\g)$ by the elements $ab\ox X -
a\ox bX$. One has
\begin{align*}
(c \ox X)(ab \ox Y - a \ox bY) &= cab \ox XY + c\a(X)(ab) \ox Y - ca
\ox XbY - c\a(X)a \ox bY
\\
&= cab \ox XY + c\a(X)a\,b \ox Y + ca\a(X)b \ox Y
\\
&\quad -ca \ox bXY - ca \ox \a(X)b\,Y - c\a(X)a \ox bY
\\
&= cab\ox XY - ca \ox bXY + c\a(X)a\,b \ox Y - c\a(X)a \ox bY
\\
&\quad + ca\a(X)b \ox Y - ca \ox \a(X)b\,Y,
\end{align*}
where in the last equality we just reordered terms; hence $J$ is
two-sided. By construction it is clear that the quotient
$\U(\A,\g):=\A\rtimes U(\g)/J$ together with the obvious quotient
morphisms $\imath_\A$ and $\imath_\g$ possesses the universal
property. Note that $\imath_\A$ is injective. A morphism
$(\phi,\psi):(\A,\g)\to(\A',\g')$ induces a morphism of algebras
$\U(\phi,\psi):\U(\A,\g)\to\U(\A',\g')$, and vice versa; this is an
equivalence of categories.

One obtains by this construction the ordinary algebra of differential
operators $\D(M)=\U(\F(M),\X(M))$. It is also clear that
$\U(\C,\g)=\U(\g)$ with trivial action of~$\g$. Just like the
enveloping algebra, the universal algebra $\U(\A,\g)$ is filtered,
with an associated graded object $\gr\,\U(\A,\g)$, which is a
commutative graded $\A$-algebra. There is also a
Poincar\'e--Birkhoff--Witt theorem for $\U(\A,\g)$ when~$\g$ is
projective over~$\A$ ---which is the case in the geometrical examples.
It claims that if $S_\A[\g]$ is the symmetric $\A$-algebra on~$\g$,
then the natural surjection
\begin{equation}
S_\A[\g] \to \gr\,\U(\A,\g),
\label{eq:Casimir-future}
\end{equation}
is an isomorphism of $\A$-algebras, rather like the $S[\g]\simeq\gr\,
\U(\g)$ effected, say, through the `symmetrization' map $\sigma:S[\g]
\to\U(\g)$. In that case $\imath_\g$ is of course injective.

\section{Coming by Hopf algebra}

We begin here a journey from the geometrical to the Hopf world. Let us
start by some well-known observations~\cite{KawskiS}. A smooth
manifold~$M$ is determined by the linear space~$\F(M)$, in the sense
that points of~$M$ are in one-to-one correspondence with a particular
class of linear functionals on~$\F(M)$, to wit, multiplicative ones.
One writes $\<x,h>:=h(x)$ to express this correspondence. As a
consequence $M\hookto\F'(M)$, where $\F'(M)$ is the space of compactly
supported distributions on~$M$. We denote by $\C M$ the subspace
of~$\F'(M)$ generated by the points of~$M$. It will sometimes be
convenient to write $Th=\<T,h>$ for the value of the distribution
$T\in\F'(M)$ at $h\in\F(M)$; accordingly we abbreviate to~$xh=h(x)$.
Many other geometrical objects can be expressed as functionals in this
way; for instance, if $v_x\in T_xM$, then $\<v_x,h>\equiv v_xh$ is the
derivative of~$h$ in the direction of the tangent vector~$v_x$ at~$x$.
Therefore $TM\hookto\F'(M)$, too. If $X$ is a smooth vector field and
$Xh:=X(h)$, then $xX$ is defined naturally by
$$
(xX)h = x(Xh).
$$
That is, $xX=X(x)$; and we may omit the parentheses in~$xXh$. An
advantage of thinking in this way is that operations in principle not
meaningful on~$M$ make sense in~$\F'(M)$. For instance, given a curve
$\ga:\R_t\to M$ with $\ga(0)=x_0$, the definition
$$
\dot\ga(0) = \lim_{\eps\downarrow0}\frac{\ga(\eps) - x_0}{\eps},
$$
which looks unacceptable on~$M$, in~$\F'(M)$ just means that for all
$h\in\F(M)$:
$$
\dot\ga(0)h := \bigl<\lim_{\eps\downarrow0}\frac{\ga(\eps) -
x_0}{\eps}, h\bigr> := \lim_{\eps\downarrow0}\frac{h(\ga(\eps)) -
h(x_0)}{\eps}.
$$
Let now $f:N\to M$ be smooth. If $S\in\F'(N)$, a corresponding
element~$S_f$ is defined in~$\F'(M)$ by
$$
S_f h := S(f^*h).
$$
Clearly, $x_f=xf=f(x)$, and so $S\to S_f$ extends~$f$ to a map
from~$\F'(N)$ to~$\F'(M)$. Somewhat rashly, one denotes the extension
by the same letter; with this notation, if $v_x$ is a tangent vector
at~$x\in N$, the tangent vector $T_xf(v_x)$ at~$f(x)\in M$ would
become~$v_x f$.

We may freely use the notation $e^{tX}$ for the flow generated by a
vector field~$X$: if $\ga_X(t,x_0)$ is the integral curve of~$X$ going
through $x_0$ at $t=0$ and $xe^{tX}:=\ga_X(t,x)$, then the identity
$\fd{}{t}\bigl(\ga_X(t,x)\bigr)= X\bigl(\ga_X(t,x)\bigr)$ acquires the
linear look $\fd{}{t}\bigl(xe^{tX}\bigr) = xe^{tX}X$. We have indeed
linearized the dynamical system associated to~$X$. In more detail: if
$xe^{tX}h:=h\bigl(\ga_X(t,x)\bigr)$, then
\begin{align*}
&\fd{}{t}\bigl(xe^{tX}\bigr)h = \sum_{i=1}^n\fd{x^i}{t}\; \frac{\pd
h}{\pd x^i}\bigl(\ga_X(t,x)\bigr) =\sum_{i=1}^n X^i
\bigl(\ga_X(t,x)\bigr) \frac{\pd h}{\pd x^i}\bigl(\ga_X(t,x)\bigr)
\\
&=Xh\bigl(\ga_X(t,x)\bigr) =: xe^{tX}Xh.
\end{align*}
Linearization works as well for non-autonomous dynamical systems.
Recall~\eqref{eq:alpha-and-omega} under the form:
\begin{equation}
\dot x = Y(t;x(t));  \qquad x(t_0) = x_0.
\label{eq:gatos}
\end{equation}
{}For~$t$ given, the vector $Y(t;x(t))$ lives in the fibre
over~$x(t)$. Denote
$$
L(t,t_0)h(x) = h(x(t)), \sepword{for $h\in\F(M)$; then}
\frac{dL(t,t_0)}{dt} = L(t,t_0)Y(t);
$$
with~$Y$ interpreted as the corresponding time-dependent vector field.
This, as announced in Section~1, is in the guise
of~\eqref{eq:menos-gatos}. The Cauchy problem
$$
\dot x = xL(t,t_0); \qquad x(t_0) = x_0
$$
has that of~\eqref{eq:gatos} as unique
solution~\cite{RusosNoDescontrolados}. The difference with the
finite-dimensional case is of course substantial; at the analytical
level this is discussed at the end of Section~11.

\smallskip

Linearization is precisely what Hopf algebra is about. Things become
really interesting when there is a symmetry group~$G$ of the
manifold~$M$. As pointed out in~\cite{HazewinkelRecom}, linearization
is then a particularly good idea; for instance, often the
action~$\Phi$ is indecomposable (think of the case $M=G$ and lateral
action) and so contains little information; whereas the linear actions
of $G$ on~$\F(M)$ and~$\F'(M)$ are generally decomposable (for
instance when $G=M=\Sf^1$). This is the point of harmonic analysis.
For a fully algebraic description of these phenomena, we try to regard
$\F'(G)$ and~$\F(G)$, eventually restricting appropriately the
functors~$\F',\F$, as Hopf algebras.

There is no trouble in recognizing an algebra structure for the
\textit{whole} of~$\F'(G)$: this is given just by
\textit{convolution}, which is a map $\F'(G)\ox\F'(G)\hookto\F'(G\x
G)\to\F'(G)$. If $S_1,S_2\in \F'(G)$, then $S_1*S_2$ is defined as the
image $\mu(S_1,S_2)$, of the extended group multiplication $\mu:G\x
G\to G$. This is an associative operation~\cite{DieudonneIII}. We have
in particular $g_1*g_2=g_1g_2$, for $g_1,g_2\in G$. The unit element
in $\F'(G)$ is~$1_G$. An augmentation on~$\F'(G)$ is given by
evaluation on the function $1\in\F(G)$:
$$
\eta(S) = S1.
$$
In particular $\eta(g)=1$ for all~$g\in G$. Clearly $\eta(S_1*S_2)=
\eta(S_1)\eta(S_2)$. A candidate antipode is the extension of the
inversion diffeomorphism~$\imath:g\mapsto g^{-1}$; certainly it is an
algebra antiautomorphism:
$$
\imath(S_1*S_2) = \imath(S_2)*\imath(S_1).
$$
In the sequel, the integral notation for convolution
$$
(S_1*S_2)h' = \int dS_1(g')\,dS_2(g)\,h'(g'g),
$$
will be handy. A locally summable function $h$ defines a distribution
by $h'\mapsto\int h(g)h'(g)\,dg$, with $dg$ a (left) invariant measure
on~$G$. For instance $1\in\F'(G)$ if $G$ is compact. Now $S*h$ is
defined by:
$$
\int dS(g')\,h(g)h'(g'g)\,dg = \int dS(g')\,h({g'}^{-1}g)h'(g)\,dg;
$$
so we identify it with the function
\begin{equation}
g \mapsto \int dS(g')\,h({g'}^{-1}g) = \int dS(gg')\,h(g).
\label{eq:the-end-of-the-beginning}
\end{equation}
Similarly, for $h*S$:
$$
\int dS(g)\,h(g')h'(g'g)\,dg' = \int dS(g)\,h(g'g^{-1})h'(g')\,
\dl(g^{-1})\,dg';
$$
where $\dl$ is the modular function, so we identify $h*S$ with
$g'\mapsto\int dS(g)\,h(g'g^{-1})\,\dl(g^{-1})$. In particular,
$$
h_1*h_2(g) = \int_G h_1(g')h_2({g'}^{-1}g)\,dg' = \int_G
h_1(gg')h_2({g'}^{-1})\,dg' = \int h_1(g{g'}^{-1})h_2(g')\,
\dl^{-1}(g')\,dg'.
$$

To give $\F'(G)$ a coalgebra structure, one might try the
following strategy. The diagonal homomorphism $d\: G\to G\x G$,
given by $g\mapsto(g,g)$, extends to~$d\:\F'(G)\to\F'(G\x G)$, by
the repeatedly used procedure. However, $\F'(G\x G)$ is vastly
bigger than $\F'(G)\ox\F'(G)$. So we look for convolution
subalgebras $O(G)$ of~$\F'(G)$ for which $O(G\x G)\simeq O(G)\ox
O(G)$. For a start, $\C G$ will do; and naturally the elements
of~$G$, when regarded as elements of $\F'(G)$, are
\textit{grouplike} in the sense of Hopf algebra theory:
for~$h_1,h_2\in\F(G)$:
$$
\<\Dl g, h_1 \ox h_2> = \<g, \mu(h_1 \circ d \ox h_2 \circ d)>
=\<g\ox g, h_1\ox h_2>.
$$
The Hopf algebra $\C G$ is too small for our purposes.
Nonetheless, recall that the tangent algebra~$\g$ of~$G$ also sits
inside~$\F'(G)$. The discussion around~\eqref{eq:for-later-use} in
Appendix~A allows us to regard the elements of~$\g$ as right
invariant differential operators on~$G$. So let us focus on the
subalgebra $\D^R(G)\subsetneq\D(G)$ of right invariant
differential operators on~$\F(G)$. There is great advantage in
regarding any $D\in\D(G)$ as extended to distributions by
$(DS)h=D(Sh)$. Now, $\D^R(G)$ can be realized as the algebra of
distributions on~$G$ with support at~$1_G$. In effect, look first
at the fundamental vector fields~$\xi_G\equiv X^R_\xi$.
From~\eqref{eq:the-end-of-the-beginning}, we see that in general
$$
\xi_G(S*h) = \xi_G(S)*h; \sepword{as $1_G*h=h$, we conclude}
\xi_G(h) = \xi_G(1_G)*h.
$$
For any element $D$ of~$\D^R(G)$ analogously $D(h)=D(1_G)*h$. Note
that $D(1_G)$ is a distribution concentrated at~$1_G$. Moreover,
$$
DD'h = D(1_G)*D'(1_G)*h,
$$
so the map $D\mapsto D(1_G)$ is a homomorphism. It is in fact an
isomorphism, as any distribution vanishing outside a point is a finite
sum of derivatives of a Dirac function; thus conversely $D$ can be
written as a polynomial in the right invariant vector fields.

Therefore we have a new subalgebra of $\F'(G)$. Let us just
write~$\xi$ for~$\xi_G(1_G)$. Furthermore, by the Leibniz rule we are
able to define the shuffle coproduct:
\begin{equation}
\Dl\xi = \xi \ox 1 + 1 \ox \xi.
\label{eq:Aquitaine}
\end{equation}
This extends to~$\D^R(G)(1_G)$ as an algebra homomorphism.
Naturally~\eqref{eq:Aquitaine} says that the elements of $\g$,
when regarded as elements of $\F'(G)$, are primitive in the sense
of Hopf algebra theory. A little more work shows that in fact
$\D^R(G)\simeq\U(\X^R(G))\equiv\U(\g)$, the algebra of right
invariant differential operators \textit{coincides} with the
enveloping algebra of the Lie algebra of fundamental vector fields
for the left action of~$G$ on itself. Also we remark here that the
equivalence of the Lie algebra structures considered on~$T_1G$ in
Appendix~A can be seen from
$$
\xi*\eta - \eta*\xi = [\xi, \eta];
$$
see~\cite{DieudonneIV}.

Our~$O(G)$ will be the convolution algebra generated by
$\U(\g)\equiv\D^R(G)(1_G)$ and $\C G$; this is a Hopf crossed
product~\cite{Quaoar}, as $g*\xi*g^{-1} = \Ad_g\xi$, and similarly
for more general elements of~$\U(\g)$; it can be also regarded as
a completion of the latter. Note $\Dl\circ i=(i\ox i)\Dl$ as well;
we invite the reader to check the rest of the expected Hopf
algebra properties. By the way, extending~$\Ad$ to~$S(\g)$ as
well, it is found that the symmetrization map mentioned
after~\eqref{eq:Casimir-future} intertwines both actions of~$G$.
The centre~$Z(\g)$ of left \textit{and} right (Casimir) invariant
differential operators is clearly a commutative algebra,
isomorphic to the algebra of the $G$-invariant elements in~$S(\g)$
---this states a strong form of the Gelfand--Harish--Chandra
theorem.

One begins to feel the power of the Hopf algebra approach:
equations~\eqref{eq:gato-encerrado} and~\eqref{eq:menos-gatos}
make sense in~$O(G)$ without further ado; we are allowed to write
for them
\begin{equation}
\dot g(t) * g^{-1}(t) = \xi(t),        \quad g(t_0) = 1_G \sepword{or}
\dot g(t)             = \xi(t) * g(t), \quad g(t_0) = 1_G;
\label{eq:auto-gato}
\end{equation}
similarly for~\eqref{eq:menos-gatos}:
\begin{equation*}
g^{-1}(t) * \dot g(t) = \eta(t),        \quad g(t_0) = 1_G \sepword{or}
\dot g(t)             = g(t) * \eta(t), \quad g(t_0) = 1_G;
\end{equation*}
The rigorous but roundabout arguments at the end of Section~3 are
simplified thereby. Moreover the possibility of considering
interpolated equations, of the form
$$
g^{-a}(t) * \dot g(t) * g^{-b}(t) = \kappa(t), \sepword{with $a + b =
1$,}
$$
opens distinctly~\cite{Eunomia}. This is uninvestigated.

Before leaving~$\F'(G)$, let us note that we refrained from pondering
distributions with point supports other than the elements of~$G$. This
would have allowed us in particular to consider $TG\hookto\F'(G)$; and
then $g*v_{g'}=gv_{g'},\,v_{g'}*g=v_{g'}g$
---see~\eqref{eq:maestros-sutiles} in Appendix~A. However, one should
not believe $v_g*v_{g'}=v_gv_{g'}$; that is, $O(TG)$ does not embed
into~$\F'(G)$.

\smallskip

It is well known that the subalgebra $\Rr(G)$ of `representative
functions' in~$\F(G)$, with its ordinary commutative multiplication,
is also a Hopf algebra. The space~$\Rr(G)$ is made of those functions
whose translates $x\mapsto h(xt)$, for all $t\in G$, generate a
finite-dimensional subalgebra of~$\F(G)$. Then also $\Rr(G)$ is
endowed with a coproduct in which
\begin{equation}
\Dl h\in\Rr(G) \ox \Rr(G) \sepword{is given by} \Dl h(x,y) :=
\bigl(h_{(1)} \ox h_{(2)}\bigr)(x,y) := h(xy);
\label{eq:dies-illa}
\end{equation}
which is not cocommutative, unless~$G$ is abelian. One has:
$$
\eta(h) = h(1); \qquad Sh(g) = h(g^{-1}).
$$
Both previous constructions of $O(G)$ and~$\Rr(G)$ are mutually dual.
Questions of duality are delicate in Hopf algebra theory; fortunately
we need not deal with them in particular detail. The main point there
is the following. Given any Hopf algebra~$H$ and an algebra $A$, one
can define~\cite{Quaoar} the \textit{algebraic convolution} of two
$\C$-linear maps $f,h \in \Hom(H,A)$ as the map $f*h \in \Hom(H,A)$
given by the composition
$$
H \xrightarrow{\Dl} H \ox H \xrightarrow{f\ox h} A \ox A
\xrightarrow{m_A} A.
$$
Here $m_A$ denotes the product map from $A\ox A$ to~$A$. Because of
coassociativity of~$\Dl$, the triple $\bigl(\Hom(H,A),*,
u_A\eta_H=:\eta_A\bigr)$ is an associative algebra with unit. Now, for
$A=\C$, the Hopf algebraic definition of convolution on~$O(G)$ as a
dual of~$\Rr(G)$ \textit{coincides} with the analytical one.

Algebra morphisms respect convolution, in the following way
\begin{equation*}
\ell(f * h) = \ell f * \ell h; \sepword{similarly} (f * h)\ell = f\ell
* h\ell,
\end{equation*}
if $\ell$ is a coalgebra morphism. Clearly the antipode~$S$ is the
inverse of the identity map~$\id$ for the convolution product of
endomorphisms of~$H$~\cite{Quaoar}. If $f\in\Hom(H,A)$ is an algebra
morphism, using the convolution product of~$\End(H)$ one finds that
its composition~$fS$ with the antipode is a convolution inverse
for~$f$:
$$
f * fS = f(\id*S) = fu_H\eta_H = \eta_A = f(S * \id) = fS * f.
$$
Denote by $\Hom_\alg(H,A)$ the convolution monoid ---with unit
element the map~$\eta_A$--- of multiplicative morphisms of~$H$ on
the algebra $A$. In general $fS$ does not belong to
$\Hom_\alg(H,A)$; but it does when the algebra~$A$ is commutative.
Moreover, if~$A$ is commutative, the convolution product of two
multiplicative maps is again multiplicative, so $\Hom_\alg(H,A)$
becomes a group, that we may call $G_H(A)$. In particular, this
happens for the set $G_H(\C)$ of scalar characters, and
for~$G_H(H)$ if $H$ is commutative. Thus we have a (representable
by definition) functor $G_H$ going from commutative Hopf algebras
to groups. We may call $G_H$ an `affine group scheme'. If we
suppose $H$ to be graded, connected (meaning that the scalars are
the only elements in degree zero) and of finite type, then
$G_H(\C)$ is a projective limit of triangular matrix groups. An
important example is studied in Appendix~D.

In general there will be only an embedding ---that can be made
continuous--- of~$G$ into the group $G_{\Rr(G)}(\R)$ of characters
of~$\Rr(G)$; under favourable circumstances (for instance, for $G$
compact, thanks to the Peter--Weyl theorem) both groups coincide. Also
if $A$ is commutative then $\Hom(H,A)$ is an $A$-algebra. Then a Lie
algebra $\g_H(A)$ can be obtained as well by considering the
elements~$L$ (`infinitesimal characters') of~$\Hom(H,A)$ satisfying
the Leibniz rule
$$
L(cd) = \eta_A(c)Ld + \eta_A(d)Lc,
$$
for all $c,d\in H$. The bracket $[L_1,L_2]:=L_1*L_2-L_2*L_1$ of two
infinitesimal characters is an infinitesimal character, and so we have
a functor~$\g_H$ from commutative algebras to Lie algebras. Needless
to say, under favourable circumstances $\g_{\Rr(G)}(\R)$ is just~$\g$.

\smallskip

To summarize, the situation is here quite different of that examined
in Section~2, whereby we showed $\D(G)\nsim\U(\X(G))$, whereas
$\D^R(G)\sim\U(\X^R(G))$. Notice that $\D^R(G)$ can be expressed
directly in Hopf theoretic terms, as follows: a derivation of the
commutative algebra~$\Rr(G)$ belongs to $\D^R(G)$ iff it is of the
form $L*\id$, with $L$ an infinitesimal character of~$\Rr(G)$. Here
the convolution of an endomorphism of~$\Rr(G)$ and an element
of~$O(G)$ is clearly well defined; and indeed
$$
L * \id(h_1h_2) = D(h_1)h_2 + h_1D(h_2),
$$
after a short calculation. Right invariance of~$L*\id$ is clear.
Reciprocally $\eta D$ is an infinitesimal character. All this is
in~\cite{TheSecondComing}.

The books~\cite{DieudonneIV,Warner,Taylor} and the review
article~\cite{Cartier2006} are good references for most of this
section.

\section{Some structure results for Hopf algebras}

Familiarity with the tensor $\Te(V)$ and cotensor (or shuffle)
$\Te^*(V)$ Hopf algebras is very convenient; we survey them here.
Consider a countable basis $B=\set{v_1,\dots,v_p,\dots}$ of the vector
space~$V$, and think of it as an alphabet, a word of this alphabet
being a finite sequence of $v$'s. We let $\Te^*(V)$ be the vector space
generated by the set of words and~$1$ (corresponding to the empty
word). The length of a word $w=v_{i_1}\cdots v_{i_n}$ is denoted by
$|w|=n$; naturally $|1|=0$. Introduce a noncocommutative
(deconcatenation) coproduct on~$\Te^*(V)$ by the formulae $\Dl1=1\ox1$
and
$$
\Dl w = \sum_{p=0}^{n} v_{i_1}\cdots v_{i_p} \ox v_{i_{p+1}}\cdots
v_{i_n},
$$
with the agreement that when all the terms are on the one side of
the tensor sign there is a~$1$ on the other side. Notice that
\begin{equation*}
(\Dl \ox \id)\Dl(v_{i_1}\cdots v_{i_n}) = \sum_{0\le p\le q\le n}
v_{i_1} \cdots v_{i_p} \ox v_{i_{p+1}} \cdots v_{i_q} \ox v_{i_{q+1}}
\cdots v_{i_n} = (\id \ox \Dl)\Dl (v_{i_1} \cdots v_{i_n}),
\end{equation*}
understanding that when $p=q$ the middle term of the summand is~$1$.
Hence $\bigl(\Te^*(V),\Dl,\eta\bigr)$, where $\eta\: \Te^*(V)\to\R$ is
defined by $\eta(1)=1$ and $\eta(w)=0$ if $|w|>0$, is indeed a
coalgebra; by the way, any commutative $\Q$-algebra at the place of
the real numbers would do here. The dual vector space consists of all
infinite series of the form $\sum \la_I v'_I$, where
$I=\set{i_1,\dots,i_n}$ and~$v'_I$ denotes the dual of $v_I:=
v_{i_1}\cdots v_{i_n}$. It becomes an algebra with product
$$
\<v'_J v'_K, v_I> := \<v'_J\ox v'_K, \Dl v_I>.
$$
Since the right hand side vanishes unless $J\cup K=I$ as ordered
sets, and in that case equals~$1$, this product is simply
concatenation:
$$
v'_J v'_K = v'_{j_1}\cdots v'_{j_m} v'_{k_1}\cdots v'_{k_l}.
$$
In other words, this dual is the algebra $\R[[B']]$ of noncommutative
formal power series in the variables $v_i$, which is the (Krull
topology) completion of the algebra $\R[B']$ of noncommutative
polynomials in the same variables ---that is the tensor algebra
$\Te(V)$, as tensor product is given by concatenation.

It is clear that $\Te(V)$ is a free associative algebra on
$B'$~\cite{Reutenauer}. Moreover, if $\L(B')$ is the free Lie algebra
on $B'$, from the universal properties of~$\L(B')$ and
of~$\U\bigl(\L(B')\bigr)$ it follows that also $\U\bigl(\L(B')\bigr)$
is a free associative algebra on~$B'$; therefore $\Te(V)=
\U\bigl(\L(B')\bigr)$. In particular, we have a Hopf algebra structure
on~$\Te(V)$, which is inherited by its completion. Their cocommutative
coproduct is given on monomials by the formula
\begin{equation*}
\Dl(v'_{i_1}\cdots v'_{i_n}) = \sum_{p=0}^n\sum_{\sg\in S_{n,p}}
v'_{\sg(i_1)}\cdots v'_{\sg(i_p)} \ox v'_{\sg(i_{p+1})}\cdots
v'_{\sg(i_n)}.
\end{equation*}
Here we deal with $(p,n-p)$-shuffles, that is, permutations~$\sg$ of
$[n]=\set{1,\dots,n}$ such that $\sg(1)<\cdots<\sg(p)$ and
$\sg(p+1)<\cdots<\sg(n)$; we write $\sg\in S_{n,p}$. From this
coproduct we obtain a product~$\pitchfork$ on $\Te^*(V)$ by dualization:
\begin{equation}
\<v'_I, v_J \pitchfork v_K> := \<\Dl v'_I, v_J \ox v_K>.
\label{eq:shuffleProduct}
\end{equation}
Explicitly, the commutative \textit{shuffle} product $\pitchfork$ is
given by
$$
v_{i_1}\cdots v_{i_p} \pitchfork v_{i_{p+1}}\cdots v_{i_n} =
\sum_{\sg\in S_{n,p}} v_{i_{\sg(1)}} \cdots v_{i_{\sg(n)}}.
$$
For instance $v_iv_j\pitchfork v_k=v_iv_jv_k+v_iv_kv_j+v_kv_iv_j$ and
\begin{equation*}
v_i v_j \pitchfork v_k v_l = v_iv_jv_kv_l + v_iv_kv_jv_l +
v_iv_kv_lv_j + v_kv_iv_jv_l + v_kv_iv_lv_j + v_kv_lv_iv_j.
\end{equation*}
It is also easy to check the following formula, which can be employed
as a recursive definition of the shuffle product:
\begin{equation*}
v_{i_1}\cdots v_{i_p} \pitchfork v_{i_{p+1}}\cdots v_{i_n} =
(v_{i_1}\cdots v_{i_p} \pitchfork v_{i_{p+1}}\cdots
v_{i_{n-1}})v_{i_n} + (v_{i_{p+1}}\cdots v_{i_n} \pitchfork
v_{i_1}\cdots v_{i_{p-1}})v_{i_p}.
\end{equation*}
The coproduct on $\Te^*(V)$ and $\pitchfork$ are compatible since they
are respectively obtained by dualization of the product and coproduct
of the Hopf algebra $\R[[B']]$. The resulting commutative, connected,
graded Hopf algebra $\Sh(V)\equiv\Te^*(V)$ is called the shuffle Hopf
algebra over~$V$. The construction does not depend on the choice of
the basis $B$, since all the algebras involved only depend on the
cardinality of $B$. The antipode on~$\Te(V)$ is given by
$$
S(v'_1\cdots v'_n)= (-1)^n v'_n\cdots v'_1;
$$
by duality the same formula holds on~$\Sh(V)$.

Every polynomial $P\in \R[B']$ can be written in the form $P=\sum_n
P_n$ where $P_n$ is the sum of all monomials of $P$ of degree $n$. It
is called a \textit{Lie element} if $P_0=0$ and each $P_n$ belongs to
the free Lie algebra generated by the~$v'_i$. The following is a
classical theorem by Friedrichs.

\begin{thm}
A polynomial $P$ is a Lie element if, and only if, it is primitive.
\end{thm}

\begin{proof}
The `only if' part follows easily by induction from the obviously true
assertion for $n=1$. To prove the converse we invoke in context the
Dynkin operator~$D$, for whose study we
recommend~\cite{Eunomia,WignerGood,FredericReut}. An abstract
definition of~$D$ is $D=S*Y$, where $Y$ is the derivation given by the
grading; equivalently $\id*D = Y$. If~$P_n$ is primitive, then so are
$$
nP_n = YP_n = \pi(\id \otimes D)(1 \otimes P_n + P_n \otimes 1) =
D(P_n);
$$
and vice versa. But $D(P_n)$ is a Lie element, as it corresponds to
the left-to-right bracketing:
$$
D(x_{i_1}\ldots x_{i_n}) = [\dots[[x_{i_1}, x_{i_2}], x_{i_3}], \dots,
x_{i_n}] \sepword{(the Dynkin--Specht--Wever theorem).}
$$
To prove the last equality, note that it is true for~$|w|=1$. Assume
that it holds for all words of degree less than~$n$, and let
$w=xx_{i_n}=x_{i_1} \cdots x_{i_n}$. Then
$$
\Dl w = \Dl x \Dl x_{i_n} = (x_{(1)} \ox x_{(2)})(x_{i_n} \ox 1 + 1\ox
x_{i_n}) = x_{(1)}x_{i_n} \ox x_{(2)} + x_{(1)} \ox x_{(2)}x_{i_n}.
$$
Since $Sx_{i_n} = -x_{i_n}$ and $\eta(x)=0$,
\begin{align*}
(S*Y)w &= S(x_{(1)}x_{i_n})Yx_{(2)} + Sx_{(1)}Y(x_{(2)}x_{i_n})
\\
&= S(x_{i_n})Sx_{(1)}Yx_{(2)} + Sx_{(1)}Yx_{(2)}x_{i_n} +
Sx_{(1)}x_{(2)}Yx_{i_n}
\\
&= -x_{i_n}Sx_{(1)}Yx_{(2)} + Sx_{(1)}Yx_{(2)}x_{i_n}
\\
&= [Dx, x_{i_n}] = [\dots[[x_{i_1}, x_{i_2}], x_{i_3}], \dots,
x_{i_n}],
\end{align*}
upon using the induction hypothesis in the last equality. (The
definition $D=Y*S$ would work the same, yielding right-to-left
bracketing.)
\end{proof}

When $V$ is finite dimensional, the previous argument of Friedrichs'
theorem goes through for formal power series, because the homogeneous
components are polynomials, and there is only a finite number of words
of a given length. In the infinite-dimensional case Lie series are
defined as those such that their projections to any finite-dimensional
subspace~$\tilde V$ are Lie series over~$\tilde V$, so the theorem
also holds for series in the infinite-dimensional
context~\cite[Section~3.1]{Reutenauer}.

Given a power series $Z$, let us denote by $(Z,w)$ the coefficient of
the word $w$ in $Z$. The topology in $\R[[B']]$ alluded above is the
weakest topology such that for each~$w$ the mapping $Z\mapsto (Z,w)$
is continuous, when $\R$ is equipped with the discrete topology. In
particular, the neighbourhoods of 0~are indexed by finite sets of
words, and correspond to those series whose coefficients vanish on all
the words of the given finite set. Thus, given a sequence of series
$(Z_n)$ such that for each neighbourhood of~0, all but a finite number
of~$Z_n$'s are in this neighbourhood, their sum $\sum_n Z_n$ is
defined as the power series~$Z$ satisfying
$$
(Z,w) = \sum_n (Z_n,w).
$$
This sum makes sense since only finitely many terms are different
>from zero for each $w$. Notice that $Z$ can be written as
$Z=\sum_w (Z,w)w$, where the sum runs over the set of words.
(Henceforth we shall no longer be fussy on `topological' matters.)

When $Z$ is a series such that $(Z,1)=0$, then the expression $\sum_n
\la_n Z^n$ has a meaning for any choice of the numbers $\la_n$. In
particular, we may define exponentials and logarithms as usual:
$$
\exp(Z) = \sum_{n=0}^\infty \frac{Z^n}{n!} \sepword{and} \log(1+Z) =
\sum_{n=1}^\infty \frac{(-1)^{n-1}}{n} Z^n,
$$
As expected
$$
\log\bigl(\exp(Z)\bigr) = Z, \sepword{and} \exp\bigl(\log(1 + Z)\bigr)
= 1 + Z.
$$
and routine calculations establish that $\exp$ is a bijection from
the set of primitive elements in the completion of~$\Te(V)$ into the
set of grouplike elements, and vice versa for~$\log$.

Equation~\eqref{eq:shuffleProduct} entails
$$
\Dl Z = \sum_{w,x} (Z,w\pitchfork x) w\ox x.
$$
Since for grouplike elements $\Dl Z = \sum_{w,x}(Z,w)(Z,x)w\ox x$,
it follows that
\begin{equation}
(Z,w\pitchfork x) = (Z,w)(Z,x),
\label{eq:Ree}
\end{equation}
for them. This of course means that the grouplike elements
of~$\R[[B']]$ are precisely those~$Z$ for which the map
$w\mapsto(Z,w)$ is an algebra homomorphism for the shuffle product.
This characterization is originally due to Ree~\cite{Ree}.

\smallskip

We collect next some elements of structure theory of commutative or
cocommutative Hopf algebras ---mostly due to
Patras~\cite{Frederic1,Frederic2}--- beginning by a `double series'
argument similar to the one in~\cite{Reutenauer} for the
shuffle-deconcatenation Hopf algebra.

Consider, for $H=\bigoplus_m^\infty H^{(m)}$ a graded connected
\textit{commutative} Hopf algebra with augmentation ideal~$H_+$ and
graded dual~$H'$, a suitable completion~$H\barox H'$ of the tensor
product~$H\ox H'$. This is a unital algebra, with product $m\ox\Dl^t$
and unit~$1\ox1$. Now by Leray's theorem ---an easy dual version of
the Cartier--Milnor--Moore theorem--- our $H$ is a symmetric algebra
over a supplement~$V$ of~$H_+^2$ in~$H_+$~\cite{Quaoar,Frederic2}. Let
$A$ index a basis for~$V$, let ${\tilde A}$ (the monoid freely
generated by~$A$) index the words~$X_u$, and let~$Z_u$ denote an
element of the dual basis in~$H'$; then the product on~$H \barox H'$
is given by the double series product:
$$
\biggl(\,\sum_{u,v\in{\tilde A}}\a_{uv}X_u\ox Z_v\biggr)
\biggl(\,\sum_{w,t\in{\tilde A}}\b_{wt}X_w\ox Z_t\biggr)
:= \sum_{u,v,w,t\in{\tilde A}} \a_{uv} \b_{wt} \,X_u X_w\ox Z_v Z_t.
$$
The linear embedding $\End H \to H \barox H'$ given by
$$
f \mapsto \sum_{u\in{\tilde A}}f(X_u)\ox Z_u,
$$
is really a convolution algebra embedding
$$
(\End H,*) \to (H \barox H',m\ox\Dl^t).
$$
Indeed,
\begin{align}
\biggl(\,\sum_{u\in{\tilde A}} f(X_u) \ox Z_u\biggr)
&\biggl(\,\sum_{v\in{\tilde A}} g(X_v) \ox Z_v\biggr)
= \sum_{u,v\in{\tilde A}} f(X_u) g(X_v) \ox Z_u Z_v
\notag \\
&= \sum_{t\in{\tilde A}} \biggl(\,\sum_{u,v\in{\tilde A}}
f(X_u) g(X_v) \,\<Z_uZ_v,X_t> \biggr) \ox Z_t
\notag \\
&= \sum_{t\in{\tilde A}} \biggl(\,\sum_{u,v\in{\tilde A}}
f(X_u) g(X_v) \,\<Z_u\ox Z_v,\Dl X_t> \biggr) \ox Z_t
\notag \\
&= \sum_{t\in{\tilde A}} f*g(X_t) \ox Z_t.
\label{eq:conv-series}
\end{align}
Notice that the identities $u\eta$ for convolution and~$\id$ for
composition in~$\End H$ correspond respectively to
$$
u\eta \mapsto 1\ox1 \sepword{and}
\id \mapsto \sum_{u\in {\tilde A}}X_u\ox Z_u.
$$
Denote
$$
\pi_1(X_w) := \sum_{k\ge1} \frac{(-1)^{k-1}}{k} \sum_{u_1, \dots,
u_k\ne1} \<Z_{u_1} \cdots Z_{u_k},X_w>\, X_{u_1} \cdots X_{u_k} =:
\log^*\id\,X_w.
$$
Using the same idea as in~\eqref{eq:conv-series}, we get
\begin{align}
\log\biggl(\,\sum_{u\in{\tilde A}} X_u\ox Z_u \biggr)
&:= \sum_{k\ge1} \frac{(-1)^{k-1}}{k}
\biggl(\,\sum_{u\ne1} X_u \ox Z_u \biggr)^k
\nn \\
&= \sum_{k\ge1}\frac{(-1)^{k-1}}{k}
\sum_{u_1,\dots,u_k\ne1}
X_{u_1} \cdots X_{u_k} \ox Z_{u_1} \cdots Z_{u_k}
\nn \\
&= \sum_{w\in{\tilde A}}\sum_{k\ge1}\frac{(-1)^{k-1}}{k}
\sum_{u_1,\dots,u_k\ne1}
\<Z_{u_1} \cdots Z_{u_k},X_w>\, X_{u_1} \cdots X_{u_k} \ox Z_w
\nn \\
&= \sum_{w\in{\tilde A}} \pi_1(X_w) \ox Z_w.
\label{eq:basic-truth}
\end{align}
We moreover consider the endomorphisms $\pi_n:=\pi_1^{*n}/n!$ so that,
by~\eqref{eq:conv-series}:
$$
\sum_{w\in\tilde A}\pi_n(X_w)\ox Z_w =
\frac{1}{n!}\biggl(\,\sum_{v\in\tilde A}\pi_1(X_v)\ox Z_v\biggr)^n.
$$
We may put $\pi_0 := u\eta$. Thus, if $a\in H$ is of order~$n$,
$\pi_m(a) = 0$ for $m > n$. Furthermore, for~$n > 0$,
\begin{equation}
\id^{*l}a = \exp^*(\log^*(\id^{*l}))\,a =
\sum_{m=1}^n\frac{(\log^*(\id^{*l}))^m}{m!}a =
\sum_{m=1}^nl^m\frac{(\log^*\id)^m}{m!}a = \sum_{m=1}^nl^m\pi_m(a).
\label{eq:greater-good}
\end{equation}
In particular $\id=\sum_{m\ge0}\pi_m$. The graded maps $\id^{*n}$ are
called the \textit{Adams operations} or characteristic endomorphisms
of~$H$; they play an important role in the (Hochschild, cyclic)
cohomology of commutative algebras~\cite{GS1,NoOffenseIntended,GS2}.
The $\pi_n$ are often called \textit{Eulerian idempotents}. We have
for them:

\begin{prop}
\label{pr:Trick}
For any integers $n$ and $k$,
\begin{equation}
\id^{*n}\,\id^{*k} = \id^{*nk} = \id^{*k}\,\id^{*n}.
\label{eq:burden of proof}
\end{equation}
and
\begin{equation}
\pi_m\pi_k = \dl_{mk}\,\pi_k.
\label{eq:make-believe}
\end{equation}
\end{prop}

\begin{proof}
The first assertion is certainly true for $k=1$ and all integers~$n$,
and if it is true for some~$k$ and all integers~$n$, then taking into
account that $\id$ is an algebra homomorphism, the induction
hypothesis gives
$$
\id^{*n}\,\id^{*k+1}=\id^{*n}(\id^{*k}*\id)
=\id^{*nk}*\id^{*n}=\id^{*n(k+1)}.
$$
Substituting the final expression of~\eqref{eq:greater-good}
in~\eqref{eq:burden of proof}, with very little work one
obtains~\eqref{eq:make-believe}. So indeed the $\pi_k$ form a family
of orthogonal projectors.
\end{proof}

Thus the space $H=\bigoplus_m^\infty H^{(m)}$ always has the
direct sum decomposition
\begin{equation}
H=\bigoplus_{n\ge0} H_n := \bigoplus_{n\ge0}\pi_n(H).
\label{eq:decomposition}
\end{equation}

Moreover, from~\eqref{eq:greater-good},
$$
\id^{*l} H_n = l^n H_n,
$$
so the $H_n$ are the common eigenspaces of the operators
$\id^{*l}$ with eigenvalues $l^n$. Thus, the
decomposition~\eqref{eq:decomposition} turns $H$ into a graded
algebra. Indeed, if $a\in H_r$ and $b\in H_s$, then
$$
\id^{*l}(ab) = \id^{*l}a\,\id^{*l}b = l^{r+s}(ab),
$$
and therefore $m$ sends $H_r\ox H_s$ into~$H_{r+s}$.
We shall denote by $\pi_n^{(m)}$ the restriction of $\pi_n$ to
$H^{(m)}$, the set of elements of degree $m$, with respect to the
original grading.

\smallskip

If $H$ is cocommutative instead of commutative, the previous arguments
go through. One then has
\begin{equation*}
\log\biggl(\,\sum_{u\in{\tilde A}}X_u \ox Z_u\biggr) =
\sum_{w\in{\tilde A}}X_w \ox \pi_1(Z_w).
\end{equation*}
Furthermore, in this case the Eulerian idempotents of~$H$ are the
transpose of the Eulerian idempotents of the graded commutative
Hopf algebra~$H'$. In particular, for $H$ cocommutative,
$\pi_1(H)=P(H)$, the Lie algebra of primitive elements in~$H$.
This is easily sharpened into the following
version~\cite{Frederic2} of the Cartier--Milnor--Moore theorem:
the inclusion $\pi_1(H)\hookrightarrow H$ extends to an
isomorphism of~$\U\bigl(\pi_1(H)\bigr)$ with~$H$.

\section{The CBHD \textit{development} and Hopf algebra}

There are three paradigmatic methods (and sundry hybrid forms) to deal
with first order non-autonomous differential equations: the iteration
formula or Dyson--Chen expansional, the Magnus expansion and the
product integral. For reasons expounded later, at the beginning of
Section 10, in this paper we look first for the Magnus
expansion~\cite{Magnus}. In the influential
paper~\cite{inflationbuster} dealing with the latter method (although
Magnus' seminal contribution is not mentioned) the famous
Campbell--Baker--Hausdorff--Dynkin (CBHD) formula in Lie algebra
theory is shown to be a special case of general formulas for the
solution of~\eqref{eq:menos-gatos}. This is scarcely surprising, as
that solution involves some kind of exponential with non-commuting
exponents; also the quest for `continuous analogues' of the CBHD
formula was a motivation for Chen's work. Conversely, a heuristic
argument for obtaining Magnus' expansion from the CBHD formula has
been known for some time~\cite{macrote,MexicoLindo}; and a routine, if
rigorous and Hopf flavoured as well, derivation of Magnus' method from
CBHD is available in~\cite{MexicoOtraVez}. Hence the interest, as a
prelude to our own derivation of the Magnus expansion from the CBHD
development (that will employ the concept of nonlinear CBHD
\textit{recursion} and Rota--Baxter theory techniques) of rendering
the proof of the CBHD expansion in Hopf algebraic terms. This was
recognized as the deeper and more natural approach to the subject some
fifteen years ago, but remains to date woefully ignored. Standard
treatments of the CBHD development can be found in good Lie group
theory books like~\cite{H-ChProphet}.

In the sequel we follow~\cite{RiauRiau} and~\cite{LodayEM}. It will be
soon clear to the reader, according to the previous discussion, that
the CBHD formulae are universal; thus we can as well return to the
case where $H$ is the Hopf tensor algebra~$\Te(V)$ and where $V$
possesses a basis $B=\set{X_1,\ldots,X_n}$. The CBHD~series
$\sum_{m\ge1} \Phi_m(X_1,\dots,X_n)$ is defined by
$$
\sum_{m\ge1}\Phi_m(X_1,\ldots,X_n)
= \log\bigl(e^{X_1}\cdots e^{X_n}\bigr),
$$
where $\Phi_m(X_1,\ldots,X_n)$ are homogeneous polynomials of
degree~$m$.

Now, if $a$ is a grouplike element in a Hopf algebra $H$, and
$f,h\in\Hom(H,A)$, where $A$ is a unital algebra, then
$$
f*h(a) = f(a)h(a).
$$
In particular
\begin{equation}
\log\bigl(e^{X_1}e^{X_2}\cdots e^{X_n}\bigr) =
\log^*\id\bigl(e^{X_1}e^{X_2}\cdots e^{X_n}\bigr) =:
\pi_1\bigl(e^{X_1}e^{X_2}\cdots e^{X_n}\bigr).
\label{eq:tabla-de-salvacion}
\end{equation}
Take first $n=2$. Then
$\Phi_m(X,Y)=\pi_1^{(m)}\bigl(e^Xe^Y\bigr)$. The Cauchy product
gives
$$
e^Xe^Y= \sum_{m\ge0}\left(\sum_{i=0}^m
\frac{X^i}{i!}\,\frac{Y^{n-i}}{(n-i)!}\right),
$$
hence
$$
\Phi_m(X,Y)=\sum_{i+j=m}\frac{1}{i!j!}\;\pi_1^{(m)}(X^iY^j).
$$
A similar argument entails the following proposition.

\begin{prop}
\label{pr:Dirty}
\begin{equation}
\Phi_m(X_1,\ldots,X_n) = \sum\frac{1}{i_1!\cdots i_n!}\,
\pi_1^{(m)}(X_1^{i_1}X_2^{i_2}\cdots X_n^{i_n}),
\label{eq:at-sea}
\end{equation}
where the sum runs over all vectors $(i_1,\dots,i_n)$, with
nonnegative coordinates, such that $i_1+\cdots+i_n=m$.
\end{prop}

Denote by $\vf_n(X_1,\ldots,X_n)$ the `multilinear' part of
$\Phi_n(X_1,\ldots,X_n)$ ---that is, the homogeneous polynomial of
degree $n$ that consist of those monomials of $\Phi_n(X_1,\ldots,X_n)$
that include all the $X_i$'s. This amounts to take $X_i^2=0$ in
$\Phi_n(X_1,\ldots,X_n)$, for all $i$. So
by~\eqref{eq:tabla-de-salvacion}
$$
\vf_n(X_1,\ldots,X_n) = \pi_1^{(n)}(X_1\cdots X_n),
$$
since in that case
\begin{align}
e^{X_1}\cdots e^{X_n} &= (1 + X_1) \cdots (1 + X_n)
\nonumber \\
&= \sum_iX_i + \sum_{i<j} X_iX_j + \sum_{i<j<k}X_iX_jX_k + \cdots +
X_1\cdots X_n.
\label{eq:useful}
\end{align}

Now, if $X_\sg:=(X_1,\dots,X_n)\cdot\sg:=
(X_{\sg(1)},\dots,X_{\sg(n)})$ denotes the standard right action of
the symmetric group~$S_n$ on~$V^{\ox n}$, then the monomials that
include all the $X_i$'s are of the form $X_\sg$, therefore
$$
\pi_1^{(n)}(X_1\cdots X_n) = \sum_{\sg\in S_n} c_\sg X_\sg,
$$
for some coefficients $c_\sg$, that we shall determine in a
moment.

\begin{prop}
\label{pr:NotQuick}
\begin{equation}
\pi_1^{(n)}(X_1\cdots X_n) = \sum_{\sg\in S_n}
\frac{(-1)^{d(\sg)}}{n}{\binom{n-1}{d(\sg)}}^{-1}X_\sg.
\label{eq:la-madre-del-cordero}
\end{equation}
where $d(\sigma)$ is the number of {\rm descents} of~$\sigma$,
that is, the number of `errors' in ordering consecutive terms in
$\sg(1),\dots,\sg(n)$.
\end{prop}

\begin{proof}
Assume that $\sg$ has $d$ descents, say in $n_0,n_0+n_1,n_0+n_1+\cdots
+n_{j-1}$, set $n_j=n-n_0-\cdots-n_{j-1}$ and let $Z=\sum_i X_i +
\sum_{i<j} X_iX_j+\cdots+X_1\cdots X_n$. By~\eqref{eq:useful},
$e^{X_1}\cdots e^{X_n}=Z+Y$, where $Y$ is a collection of terms that
contains at least one factor of the form $X^2_i$, therefore they will
not contribute to the coefficient of $X_\sg$, and we neglect them.
Now, since $\log(1+Z)=\sum\frac{(-1)^j}{j}Z^j$ we have to compute the
contribution $c(j)$ from each power $Z^j$.

Suppose that the monomial $X_{\sg(1)}\cdots X_{\sg(n_0)}$ is built
>from $j_1$ monomials of $Z$,  and in general that each monomial
$X_{\sg(n_0+\cdots+n_{i-1}+1)}\cdots X_{\sg(n_0+\cdots+n_i)}$ is
the product of $j_i$ monomials of $Z$. Notice that there are
${n_i-1\choose j_i-1}$ manners to construct each monomial, in such
a way, because $X_{\sg(n_0+\cdots+n_{i-1}+1)}$ is always in the
first monomial, and once the first $j_i-1$ monomials are chosen,
the last monomial is fixed since $\sg$ is increasing in each
segment. Thus
$$
c(j)=\sum_{(j_0,\dots,j_d)}{n_0-1\choose j_0-1}
{n_1-1\choose j_1-1}\cdots{n_d-1\choose j_d-1}\;,
$$
where the sum extends over all vectors $(j_0,\dots,j_d)$ satisfying
$j_0+\cdots+j_d=j$. Since ${n_k-1\choose j_k-1}$ is the coefficient of
$x^{j_k-1}$ in the binomial expansion of $(1+x)^{n_k-1}$, and
$\sum_{i=0}^d(j_i-1)=j-d-1$, $c(j)$ is the coefficient of $x^{j-d-1}$
in
$$
\prod_{i=0}^d(1+x)^{n_k-1}=(1+x)^{\sum_{i=0}^d(n_i-1)}
=(1+x)^{n-d-1},
$$
we therefore conclude that
$$
c(j)={n-d-1\choose j-d-1}\;.
$$
Now, we have $j\le n$ since $X_\sg$ has $n$ letters. Also $j\ge d+1$
as $X_\sg$ is broken in $d+1$ parts. Therefore
$$
c_\sg=\sum_{j=d+1}^n\frac{(-1)^{j-1}}{j}{n-d-1\choose j-d-1}
=(-1)^d\sum_{i=0}^m\frac{(-1)^i}{i+d+1}{m\choose i},
$$
where $m=n-d-1$. Now, from the binomial
identity
$$
\int_0^1\!(1-x)^mx^d\,dx
=\sum_{i=0}^m(-1)^i{m\choose i}\int_0^1\!x^{i+d}\,dx
=\sum_{i=0}^m\frac{(-1)^i}{i+d+1}{m\choose i}.
$$
Finally, a simple induction, using integration by parts, gives
$$
\int_0^1\!(1-x)^mx^{d}\,dx=\frac{d!\,m!}{(m+d+1)!}
=\frac{1}{m+d+1}{d+m\choose d}^{-1}
=\frac{1}{n}{n-1\choose d}^{-1}.
$$
Our task is over. But the number of descents will reappear soon
enough.
\end{proof}

This construction performed here is arguably more elegant and simpler
than the standard treatments of the CBHD development by purely Lie
algebraic methods. We came in by the backdoor, using the bigger free
associative algebra, knowing that $\log\bigl(e^{X_1}\cdots
e^{X_n}\bigr)$ ---and each of its homogeneous parts--- is primitive,
i.e., a Lie element; and that we have the Dynkin operator to rewrite it
in terms of commutators.

Let us exemplify with the case $n=2$. Obviously we have
$$
\Phi_1(X,Y) = X + Y; \qquad \Phi_2(X,Y) = \frac{1}{2}[X,Y].
$$
Now,
$$
\pi_1^{(3)}(X_1X_2X_3)=\frac{1}{3}X_{(123)}
-\frac{1}{6}\left(X_{(132)}+X_{(213)}+X_{(231)}+X_{(312)}\right)
+\frac{1}{3}X_{(321)}.
$$
Therefore
\begin{align*}
\Phi_3(X,Y) &=\frac{1}{2}(\pi_1^{(3)}(X^2Y)+\pi_1^{(3)}(XY^2)) \\
&=\frac{1}{2}\left(\frac{1}{6}X^2Y-\frac{1}{3}XYX+\frac{1}{6}YX^2
+\frac{1}{6}XY^2-\frac{1}{3}YXY+\frac{1}{6}Y^2X\right) \\
&=\frac{1}{12}\left([[X,Y],Y]-[[X,Y],X]\right).
\end{align*}
Both cubic Lie elements appear in~$\Phi_3$. Similarly
\begin{align*}
\pi_1^{(4)}(X_1X_2X_3X_4)&=\frac{1}{4}X_{(1234)}
-\frac{1}{12}\biggl(X_{(1243)}+X_{(1324)}+X_{(1342)}
+X_{(1423)}+X_{(2134)} \\
&\quad+X_{(2314)}+X_{(2341)}+X_{(2413)}+X_{(3124)}+X_{(3412)}
+X_{(4123)}\biggr) \\
&\quad+\frac{1}{12}\biggl(X_{(1432)}+X_{(2143)}+X_{(2431)}+X_{(3142)}
+X_{(3214)}+X_{(3241)}  \\
&\quad+X_{(3421)}+X_{(4132)}+X_{(4213)}+X_{(4231)}+X_{(4312)}\biggr)
-\frac{1}{4}X_{(4321)}
\end{align*}
We concentrate on $X^2Y^2$, as it is clear that most terms coming
from~$X^3Y$ or~$XY^3$ will vanish; and in fact the corresponding
contributions \textit{in toto} come to naught. We obtain
\begin{align*}
\Phi_4(X,Y) &= -\frac{1}{192}\bigl(4XYXY + 2XY^2X +2YX^2Y - 3YXYX
- 2XY^2X - 3YX^2Y\bigr)  \\
&=-\frac{1}{24}[[[X,Y],X],Y].
\end{align*}
The identity of Jacobi has been used, under the form
$$
[[[X,Y],X],Y] = [[[X,Y],Y],X].
$$
It is remarkable that the other quartic Lie elements, $[[[X,Y],X],X]$
and $[[[X,Y],Y],Y]$, do not appear in the fourth degree term.

\section{Rota--Baxter maps and the algebraization of integration}

This paper draws inspiration partly from~\cite{TheSecondComing}, where
Connes and Marcolli have introduced logarithmic derivatives in the
context of Hopf algebras. Our intent and methods are different; but it
is expedient to dwell here a bit on their considerations. Given~$H$
and~$A$ commutative as in the last part of Section~4, and a
derivation~$\dl$ on~$A$, for a multiplicative map $\phi \in G_H(A)$
Connes and Marcolli define two maps in $\Hom(H,A)$ by $\dl(\phi):=\dl
\circ \phi$, and then
$$
D_\dl(\phi) := \phi^{-1} * \dl(\phi)
$$
This yields an $A$-valued infinitesimal character. Indeed, using
Sweedler's notation and multiplicativity of~$\phi\in G_H(A)$, one
has
\begin{eqnarray*}
D_\dl(\phi)[cd] &=& \phi^{-1} * \dl(\phi)[cd] = m_A(\phi^{-1} \ox
\dl(\phi))\Delta(cd) =
\phi^{-1}(c_{(1)}d_{(1)})\dl(\phi(c_{(2)}d_{(2)}))
\\
&=&
\phi^{-1}(c_{(1)})\phi^{-1}(d_{(1)})\bigl(\dl(\phi(c_{(2)}))\phi(d_{(2)})
+ \phi(c_{(2)})\dl(\phi(d_{(2)}))\bigr)
\\
&=&\phi^{-1}(c_{(1)})\dl(\phi(c_{(2)}))\phi^{-1}(d_{(1)})\phi(d_{(2)}) +
\phi^{-1}(c_{(1)})\phi(c_{(2)})\phi^{-1}(d_{(1)})\dl(\phi(d_{(2)}))
\\
&=& D_\dl(\phi)[c] \eta_A(d) + \eta_A(c) D_\dl(\phi)[d].
\end{eqnarray*}
Therefore $D_\dl(\phi)$ belongs to~$\g_H(A)$.

The Dynkin operator appearing in Section~5 ---one of the fundamental
Lie idempotents in the theory of free Lie
algebras~\cite{Reutenauer,FredericReut}--- is a close cousin of the
logarithmic derivative $D_\dl(g)$. Consider~$G_H(H)$, for $H$
connected and graded. The grading operator~$Y$ is a derivation of~$H$
$$
Y(hh') = Y(h)h'+ hY(h') =: |h|hh' + hh'|h'|.
$$
The map~$Y$ extends naturally to a derivation on $\End(H)$. With~$f,g
\in\End(H)$ and $h\in H$ we find
\allowdisplaybreaks{
\begin{eqnarray*}
Y(f \ast g)(h) &:=& f * g\,(Y(h)) = |h|(f * g)(h) =
|h|f(h^{(1)})g(h^{(2)})
\\
&=& |h^{(1)}|f(h^{(1)})g(h^{(2)}) + |h^{(2)}|f(h^{(1)})g(h^{(2)})
= Yf* g\,(h) + f * Yg\,(h),
\end{eqnarray*}}
where we used that $\Delta (Y(h))=|h|\Delta (h)=\bigl(|h^{(1)}|+
|h^{(2)}|\bigr)\,h^{(1)}\ox h^{(2)}$. Now, as before, convolution
of the antipode~$S$ with the derivation~$Y$ of~$H$ defines a
Dynkin operator, to be interpreted as an $H$-valued infinitesimal
character~\cite{Eunomia}.

Suppose we have a smooth map $t \mapsto L(t)$ from~$\R_t$
to~$\g_H(A)$. We could say that one of the main aims of this paper
is to solve for~$g(t)$ the initial value scheme
\begin{equation}
D_{d/dt}\bigl(g(t)\bigr) = L(t); \qquad g(0) = \eta_A,
\label{eq:donya-toda}
\end{equation}
at least for (real and) complex points. Now, both the classical
notions of derivation and integration have interesting
generalizations. It would then be a pity to limit ourselves to the
classical framework; and so we now jump onto a somewhat more
adventurous path.

For integration, one lacks a good algebraic theory similar to the one
developed in~\cite{Kolchin}, say. Next we elaborate on a somewhat
unconventional presentation of the integration-by-parts rule using the
algebraic notion of the weight-$\theta$ Rota--Baxter relation
corresponding to the generalization of the Leibniz rule in terms of
weight-$\theta$ \textit{skewderivations}. One should strive for
nothing less ambitious than developing Rota's program, beautifully
outlined in~\cite{Rota98} in the context of Chen's
work~\cite{Chen}, of establishing an algebraic theory of integration
in terms of generalizations of the integration-by-parts rule.

Let us recall first the integration-by-parts rule for the Riemann
integral map. Let $A:=C(\R)$ be the ring of real continuous functions.
The indefinite Riemann integral can be seen as a linear map on~$A$
\begin{equation}
I: A \to A, \qquad  I(f)(x) := \int_0^x f(t)\,dt.
\label{eq:Riemann}
\end{equation}
Then, integration-by-parts for the Riemann integral can be written
as follows. Let
$$
F(x) := I(f)(x) = \int_0^x f(t)\,dt, \qquad G(x) := I(g)(x) = \int_0^x
g(t)\,dt;
$$
then
$$
\int_0^x F(t)\frac{d}{dt}\bigl(G(t)\bigr)\,dt = F(x)G(x) -
\int_0^x \frac{d}{dt}\bigl(F(t)\bigr)G(t)\,dt.
$$
More compactly, this well-known identity is written
\begin{equation}
I(f)(x)I(g)(x) = I\bigl( I(f) g \bigr)(x) + I\bigl( fI(g)
\bigr)(x), \label{eq:integ-by-parts}
\end{equation}
dually to the Leibniz rule.

Now, we introduce so-called skewderivations of weight~$\theta\in\R$ on
an algebra~$A$~\cite{JosephA}. A skewderivation is a linear map $\dl:A
\to A$ fulfilling the condition
\begin{equation}
\dl(ab) = a\dl(b) + \dl(a)b - \theta\dl(a)\dl(b).
\label{eq:crash}
\end{equation}
We call skewdifferential algebra a double $(A,\dl;\theta)$ consisting
of an algebra~$A$ and a skewderivation~$\dl$ of weight~$\theta$. A
skewderivation of weight~$\theta=0$ is just an ordinary derivation. An
induction argument shows that if~$A$ is commutative we have
$$
\dl(a^n) = \sum_{i=1}^{n}{n \choose i}(-\theta)^{i-1}a^{n-i} \dl(a)^{i}.
$$
Also
$$
\dl^{n}(ab) = \sum_{i=0}^{n}{n\choose
i}\sum_{j=0}^{n-i}{n-i\choose
j}(-\theta)^{i}\dl^{n-j}(a)\dl^{i+j}(b).
$$
Both formulae generalize well-known identities for an ordinary
derivation. We mention examples. First, on a suitable function
algebra~$A$ the simple finite difference operation $\dl: A\to A$
of step~$\la$,
\begin{equation}
\dl(f)(x) := \frac{f(x - \la)-f(x)}{\la},
\label{skewDiffexample1}
\end{equation}
satisfies identity~\eqref{eq:crash} with $\theta=-\lambda$.
See~\cite{Zudilin} for an interesting application of the $\la=1$ case
in the context of multiple zeta values. A closely related, though at
first sight different, example is provided by the $q$-difference
operator
\begin{equation}
\dl_qf(x) := \frac{f(qx)-f(x)}{(q-1)x}
\label{qDiff}
\end{equation}
which satisfies the $q$-analog of the Leibniz rule,
\begin{equation*}
\dl_q(fg)(x) = \dl_q f(x) g(x) + f(qx)\dl_q g(x) = \dl_q f(x)g(qx) +
f(x)\dl_q g(x).
\end{equation*}
This corresponds to relation~\eqref{eq:crash} for $\theta=(1-q)$,
modulo the identity
\begin{equation*}
\dl_q(fg)(x) = \dl_q f(x)g(x) + f(x)\dl_q g(x) + x(q-1)\dl_q f(x)
\dl_q g(x);
\end{equation*}
defining now $\bar\dl_q=x\dl_q$, it is a simple matter to check that
$\bar\dl_q$ is a skewderivation of weight~$1-q$.

We may ask for an \textit{integration operator} corresponding to the
skewderivation in~\eqref{skewDiffexample1}. On a suitable class of
functions, we define the summation operator
\begin{equation}
Z(f)(x) := \sum_{n\geq 1} \theta f(x + \theta n).
\label{eq:le-clou}
\end{equation}
{}For $\dl$ being the finite difference map of step~$\theta$,
\allowdisplaybreaks{
\begin{eqnarray*}
Z\dl(f)(x) &=& \sum_{n\geq 1} \theta\dl(f)(x + \theta n) = \sum_{n\geq
1}\theta\,\frac{f(x +\theta n - \theta) - f(x + \theta n)}{\theta}
\\
&=& \sum_{n\geq 1} f\bigl(x +\theta (n - 1)\bigr) - f(x + \theta
n) = \sum_{n\geq 0} f(x +\theta n) - \sum_{n\geq 1}f(x + \theta n)
= f(x).
\end{eqnarray*}}
As $\dl$ is linear we find as well $\dl Z(f)=f$. Observe, moreover,
that
\begin{align}
&\biggl(\sum_{n\geq 1} \theta f(x + \theta
n)\biggr)\biggl(\sum_{m\geq 1} \theta g(x + \theta m)\biggr) =
\sum_{n\geq 1, m\geq 1}\theta^2 f(x+\theta n)g(x+\theta m)
\nonumber \\
&= \biggl( \sum_{n > m \geq 1} + \sum_{m > n \geq 1} + \sum_{m = n
\geq 1}\biggr) \theta^2 f(x+\theta n) g(x+\theta m)
\nonumber \\
&= \sum_{m \geq 1} \biggl( \sum_{k \geq 1} \theta^2 f\bigl(x +
\theta(k+m)\bigr)\biggr) g(x +\theta m) + \sum_{n\geq 1}
\biggl(\sum_{k\geq 1} \theta^2 g\bigl(x + \theta
(k+n)\bigr)\biggr) f(x +\theta n)
\nonumber \\
&+ \sum_{n\geq 1} \theta^2 f(x+\theta n)g(x + \theta n) =
Z\bigl(Z(f)g\bigr)(x) + Z\bigl(fZ(g)\bigr)(x) + \theta Z(fg)(x).
\label{eq:parto-de-los-montes}
\end{align}

Related to the $q$-difference operator~\eqref{qDiff} there is the
\textit{Jackson integral}
\allowdisplaybreaks{
\begin{eqnarray*}
J[f](x) := \int_{0}^{x}f(y)\,d_qy = (1-q) \sum_{n \ge 0} f(xq^n) xq^n
\qquad (0<q<1).
\end{eqnarray*}}
This can be written in a more algebraic way, using the operator
$P_q[f] := \sum_{n>0}E_q^{n}[f]$, with the algebra endomorphism
($q$-dilatation) $E_q[f](x):=f(qx)$, for $f\in A$. The map~$P_q$
is a Rota--Baxter operator of weight~$-1$ and hence,
$\id+P_q=:\hat{P}_q$ is of weight~$+1$, see \cite{Rota2}.
Jackson's integral is given in terms of the above operators $P_q$
and the multiplication operator $M[f](x):= xf(x),\,f \in A$, by
$J[f](x) = (1-q) \hat{P}_q M[f](x)$. The modified Jackson
integral~$\bar{J}$, defined by $\bar{J}[f](x) =
(1-q)\hat{P}_q[f](x)$, satisfies the relation
\begin{eqnarray*}
\bar{J}[f]\,\bar{J}[g] + (1-q)\bar{J}[f\,g] = \bar{J}\big[ f\,
\bar{J}[g] \big] + \bar{J}\big[ \bar{J} [f] \, g\big].
\end{eqnarray*}
For motivational reasons we remark that the map $\hat{P}_q$ is of
importance in the construction of $q$-analogs of
multiple-zeta-values. The examples motivate the generalization of
the dual relation between the integration-by-parts rule and the
Leibniz rule for the classical calculus.

\begin{defn}
A \textit{Rota--Baxter map}~$R$ of weight~$\theta\in\R$ on a not
necessarily associative algebra~$A$, commutative or not, is a
linear map $R:A\to A$ fulfilling the condition
\begin{equation}
R(a)R(b) = R(R(a)b) + R(aR(b)) - \theta R(ab), \qquad a,b \in A.
\label{eq:brokeback}
\end{equation}
The reader will easily verify that $\tilde{R}:=\theta\,\id-R$ is a
Rota--Baxter map of the same weight, as well. We call a pair $(A,R)$,
where $A$ is an algebra and $R$ a Rota--Baxter map of weight~$\theta$,
a \textit{Rota--Baxter algebra} of weight~$\theta$. The indication
`not necessarily associative' is indispensable in this paper, as we
soon meet Rota--Baxter algebras that are neither Lie nor associative.
\end{defn}

We state a few simple observations, which will be of use later. The
so-called \textit{double Rota--Baxter} product
\begin{equation}
x *_R y := xR(y) + R(x)y - \theta xy, \qquad x,y \in A,
\label{def:doubleRBprod}
\end{equation}
endows the vector space underlying $A$ with another Rota--Baxter
algebra structure, denoted by $(A_R,R)$. In fact, $R$ satisfies
the Rota--Baxter relation for the new product. One readily
shows, moreover:
\begin{equation}
R(x *_R y) = R(x)R(y) \sepword{and} \tilde{R}(x *_R y) = -
\tilde{R}(x)\tilde{R}(y),\qquad x,y \in A.
\label{eq:doubleRBhom}
\end{equation}
This construction may be continued, giving a hierarchy of
Rota--Baxter algebras.

\begin{prop}
Let $(A,R)$ be an associative Rota--Baxter algebra of weight
$\theta \in \R$. The Rota--Baxter relation extends to the Lie
algebra $A$ with the commutator $[x,y]:=xy-yx$,
$$
[R(x), R(y)] + \theta R\bigl([x, y]\bigr) = R\bigl([R(x), y] + [x,
R(y)]\bigr)
$$
making $(A,[.,.],R)$ into a Rota--Baxter Lie algebra.
\label{prop:LieRB}
\end{prop}

This is a mere algebra exercise. A more exotic result coming next will
prove to be important in the context of Magnus' expansion and beyond.

\begin{prop}
Let $(A,R)$ be an associative Rota--Baxter algebra of weight
$\theta \in \R$. The binary composition
\begin{equation}
a \cdot_R b := [a,R(b)] + \theta ba
\label{eq:preLie}
\end{equation}
defines a right pre-Lie (or Vinberg) product such that $A$ becomes a
Rota--Baxter right pre-Lie algebra.
\label{prop:preLieRB}
\end{prop}

\begin{proof}
Recall that for a pre-Lie algebra $(A,\cdot)$ the (right) pre-Lie
property is weaker than associativity
$$
a \cdot (b \cdot c) - (a \cdot b) \cdot c = a \cdot( c \cdot b) - (a
\cdot c) \cdot b ,\quad \forall a,b,c \in A,
$$
As the Jacobiator is the total skewsymmetrization of the associator,
the pre-Lie relation is enough to guarantee that the commutator $[a,
b]:=a\cdot b-b\cdot a$ satisfies the Jacobi identity. For the sake of
brevity we verify only the weight-zero case and leave the rest to the
reader.
\begin{align*}
a \cdot_R (b \cdot_R c) &- (a \cdot_R b) \cdot_R c
=\big[a, R([b, R(c)])\big] - \big[[a,R(b)], R(c)\big] =
\\
&=\big[a, R([b, R(c)])\big] + \big[[R(c), a], R(b)\big] + \big[a,
[R(c), R(b)]\big]
\\
&=\big[a ,R([c, R(b)])\big] - \big[[a, R(c)], R(b)\big] =: a \cdot_R (c
\cdot_R b) - (a \cdot_R c) \cdot_R b;
\\
\sepword{and} R(a) \cdot_R R(b) &= [R(a), R(R(b))] = R\bigl([R(a),
R(b)]\bigr) + R\bigl([a, R(R(b))]\bigr)
\\
&= R\bigl( R(a) \cdot_R b\bigr) + R\bigl(a \cdot_R R(b)\bigr).
\end{align*}
Here we used Proposition~\ref{prop:LieRB} as well as the Jacobi
identity.
\end{proof}

The Lie algebra bracket corresponding to the double Rota--Baxter
product~(\ref{def:doubleRBprod}) is the double Rota--Baxter Lie
bracket $[a,b]_R:= a *_R b - b *_R a = a \cdot_R b - b \cdot_R a$,
known since the work of Semenov-Tian-Shansky~\cite{STS83}. We should
mention that these little calculations become more transparent using
the link between associative Rota--Baxter algebras and Loday's
dendriform algebras~\cite{Loday01,KEF}.

As a corollary to the last propositions we add the following
identity which will also be useful later
\begin{align}
R(a \cdot_R b) &= R([a,R(b)]) + \theta ba) = R([b,R(a)]) + [R(a),R(b)]
+ \theta R(ab)
\nonumber \\
&=R(b \cdot_R a) + [R(a),R(b)],
\label{pre-Lie-wow}
\end{align}
which is another way of saying that
$$
R([a,b]_R) = R(a *_R b - b *_R a) =[R(a),R(b)] =R(a \cdot_R b - b
\cdot_R a).
$$

\smallskip

The triple $(A,\dl,R;\theta)$ will denote an algebra~$A$ endowed with
a skewderivation~$\dl$ and a corresponding Rota--Baxter map~$R$, both
of weight $\theta$, such that $R\dl a=a$ for any $a\in A$ such that
$\dl a\ne0$, as well as $\dl Ra=a$ for any $a\in A,Ra\in0$. We check
consistency of the conditions~\eqref{eq:brokeback}
and~\eqref{eq:crash} imposed on $R,\dl$. Respectively
\begin{align*}
\theta\dl R(ab) &= R(a)b + aR(b) - \dl(R(a)R(b)) = R(a)b + aR(b) -
R(a)b - aR(b) + \theta ab = \theta ab;
\\
R\dl(ab) & = R(a\dl(b)) + R(\dl(a)b) - \theta R(\dl(a)\dl(b)) =
R(a\dl(b)) + R(\dl(a)b)
\\
&- R(a\dl(b)) - R(\dl(a)b) + ab = ab.
\end{align*}
The moral of the story is that Rota--Baxter maps are
\textit{generalized integrals}, skewderivations and Rota--Baxter
operators being natural (partial) inverses. As an example we certainly
have $(C(\R),d/dt,\int;0)$, with~$\dl=$ the derivative (with only the
scalars in its kernel). Another example is given by the aforementioned
triple $(A,\dl,Z;-\theta)$ of the finite difference map~$\dl$ of
step~$\theta$ and the summation~$Z$ in~\eqref{eq:le-clou}.

\smallskip

Rota--Baxter algebras have attracted attention in different contexts,
such as perturbative renormalization in quantum field theory (see
references further below) as well as generalized shuffle relations in
combinatorics~\cite{EbLe}. A few words on the history of the
Rota--Baxter relation are probably in order here. In the 1950's and
early 1960's, several interesting results were obtained in the
fluctuation theory of probability. One of the better known is
Spitzer's classical identity~\cite{Spitzer} for sums of independent
random variables. In an important 1960~paper~\cite{Baxter}, the
American mathematician G.~Baxter developed a combinatorial point of
view on Spitzer's result, and deduced it from the above operator
identity~\eqref{eq:brokeback}, in the context where the algebra~$A$ is
associative, unital and commutative. Then G.-C. Rota started a careful
in depth elaboration of Baxter's article in his~1969
papers~\cite{Rota1,Rota2}, where he solved the crucial ``word
problem'', and in~\cite{RotaSmith}, where he established several
important results. During the~1960's and~1970's, further algebraic,
combinatorial and analytic aspects of Baxter's identity were studied
by several people, see~\cite{Kingman,CartierBaxter,fields} for more
references. Recently, the Rota--Baxter relation became popular again
as a key element of the
Connes--Kreimer~\cite{DirkChen,ConnesKrRHI,EbKr} algebraic approach to
renormalization.

At an early stage the mathematician F.~V.~Atkinson made an important
contribution, characterizing such algebras by a simple decomposition
theorem.

\begin{thm}{\rm{(Atkinson~\cite{Atkinson})}}
Let $A$ be an algebra. A linear operator $R:A\to A$ satisfies the
Rota--Baxter relation~\eqref{eq:brokeback} if and only if the
following two statements are true. First, $A_+:=R(A)$ and $A_-:=
(\theta\,\id-R)(A)$ are subalgebras in~$A$. Second, for $X,Y,Z\in A$,
$R(X)R(Y)=R(Z)$ implies $(\theta\,\id-R)(X)(\theta\,\id-R)(Y) =
-(\theta\,\id-R)(Z)$.
 \label{Atkinson1}
\end{thm}

This result degenerates in the case $\theta=0$, whereby
$R=-\tilde{R}$. A trivial observation is that every algebra is a
Rota--Baxter algebra (of weight~1); in fact, the identity map and
the zero map are a natural Rota--Baxter pair. The case of an
idempotent Rota--Baxter map implies $\theta=1$ and, more
importantly, $A_{-}\cap A_{+}=\{0\}$, corresponding to a direct
decomposition of~$A$ into the image of~$R$ and~$\tilde{R}$.

Atkinson made another observation, formulating the following theorem,
which describes a multiplicative decomposition for associative unital
Rota--Baxter algebras.

\begin{thm}
\label{Atkinson2}
Let $A$ be an associative complete filtered unital Rota--Baxter
algebra with Rota--Baxter map~$R$. Assume~$X$ and~$Y$ in~$A$ to solve
the equations
\begin{equation}
X = 1_{A} + R(a\, X) \sepword{and} Y = 1_{A} + \tilde{R}(Y\, a),
\label{atkinsonEqs}
\end{equation}
for $a\in A^1$. Then we have the following factorization
\begin{equation}
Y (1_{A} - \theta a) X = 1_{A}, \sepword{so that } 1_A - \theta a
= Y^{-1}X^{-1}. \label{eq:Atkinsonfact1}
\end{equation}
For an idempotent Rota--Baxter map this factorization is unique.
\end{thm}

\begin{proof}
First recall that a complete filtered algebra $A$ has a decreasing
filtration $\set{A^n}$ of sub-algebras
$$
A=A^0 \supset A^1 \supset \dots \supset A^n \supset \dots
$$
such that $A^mA^n\subseteq A^{m+n}$ and $A\cong\invlim A/A^n$,
that is, $A$ is complete with respect to the topology determined
by the $\set{A^n}$. Also, note that
$$
R(a)\tilde{R}(b) = \tilde{R}(R(a)b) + R(a\tilde{R}(b)),
$$
and similarly exchanging $R$ and~$\tilde{R}$. Then the product
$YX$ is given by \allowdisplaybreaks{
\begin{eqnarray*}
YX &=&\bigl(1_A + \tilde{R}(Y\, a)\bigr)  \bigl(1_A +
R(a\,X)\bigr) = 1_A + R(a\,X) + \tilde{R}(Y\,a) +
\tilde{R}(Y\,a)R(a\,X)
\\
&=& 1_A + \tilde{R}\Bigl(Y\,a\,\bigl(1_A + R(a\,X)\bigr)\Bigr) +
R\Bigl(\bigl(1_A + \tilde{R}(Y\,a)\bigr)\,a\,X\Bigr)
\\
&=& 1_A + R(Y\,a\,X) + \tilde{R}(Y\,a\,X) = 1_A + \theta X\,a\,Y.
\end{eqnarray*}}
Hence we obtain the factorization~\eqref{eq:Atkinsonfact1}. Uniqueness
for idempotent Rota--Baxter maps is easy to show~\cite{EGMbch06}.
\end{proof}

In summary, finite difference as well as $q$-difference equations play
a role in important applications; thus it is useful to consider
generalizations of the classical apparatus for solving differential
equations. In the next section, by exploiting and complementing the
CBHD development of the previous one, we make preparations to extend
the work by Magnus on exponential solutions for non-autonomous
differential equations to general Rota--Baxter maps, beyond the
Riemann integral.

\section{The Spitzer identities and the CBHD \textit{recursion}}

In the last section we mentioned Spitzer's classical identity as a
motivation for Baxter's work. Now we spell out what that is.
Spitzer's identity can be seen as a natural generalization of the
solution of the simple initial value
problem~\eqref{eq:gato-encerrado} on the commutative algebra $A$
of continuous functions over $\R$,
\begin{equation}
\frac{df(t)}{dt} = a(t)f(t), \quad f(0) = 1, \quad a \in A.
\label{eq:IVP}
\end{equation}
This has, of course, a unique solution $f(t)=\exp\left(\int_0^t
a(u) \,du\right)$. Transforming the differential equation into an
integral equation by application of the Riemann integral~$I: A \to
A$ to~\eqref{eq:IVP},
\begin{equation}
f(t) = 1 + I(af)(t),
 \label{eq:integralIVP1}
\end{equation}
we arrive naturally at the not-quite-trivial identity
\begin{equation}
\exp\left(\int_0^t a(u)\,du\right) = \exp\bigl(I(a)(t)\bigr) = 1 +
\sum_{n=1}^\infty \underbrace{I\Bigl(a I\bigl(a\cdots I(a)\cdots
\bigr)\Bigr)}_{n\mbox{\rm -times}}(t).
 \label{eq:exp-sol}
\end{equation}
Taking into account the weight-zero Rota--Baxter
rule~\eqref{eq:integ-by-parts} for~$I$, the last identity follows
simply from
\begin{equation}
\bigl(I(a)(t)\bigr)^n = n!\underbrace{I\Bigl(a I\bigl(a\cdots
I(a)\cdots \bigr)\Bigr)}_{n\mbox{\rm -times}}(t).
 \label{eq:BSzero}
\end{equation}

Let now $(A,R)$ to be a \textit{commutative} Rota--Baxter algebra
of weight $\theta\neq0$. We formulate Spitzer's finding in
the~ring of power series $A[[t]]$, which is a complete filtered
algebra with the decreasing filtration given by the powers
of~$t,\,A^{n}:=t^nA[[t]],\, n\geq 0$. Notice that the power series
algebra~$A[[t]]$ with the operator $\mathcal{R}:A[[t]]\to A[[t]]$
acting on a series via~$R$ through the coefficients,
$\mathcal{R}\left(\sum_{n\ge 0}a_nt^n\right):=\sum_{n\ge
0}R(a_n)t^n$, is Rota--Baxter as well. Then we have

\begin{thm}{\rm(Spitzer's identity)}
Let $(A,R)$ be a unital commutative Rota--Baxter algebra of weight
$\theta\ne0$. Then for $a\in A$,
\begin{equation}
\exp\left(\!-R\biggl(\frac{\log(1 - a\theta
t)}{\theta}\biggr)\!\!\right) =
\sum_{n=0}^\infty(t)^n\underbrace{R\Bigl(\!a R\bigl(a\cdots
R(a)\cdots \bigr)\!\Bigr)}_{n\mbox{\rm -times}}
\label{eq:clSpitzer}
\end{equation}
in the ring of power series $A[[t]]$.
\label{clSpitzer}
\end{thm}

Analytic as well as algebraic proofs of this identity can be found
in the literature, see for instance~\cite{RotaSmith,EbKr}; and
anyway it is a corollary of our work further below. Observe that
$-\theta^{-1} \log(1 - a\theta t)\xrightarrow{\theta\downarrow
0}at$. Thus indeed~\eqref{eq:clSpitzer}
generalizes~\eqref{eq:exp-sol}.

Moreover, identity~\eqref{eq:BSzero} generalizes to the
\textit{Bohnenblust--Spitzer formula}~\cite{RotaSmith} of
weight~$\theta$. This is as follows. Let $(A,R)$ be a commutative
Rota--Baxter algebra of weight $\theta$ and fix $s_1,\dots,s_n\in A,
\,n>0$. Let~$S_n$ be the set of permutations of $\{1,\dots, n\}$. Then
\begin{equation}
\sum_{\sigma\in S_n}R\Bigl(s_{\sigma(1)}R\bigl(s_{\sigma(2)}\cdots
R(s_{\sigma(n)})\cdots\bigr)\Bigr) =
\sum_{\mathcal{T}\in\Pi_n}\theta^{n-|\mathcal{T}|}\prod_{T\in
\mathcal{T}}(|T|-1)!\,R\Bigl(\prod_{j\in T}s_j\Bigr),
\label{eq:BohnenblustSp}
\end{equation}
Here $\mathcal{T}$ runs through all unordered set partitions of
$\{1,\dots,n\}$; by~$|\mathcal{T}|$ we denote the number of blocks
in~$\mathcal{T}$; by~$|T|$ the size of the particular block~$T$.
The Rota--Baxter relation itself appears as a particular case
for~$n=2$. The weight $\theta=0$ case reduces the sum
over~$\mathcal{T}$ to~$|\mathcal{T}|=n$:
\begin{equation}
\sum_{\sigma\in S_n}R\Bigl(s_{\sigma(1)}R\bigl(s_{\sigma(2)}\cdots
R(s_{\sigma(n)})\cdots\bigr)\Bigr) = \prod_{j=1}^{n}
R\bigl(s_j\bigr).
\end{equation}
Also, for $n>0$ and $s_1=\dots=s_n=x$ we find
in~\eqref{eq:BohnenblustSp}:
\begin{equation}
R\Bigl(x R\bigl(x\cdots R(x)\cdots \bigr)\Bigr) =
\frac{1}{n!}\sum_{\mathcal{T}\in\Pi_n}\theta^{n - |\mathcal{T}|}
\prod_{T \in \mathcal{T}}(|T|-1)!\, R\bigl( x^{|T|}\bigr).
\label{eq:BohnenblustSp2}
\end{equation}
Relation~\eqref{eq:BohnenblustSp} follows from Spitzer's
identity~\eqref{eq:clSpitzer} by expanding the logarithm and the
exponential on the left hand side, and comparing order by order the
infinite set of identities in $A[[t]]$.

Spitzer's classical identity constitutes therefore an interesting
generalization of the initial value problem~\eqref{eq:IVP},
respectively the integral equation~\eqref{eq:integralIVP1}, to
more general integration-like operators~$R$, satisfying the
identity~\eqref{eq:brokeback}. Again we refer the reader
to~\cite{RotaSmith,fields} for examples of such applications in
the context of renormalization in perturbative quantum field
theory, $q$-analogs of classical identities, classical integrable
systems and multiple zeta values. Also, Atkinson's factorization
Theorem~\ref{Atkinson2} is obvious from Spitzer's identity. The
right hand side of identity~\eqref{eq:clSpitzer} is a solution to
$X=1_A+tR(a\, X)$ in~$A[[t]]$ corresponding to the factorization
of the element~$1_A - \theta a t$. (One ought to be careful here,
since Spitzer's identity as well as~\eqref{eq:BSzero} are only
valid for commutative Rota--Baxter algebras of weight~$\theta$,
whereas Atkinson's factorization result applies to general
associative unital Rota--Baxter algebras.)

\smallskip

Let us adopt an even more general point of view. For functions
with image in a noncommutative algebra, say $n\x n$ matrices with
entries in~$\R$, relation~\eqref{eq:exp-sol} is not valid anymore
as a solution to~\eqref{eq:IVP}; nor is identity~\eqref{eq:BSzero}
valid. From our present perspective, however, it does seem quite
natural to approach the problem of finding a solution to the
initial value problem, as well as relations~\eqref{eq:exp-sol}
and~\eqref{eq:BSzero}, on a noncommutative function algebra~$A$ by
looking for a generalization of Spitzer's identity to
noncommutative unital associative Rota--Baxter algebras of
weight~$\theta$. This latter problem was finally solved
in~\cite{E-G-K1,E-G-K2} ---see also~\cite{EGMbch06}, where the
reader may find more detail and earlier references. We will review
briefly those results, prior to extend our findings by indicating
a noncommutative generalization of the Bohnenblust--Spitzer
formula.

We then take the first steps towards the \textit{noncommutative
Spitzer} identity. Let~$A$ be a complete filtered associative
algebra. Bring in from Section~6 the
Campbell--Baker--Hausdorff--Dynkin~(CBHD) formula for the product
of exponentials of two non-commuting objects~$x,y$
$$
\exp(x)\exp(y) = \exp\bigl(x + y + \CBHD (x,y)\bigr),
\sepword{where} \sum_{m\ge 2}\Phi_m(x,y) =: \CBHD (x,y).
$$
Now let $P:A\to A$ be \textit{any} linear map preserving the
filtration and $\tilde{P}=\theta\,\id-P$, with $\theta$ an arbitrary
nonzero complex parameter. For $a\in A^1$, define the nonlinear map
$$
\chi^{\theta,\tilde{P}}(a) =
\lim_{n\to\infty}\chi^{\theta,\tilde{P}}_{(n)}(a)
$$
where $\chi^{\theta,\tilde{P}}_{(n)}(a)$ is given by the so-called
CBHD \textit{recursion},
\allowdisplaybreaks{
\begin{eqnarray}
\chi^{\theta,\tilde{P}}_{(0)}(a) &:=& a, \nonumber
\\
\chi^{\theta,\tilde{P}}_{(n+1)}(a) &=& a - \frac{1}{\theta}\,\CBHD
\bigl(\tilde{P}(\chi^{\theta,\tilde{P}}_{(n)}(a)),
P(\chi^{\theta,\tilde{P}}_{(n)}(a))\bigr), \label{eq:chik}
\end{eqnarray}}
and where the limit is taken with respect to the topology given by
the filtration. Then the map $\chi^{\theta,\tilde{P}}:A^1\to A^1$
satisfies
\begin{equation}
\chi^{\theta,\tilde{P}}(a) = a - \frac{1}{\theta}\,
\CBHD\bigl(\tilde{P}(\chi^{\theta,\tilde{P}}(a)),
P(\chi^{\theta,\tilde{P}}(a))\bigr). \label{eq:BCHrecursion1}
\end{equation}
We call $\chi^{\theta,\tilde{P}}$ the CBHD recursion of
weight~$\theta$, or just the $\theta$-CBHD recursion. In the following
we do not index the map $\chi^\theta(a):= \chi^{\theta,\tilde{P}}$ by
the operator~$\tilde{P}$ involved in its definition, when it is
obvious from context. One readily observes that $\chi^\theta$ reduces
to the identity for commutative algebras.

The following theorem states a general decomposition on the
algebra~$A$ implied by the CBHD recursion. It applies to associative
as well as Lie algebras.

\begin{thm}
Let $A$ be a complete filtered associative (or Lie) algebra with a
linear, filtration preserving map $P:A\to A$ and $\tilde{P}:=
\theta\,\id-P$. For any $a\in A^1$, we have
\begin{equation}
\exp({\theta}a) = \exp\bigl(\tilde{P}(\chi^{\theta}(a))\bigr)
\exp\bigl(P(\chi^{\theta}(a))\bigr). \label{eq:bch}
\end{equation}
Under the further hypothesis that the map~$P$ is idempotent (and
$\theta=1$), we find that for any $x\in 1_A+A^1$ there are unique
$x_{-}\in \exp\bigl(\tilde{P}(A^{1})\bigr)$ and $x_{+}\in
\exp\bigl(P(A^{1})\bigr)$ such that $x=x_-\,x_+$. \label{thm:bch}
\end{thm}

For proofs we refer the reader to~\cite{EGMbch06}. Using this
factorization one simplifies~\eqref{eq:BCHrecursion1} considerably.

\begin{lema}
\label{simpleCHI}
Let $A$ be a complete filtered algebra and $P:A \to A$ a linear map
preserving the filtration, with $\tilde{P}$ as above. The map
$\chi^\theta$ in~\eqref{eq:BCHrecursion1} solves the following
recursion
\begin{equation}
\chi^{\theta}(u) := u + \frac{1}{\theta}\CBHD \bigl(\theta u,-
P(\chi^{\theta}(u))\bigr), \quad u \in A^1.
\label{eq:BCHrecursion2}
\end{equation}
\end{lema}

The convolution algebra $(\Hom(H,A),*)$, for $H$ a connected graded
commutative Hopf algebra, will also be complete filtered. We may
immediately apply the above factorization theorem, giving rise to a
factorization of the group $G_H(A)$ of $A$-valued characters, upon
choosing any filtration-preserving linear map on $\Hom(H,A)$. In fact,
we find for $\theta=1$ in the definition of~$\chi$ the following
result.

\begin{prop}
\label{prop:factHopf}
Let~$A$ be a commutative algebra and~$H$ a connected graded
commutative Hopf algebra. Let $P$ be any filtration preserving linear
map on $\Hom(H,A)$. Then we have for all $\phi \in G_H(A)$ and
$Z:=\log(\phi) \in g_H(A)$ the characters $\phi^{-1}_{-}:=
\exp\bigl(\tilde{P}(\chi(Z))\bigr)$ and $\phi_{+}:=
\exp\bigl(P(\chi(Z))\bigr)$ such that
\begin{equation}
\phi = \phi^{-1}_{-} * \phi_{+}.
\label{eq:fact}
\end{equation}
If $P$ is idempotent this decomposition is unique.
\end{prop}

A natural question is whether one can find closed expressions for
the map $\chi^\theta$. The answer is certainly affirmative in some
non-trivial particular cases~\cite{EGMbch06}.

\begin{corl}
In the setting of the last proposition we find for the particular
choice of $P=\pi_-: H \to H$ being the projection to the odd
degree elements in $H$ (hence $\theta=1$)
\begin{equation*}
\chi(Z) = Z + \CBHD\Big(Z,-\pi_{-}(Z) - \thalf\CBHD\big(Z,Z -
\pi_{-}(Z) \big)\Big), \quad Z \in g_H(A).
\end{equation*}
\end{corl}

Before proceeding, we must underline that the factorization is due
solely to the map $\chi^\theta$; in fact, the map~$P$
---respectively $\tilde{P}$--- involved in its definition has
only to be linear and filtration preserving. The role played by
this map is drastically altered when we assume it moreover to be
\textit{Rota--Baxter} of weight~$\theta$, on a complete filtered
Rota--Baxter algebra. This we do next, to rederive and generalize
the Magnus expansion. Also, it will soon become clear
what~$\chi^0$ is. One of the aims of this paper is to attack the
solution of the CBHD recursion when $P$ is Rota--Baxter.

\smallskip

We noted earlier Atkinson's multiplicative decomposition of
associative complete filtered Rota--Baxter algebras. Let from now on
$(A,R)$ denote one such, of weight~$\theta\neq0$. Observe the
useful identity
\begin{equation}
\theta \prod_{i=1}^{n}R(x_i) = R\Bigl(\prod_{i=1}^{n} R(x_i) -
\prod_{i=1}^{n} \tilde{R}(-x_i)\Bigr), \sepword{for} x_i\in
A,\;i=1,\dots,n. \label{eq:theIdentity}
\end{equation}
This comes from the Rota--Baxter relation~\eqref{eq:brokeback}.
The case $n=2$ simply returns it. The reader should check it with
the help of the double Rota--Baxter
product~(\ref{def:doubleRBprod}). In the following we consider
$\chi^{\theta}:=\chi^{\theta,\tilde{R}}$ on~$A^1$.
Using~\eqref{eq:theIdentity}, for $\theta^{-1}\log(1_A-\theta
a)=:u \in A^1$ one readily computes \allowdisplaybreaks{
\begin{eqnarray*}
&&\exp\bigl(-R(\chi^\theta (u))\bigr) = 1_A + \sum_{n>0}
\frac{R(-\chi^\theta (u))^{n}}{n!}
\\
&&= 1_A + \sum_{n>0} \frac{(-1)^n}{n!\theta}R\bigl(R(\chi^\theta
(u))^n - (-\tilde{R}(\chi^\theta (u))^n)\bigr)
\\
&&= 1_A + \frac{1}{\theta} R\Bigl(\exp\bigl(-R(\chi^\theta
(u)\bigr) - \exp\bigl(\tilde{R}(\chi^\theta (u)\bigr)\Bigr) = 1_A
+ R\Bigl(a \exp\bigl(-R(\chi^\theta (u)\bigr)\Bigr).
\end{eqnarray*}}
In the last step we employed the factorization
Theorem~\ref{thm:bch} corresponding to~$\chi^{\theta}$. Therefore,
on the one hand we have found that
$X:=\exp\bigl(-R(\chi^\theta(u))\bigr) \in1_A+A^1$ solves
$X=1_A+R(a\, X)$, one of Atkinson's recursions in
Theorem~\ref{Atkinson2}. On the other hand, a solution to this
recursion follows from the iteration
\begin{equation}
X = 1_A + \sum_{n>0}\underbrace{R\bigl(a R(a R(a\cdots
R(a)}_{n\;\mathrm{times}})\dots )\bigr). \label{eq:Chen1}
\end{equation}
Hence:

\begin{thm}
\label{OurSpitzer} The natural generalization of Spitzer's
identity~\eqref{eq:clSpitzer} to noncommutative complete filtered
Rota--Baxter algebras of weight $\theta\neq0$ is given by
\begin{equation}
\exp\biggl(- R\biggl(\chi^{\theta}\biggl(\frac{\log(1_A - \theta
a)}{\theta}\biggr)\biggr)\biggr) = \sum_{n=0}^\infty
\underbrace{R\bigl(a R(a R(a\cdots R(a)}_{n\;\mathrm{times}})\dots
)\bigr), \label{eq:NC-SpitzerId-theta}
\end{equation}
for $a\in A^1$. Recall that $\chi^{\theta}$ reduces to the identity
for commutative algebras, yielding Spitzer's classical identity.
\end{thm}

So far we have achieved the following. First we derived the
general factorization Theorem~\ref{thm:bch} for complete filtered
algebras, upon the choice of an arbitrary linear filtration
preserving map. Specifying the latter to be Rota--Baxter of
weight~$\theta$, that is, identity~\eqref{eq:theIdentity}, we have
been able to show that Atkinson's recursion equations in
Theorem~\ref{Atkinson2} have exponential solutions. It is now
natural to ask whether the Bohnenblust--Spitzer
formula~\eqref{eq:BohnenblustSp} valid for weight-$\theta$
commutative Rota--Baxter algebras can be generalized to
noncommutative ones. The answer is yes! We outline next this
generalization, postponing detailed proof to the
forthcoming~\cite{KuruJosePatras}, to keep this long work within
bounds. First, by using the pre-Lie~\eqref{eq:preLie} and the
double~(\ref{def:doubleRBprod}) Rota--Baxter products, we find
$$
R\bigl(x_1R(x_2)\bigr) + R\bigl(x_2R(x_1)\bigr) = R(x_1)R(x_2) +
R(x_2 \cdot_R x_1)= R(x_1 *_R x_2) + R(x_2 \cdot_R x_1).
$$
Recall the relations~\eqref{eq:doubleRBhom}. One may now check by a
tedious calculation that
\allowdisplaybreaks{
\begin{eqnarray}
&&\sum_{\sigma\in S_3}
R\Bigl(x_{\sigma_1}R\bigl(x_{\sigma_2}R(x_{\sigma_3}) \bigr)\Bigr)
= R(x_1 *_R x_2 *_R x_3) + R((x_2 \cdot_R x_1) *_R x_3)
\nonumber \\
\!\!&&\!\! + R((x_3 \cdot_R x_1)*_R x_2) + R(x_3 \cdot_R (x_2 \cdot_R
x_1)) + R(x_1 *_R (x_3 \cdot_R x_2)) + R(x_2 \cdot_R (x_3 \cdot_R
x_1))
\nonumber \\
\!\!&&\!\! = R(x_1)R(x_2)R(x_3) + R(x_2 \cdot_R x_1)R(x_3) + R(x_3
\cdot_R x_1)R(x_2)
\nonumber \\
\!\!&&\!\! + R(x_3 \cdot_R (x_2 \cdot_R x_1)) + R(x_1)R(x_3
\cdot_R x_2) + R(x_2 \cdot_R (x_3 \cdot_R x_1)).
\label{eq:Ave-Maria-Santisima}
\end{eqnarray}}
We obtain special cases of the above when $x_1=x_2=x_3=x$
\begin{eqnarray*}
\!\!\!&&\!\!\!3!R\Bigl(x R\bigl(x R(x)\bigr)\Bigr) = R(x *_R x *_R
x) + 2 R((x \cdot_R x) *_R x) + 2R(x \cdot_R ( x \cdot_R x))
\\
\!\!\!&+&\!\!\! R(x*_R(x \cdot_R x)) = R(x)^3 + 2 R(x \cdot_R x)R(x) +
2R(x \cdot_R ( x \cdot_R x)) + R(x) R(x \cdot_R x).
\end{eqnarray*}
Equation~\eqref{eq:Ave-Maria-Santisima} is an instance of the
following result, that seems to be new:

\begin{thm}
\label{prop:factorization}
Let $(A,R)$ be an associative Rota--Baxter algebra of weight
$\theta$. For $x_i\in A,\,i=1,\dots,n$, we have
\begin{eqnarray*}
\sum_{\sigma\in S_n}R\Bigl(x_{\sigma_1} R\bigl(x_{\sigma_2} \dots
R(x_{\sigma_n})\dots\bigr)\Bigr) = \sum_{\sigma\in S_n}
R\Bigl(x_{\sigma_1} \diamond_1 x_{\sigma_2} \diamond_2 \dots
\diamond_n x_{\sigma_n}\Bigr)
\end{eqnarray*}
where
\begin{eqnarray*}
x_{\sigma_i} \diamond_i x_{\sigma_{i+1}} = \begin{cases} x_{\sigma_i}
*_R x_{\sigma_{i+1}}, & \sepword{$\sigma_i < {\rm{min}}(\sigma_j|
i<j)$}
\\
x_{\sigma_i} \cdot_R x_{\sigma_{i+1}}, & \sepword{otherwise;}
\end{cases}
\end{eqnarray*}
furthermore consecutive $\cdot_R$ products should be performed
>from right to left, and always before the $*_R$ product.
\end{thm}
This is the \textit{noncommutative Bohnenblust--Spitzer} formula.
Obviously, in the commutative case, i.e., when $a\cdot_R b=\theta ab$,
we just recover the classical Bohnenblust--Spitzer
identities~\eqref{eq:BohnenblustSp} and~\eqref{eq:BohnenblustSp2}. On
the other hand, anticipating on coming sections, the case $\theta=0$
reduces to Lam's factorization theorem~\cite{Lam}, stated in the
context of a weight-zero Rota--Baxter algebra of operator valued
functions $(B,\int;0)$. Let us adopt Lam's notation for the $n$-fold
right bracketed pre-Lie product by
\begin{equation}
C^R_n:= C^R_n(x) := R(x \cdot_R (x \cdot_R \dots (x \cdot_R
x)\dots)). \label{eq:LamsC}
\end{equation}
Also, we introduce a notation for the so-called Rota--Baxter
words:
\begin{equation}
(Rx)^{[n+1]} = R\bigl(x(Rx)^{[n]}\bigr), \label{eq:RBwords}
\end{equation}
with the convention that $(Rx)^{[0]}=1$.
Then we obtain the general expression
\begin{eqnarray}
(Rx)^{[n]} = \sum_{l=1}^n
\sum_{\substack{k_1,\dots ,k_l \in \mathbb{N}^\ast \\ k_1 + \dots +
k_l = n}}\,\frac{C^R_{k_1}\cdots C^R_{k_l}}{k_l(k_{l-1} + k_l) \cdots
(k_1 + \dots + k_l)}.
\label{eq:RBlevel}
\end{eqnarray}
That is to say, we sum over the compositions of~$n$. The simplest
cases already examined now are written
\begin{eqnarray}
2!(Rx)^{[2]} := 2!R\bigl(xR(x)\bigr) &=& (C^R_1)^2 + C^R_2,
\nonumber \\
3!(Rx)^{[3]} := 3!R\Bigl(xR\bigl(xR(x)\bigr)\Bigr) &=& (C^R_1)^3 +
2 C^R_2C^R_1 + C^R_1C^R_2 + 2C^R_3.
\label{eq:Lam-First}
\end{eqnarray}
Later we make use of those expansions in relation with the Magnus
expansion and the Dyson--Chen series. In fact, the left hand side of
the above expressions are the second and third order terms in the
path- or time-ordered expansion, in the context when the map~$R$ is
the Riemann integral. We just generalized this to general-weight
Rota--Baxter operators.

\section{The zero-weight recursion}

Let us come back to the CBHD recursion~$\chi^{\theta}$. The question
of the limit~$\theta\downarrow0$ becomes subtler than in the
commutative case, due to the particular properties of
relation~\eqref{eq:BCHrecursion1}. Now, in general we may write
$\Phi(a,b) = \Phi_1(a,b) + \CBHD (a,b)$ as a sum
$$
\Phi(a,b) = \sum_{n\ge1}H_n(a,b),
$$
where $H_n(a,b)$ is the part of $\Phi(a,b)$ which is homogenous of
degree $n$ with respect to $a$. For $n=1$ we
have~\cite{Reutenauer}:
\begin{equation}
H_1(a,b) = \frac{\hbox{ad}\,b}{ e^{\hbox{\eightrm ad}\,b} - 1_A}(a) =
\frac{\hbox{ad}\,b}{2}\Bigl(\coth\frac{\hbox{ad}\,b}{2} - 1\Bigr)(a).
\label{eq:manes-de-Bernouilli}
\end{equation}
In the limit $\theta\downarrow0$ all higher order terms $H_{n>1}$
vanish and from~\eqref{eq:BCHrecursion2} we get a nonlinear
map~$\chi^0$ inductively defined on~$A^1$ by the formula
\allowdisplaybreaks{
\begin{equation}
\chi^0(a) = \frac{\hbox{ad} P\bigl(\chi^0(a)\bigr)}
{{\rm{e}}^{\hbox{\eightrm ad} P(\chi^0(a))} - 1_A}(a) = \biggl(1_A
+ \sum_{n>0}b_n\Big[\hbox{\rm ad} P\bigl(\chi^0(a)\bigr)\Big]^{n}
\biggr)(a).
\label{eq:BCH-recur0}
\end{equation}}
We may call this the \textit{weight-zero {\rm CBHD} recursion}. The
coefficients are $b_n:=B_n/n!$ with~$B_n$ the Bernoulli numbers. For
$n=1,2,4$ we find $b_1=-1/2,b_2=1/12$ and $b_4=-1/720$. We have
$b_3=b_5=\cdots=0$. The first terms in~\eqref{eq:BCH-recur0} are then
easily written down:
\begin{align}
&\chi^0(a) = a - \frac 12[P(a),\,a] + \frac14 \big[P\bigl([P(a),
a]\bigr), \,a\big] + \frac{1}{12}\big[P(a), [P(a), a]\big]
\label{eq:pre-Magnus} \\
& -\frac{1}{24}P\Bigl(\bigl[P([P(a),[P(a),a]]),a\bigr] +
\bigl[P(a),[P([P(a),a]),a]\bigr] +
\bigl[[P([P(a),a]),[P(a),a]\bigr ]\Bigr)
\nonumber \\
& -\frac{1}{8} P\Bigl(\bigl[P([P([P(a),a]),a]),a\bigr]\Bigr) + \cdots.
\nonumber
\end{align}
We pause here to note that~\eqref{eq:manes-de-Bernouilli} is but
an avatar of the formula
$$
D(e^a) = e^a\,\frac{1 - e^{-\hbox{\eightrm ad}\,a}}{\hbox{ad}\,a}Da,
\sepword{for $D$ a derivation;}
$$
a noncommutative chain rule familiar from linear group theory.
See~\cite[Chapter~1]{Rossmann} for instance. Apparently this is due
to~F.~Schur~(1891), and was taken up later by Poincar\'e and
Hausdorff. One may also consult the charming account of the
determination of a local Lie group from its Lie algebra, using
canonical coordinates of the first kind,
in~\cite[Chapter~13]{SudarshanM}. The appearance of the Bernoulli
numbers is always fascinating. The deep reason for it is that we are
trying to express elements of the enveloping algebra in terms of the
symmetric algebra.

Now suppose $P$ is a weight-zero Rota--Baxter operator, denoted~$R$
henceforth. The noncommutative generalization of Spitzer's identity in
the case of vanishing weight is captured in the following corollary.

\begin{corl}
\label{corl:0-ncSpitzer}
Let $(A,R)$ be a complete filtered Rota--Baxter algebra of weight
zero. For $a$ in~$A^1$ the weight-zero $\CBHD$ recursion $\chi^0:A^1\to
A^1$ is given by equation~\eqref{eq:BCH-recur0}:
\begin{equation*}
\chi^0(a) = \frac{\hbox{\rm ad}R\bigl(\chi^0(a)\bigr)}
{{\rm{e}}^{\hbox{\eightrm ad}R(\chi^0(a))} - 1_A}(a).
\end{equation*}
Moreover:
\begin{enumerate}
\item\label{eq:exp1o} The equation $x=1+R(a\, x)$ has a unique
solution $x=\exp\bigl(R(\chi^0(a))\bigr)$. \item\label{eq:exp2o}
The equation $y=1-R(y\, a)$ has a unique solution
$y=\exp\bigl(R(-\chi^0(a))\bigr)$.
\end{enumerate}
\end{corl}

We will see pretty soon that the $0$-CBHD recursion gives Magnus'
expansion. The diagram~\eqref{eq:diag-general} further below
summarizes the foregoing relations, generalizing the simple initial
value problem~\eqref{eq:IVP} in a twofold manner. First we go to the
integral equation~\eqref{eq:integralIVP1}. Then we replace the Riemann
integral by a general Rota--Baxter map and assume a noncommutative
setting. That is, we start with a complete filtered noncommutative
associative Rota--Baxter algebra $(A,R)$ of nonzero weight~$\theta$ in
the appropriate field. The top of the diagram~\eqref{eq:diag-general}
contains the solution to the equation
\begin{equation}
X = 1_A + R(a\,X), \sepword{for} a\in A^1, \label{eq:recursion}
\end{equation}
generalized to associative, otherwise arbitrary Rota--Baxter
algebras~\eqref{eq:NC-SpitzerId-theta},
\begin{equation}
X = \exp\biggl(-R\biggl(\chi^{\theta}\biggl(\frac{\log(1_A - \theta
a)} {\theta}\biggr)\biggr)\biggr). \label{eq:solution}
\end{equation}
The $\theta$-CBHD recursion $\chi^{\theta}$ is given
in~\eqref{eq:BCHrecursion1}. The left wing of~\eqref{eq:diag-general}
describes the case when first the weight $\theta$ goes to zero, hence
reducing $\chi^\theta\to\chi^0$; see~\eqref{eq:BCH-recur0}. This is
the algebraic structure underlying Magnus's~$\Omega$ series of the
next section. Then, we let the algebra~$A$ become commutative, which
reduces~$\chi^0$ to~$\id$. The right wing of
diagram~\eqref{eq:diag-general} just describes the alternative
reduction, i.e., we first make the algebra commutative, which gives
the classical Spitzer identity for nonzero weight commutative
Rota--Baxter algebras, see~\eqref{eq:clSpitzer}. Then we take the
limit~$\theta\downarrow0$.
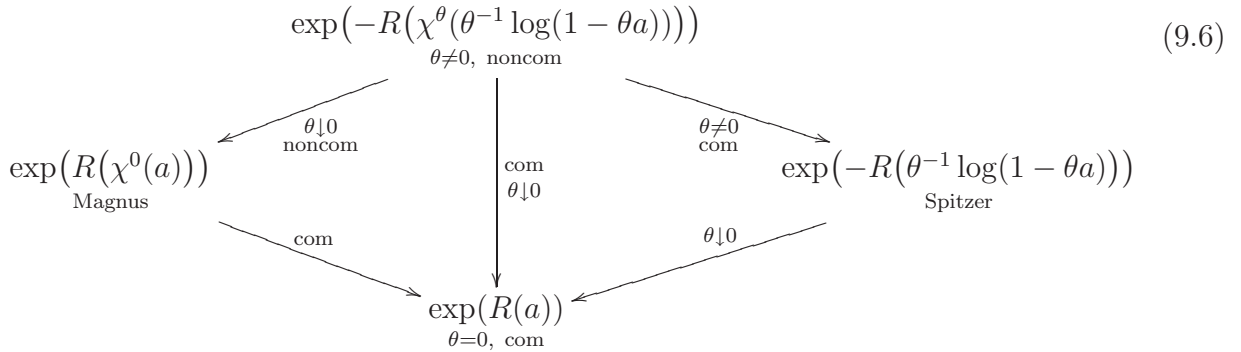
\begin{equation}
\xymatrix{
&
\underset{\theta \neq 0,\ \mathrm{noncom}}
         {\exp\bigl( -R\bigl( \chi^\theta(
          \theta^{-1} \log(1 - \theta a) ) \bigr) \bigr) }
\ar[dd]^{\substack{\mathrm{com} \\[\jot] \theta \downarrow 0}}
\ar[rd]_{\substack{\theta \neq 0 \\ \mathrm{com}}}
\ar[ld]^{\substack{\theta \downarrow 0 \\ \mathrm{noncom}}}
\\
\underset{\mathrm{Magnus}}
         {\exp \bigl( R\bigl(\chi^0(a)\bigr) \bigr) }
\ar[rd]^{\mathrm{com}}
& &
\underset{\mathrm{Spitzer}}
         {\exp \bigl( -R \bigl(\theta^{-1} \log(1 - \theta a)
          \bigr) \bigr) }
\ar[ld]_{\theta \downarrow 0}
\\
&
\underset{\theta = 0,\ \mathrm{com}}{\exp(R(a))}
}
\label{eq:diag-general}
\end{equation}
Both paths eventually arrive at the elementary fact that
equation~\eqref{eq:recursion} is solved by a simple exponential in a
commutative, weight-zero Rota--Baxter setting. We have succeeded in
finding the general algebraic structure underlying the initial value
problem for generalized integrals, that is, Rota--Baxter operators.

\section{On Magnus' commutator series}

It is high time for us to declare why we choose to deal with first
order non-autonomous differential equations primarily via the Magnus
expansion method. The latter has a somewhat chequered history. To
attack the initial value problem of the
type~\eqref{eq:gato-encerrado}:
\begin{equation}
\frac{d}{dt}F(t) = a(t)F(t), \qquad F(0) = 1,
\label{eq:gato-con-siete-vidas}
\end{equation}
with $F$ a matrix-valued function, say, Magnus proposed the
exponential Ansatz
$$
F(t) = \exp\bigl(\Omega[a](t)\bigr),
$$
with $\Omega[a](0)=0$. He found a series for~$\Omega[a]$:
\begin{equation}
\Omega[a](t) = \sum_{n>0}\Omega_{n}[a](t),
\label{eq:Magnus-Exp}
\end{equation}
in terms of multiple integrals of nested commutators, and provided
a differential equation which in turn can be easily solved
recursively for the terms $\Omega_{n}[a](t)$
\begin{equation}
\frac{d}{dt}\Omega[a](t) = \frac{\hbox{ad}\,
\Omega[a]}{{\rm{e}}^{\hbox{\eightrm ad}\,\Omega[a]} - 1}(a)(t).
\end{equation}
It is worth indicating that originally Magnus was motivated by
Friedrich's theorem of our Section~5. We already mentioned,
however, that one of the papers most influential on the
subject~\cite{inflationbuster} was written without knowledge of
Magnus' paper. In the 1990's, several mathematicians interested in
approximate integrators for differential equations developed the
discipline of geometrical integration. Originally also unaware of
Magnus' work, they derived anew Magnus' expansion. The point was to
make sure that the approximate solutions evolve in the Lie group
if~$\xi(t)$ in~\eqref{eq:auto-gato} remains in the Lie algebra.
This is \textit{not} true of the iterative Dyson--Chen method
---no finite truncation of the latter is the exact solution of any
approximating system. By construction $\chi^0$ respects the Lie
algebra structure in~\eqref{eq:auto-gato}, and thus truncations of
the series are sure to remain in the Lie group. This is one reason
why ---in view also of the considerations in Appendix~C--- we give
priority to the Magnus method. For geometrical integration,
consult~\cite{IsNo,IserlesAMS}.

Comparison with~\eqref{eq:BCH-recur0} settles the matter of the link
between Magnus series and the CBHD recursion in the context of
vanishing Rota--Baxter weight. Namely,

\begin{corl}
Let $A$ be a function algebra over~$\R_t$ with values in an operator
algebra. Let~$R$ denote the indefinite Riemann integral operator.
Magnus'~$\Omega$ expansion is given by the formula
\begin{equation*}
\Omega[a](t) = R\bigl(\chi^0(a)\bigr)(t).
\end{equation*}
\end{corl}

In conclusion, we could say that the $\theta$-CBHD
recursion~\eqref{eq:BCHrecursion1} generalizes Magnus' expansion
to general weight $\theta\neq0$ Rota--Baxter operators~$R$ by
replacing the weight-zero Riemann integral in $F=1+R\{aF\}$.
Corollary~\ref{corl:0-ncSpitzer} represents a more modest
generalization, to zero-weight Rota--Baxter operators different
>from the ordinary integral.

\smallskip

Let us write explicitly the first few terms of the Magnus
expansion using~\eqref{eq:pre-Magnus}, when $R$ is the Riemann
integral operator. The function $a=a(t)$ is defined over $\R$ and
takes values in a noncommutative algebra, say of matrices of size
$n\x n$. We obtain
\allowdisplaybreaks{
\begin{align}
&R(a)(t) = \int_0^t a(t_1)\,dt_1,
\label{eq:chapucilla} \\
&-\frac 12 R\bigl([R(a),\,a]\bigr)(t) = \frac
12\int_0^t\int_0^{t_1}[a(t_1), a(t_2)]\,dt_2\,dt_1,
\nonumber \\
&\frac 14 R\Bigl(\big[R\bigl([R(a), \,a]\bigr),\,a\big]\Bigr)(t) =
\frac 14 \int_0^t\int_0^{t_1}\int_0^{t_2}
\big[[a(t_3),\,a(t_2)],\,a(t_1)\big] dt_3\,dt_2\,dt_1,
\nonumber \\
&\frac{1}{12} R\Bigl(\big[R(a), [R(a), a]\big]\Bigr)(t) = \frac
{1}{12} \int_0^t\int_0^{t_1}\int_0^{t_1} \big[a(t_3), [a(t_2),
a(t_1)]\big]\,dt_3\,dt_2\,dt_1. \nonumber
\end{align}}
This gives indeed the first terms of the expansion \textit{precisely}
in the form that Magnus derived it. However, in later
works~\cite{inflationbuster,MexicoLindo,MexicoOtraVez,Wilcox,Salzman}
the terms in the Magnus' expansion are presented as iterated
commutator brackets of strictly time-ordered Riemann integrals.
Especially in~\cite{inflationbuster} Strichartz succeeded in giving a
closed solution to Magnus' expansion ---and hence to our recursion
$\chi^{0,R}$ when $R$ is the Riemann integral. With the notation of
Proposition~\ref{pr:NotQuick}, he found
\begin{align}
\Omega[a](t) &= \sum_{n > 0}\Omega_n[a](t), \sepword{with}
\label{eq:StrichartzMagnus} \\
\Omega_n[a](t) &= \sum_{\sg\in S_n}
\frac{(-1)^{d(\sg)}}{n^2\binom{n-1}{d(\sg)\,}}
\int_0^t\int_0^{t_1}\dots \int_0^{t_{n-1}} \big[[\dots
[a(t_1),a(t_2)]\dots], a(t_n)\big]\,\,dt_n \dots dt_2\,dt_1.
\nonumber
\end{align}
This formula clearly points to the close relation between the CBHD
expansion and Magnus' series, although the appearance of the number of
descents is still `unexplained'. However, it is not to everyone's
taste. It is immediate from the formula that
$$
\Omega_2[a](t) = \frac 12\int_0^t\int_0^{t_1}[a(t_1), a(t_2)]
\,dt_2\,dt_1,
$$
coincident with the second term in~\eqref{eq:chapucilla}; and clear enough
that
\begin{equation}
\Omega_3[a](t) =\frac{1}{6}\int_0^t \int_0^{t_1}\int_0^{t_{2}}
\Bigl(\big[[a(t_1), a(t_2)], a(t_3)\big] - \big[[a(t_2),a(t_3)],
a(t_1)\big] \Bigr)\,dt_3\,dt_2\,dt_1;
\label{eq:StrichartzMagnusThree}
\end{equation}
however the number of terms grows menacingly with~$n!$, it is
never evident when we will find cancellations, and one quickly
concludes that the beauty of~\eqref{eq:StrichartzMagnus} hides its
computational complexity. Nor is it entirely obvious, although of
course it is true, that~\eqref{eq:StrichartzMagnusThree} coincides
with the sum of the third and fourth terms
in~\eqref{eq:chapucilla}.

The best policy, in our opinion, is to invoke the alternative
Dyson--Chen solution at this point. This attacks three problems:
systematic writing of the Magnus series simplifies; the zero-weight
recursion is solved; and the comparison between different expressions
for the same terms is made easier.

\section{Enter the Dyson--Chen series}

The first order initial value problem~\eqref{eq:gato-con-siete-vidas},
respectively the corresponding recursion
$$
F(t) = 1 + R(aF)(t),
$$
where $R$ is the Riemann integral operator, possess a natural
solution in terms of iteration, see~\eqref{eq:Chen1}. The
resulting infinite series is called here Dyson--Chen integral. In
the physics literature those series are often referred to as
time-ordered exponentials or path-ordered integrals; their
importance can hardly be overstated. To reflect such nomenclature
in the notation, write
$$
\T e^{\int_0^t a(t_1)\,dt_1} = \T e^{R(a)(t)} := 1 +
\sum_{n>0}\underbrace{R\bigl(aR(aR(a\cdots
R(a))}_{n\;\mathrm{times}}\dots)\bigr)(t).
$$
The operator $\T$ implies the strict iteration of the integral
corresponding to the `time ordering'. A short presentation of
Chen's work on this kind of integrals can be found in~\cite{ShSt};
the findings of Magnus and Chen played a decisive role, especially
for Rota and his followers. Directly from the group property of
the flow, we have for the Dyson--Chen integral the factorization
$$
\T e^{\int_0^t a(t_1)dt_1} = \T e^{\int^{t'}_0 a(t_1)\,dt_1}
\,\T e^{\int_{t'}^t a(t_1)\,dt_1},
$$
giving rise to many identities of integrals and concatenation
products of series, which we need not go into. This factorization
might be compared with the quite different decomposition induced
by the CBHD recursion~\eqref{eq:bch}. The major result of the
theory is the following theorem.

\begin{thm}
\label{thm:Wanderlust}
The logarithm of a Chen series is a Lie series.
\end{thm}

The direct proof of this statement uses Hopf algebra, to wit, the
shuffle product algebra of our Section~5. It is just a matter of
verifying Ree's condition~\eqref{eq:Ree} inductively. In our present
context, the Dyson--Chen expansional is the solution to Atkinson's
recursion~\eqref{eq:recursion}, and the theorem scarcely needs
justification.

Simply by taking the logarithm in
$$
\exp(\Omega[a](t)) = \T e^{R(a)(t)},
$$
we obtain
\begin{equation}
\Omega_n[a] = \sum_{k=1}^{n}\frac{(-1)^{k+1}}{k} \sum_{\substack{l_1,
\dots ,l_k \in \mathbb{N}^\ast \\ l_1 + \cdots + l_k=n}} (Ra)^{[l_1]}
\dots (Ra)^{[l_k]}.
\label{eq:Chen-Magnus}
\end{equation}
This was of course known to the practitioners
---see~\cite{Salzman,KlOt} and references there.
It is derived in~\cite{Quaoar} by use of the Fa\`a di Bruno Hopf
algebra. Inverting these relations, one finds the $(Ra)^{[n+1]}$'s
in terms of the~$\Omega_m[a]$'s
\begin{equation}
(Ra)^{[n]} = \sum_{k=1}^{n} \frac{1}{k!} \sum_{\substack{l_1,\dots,
l_k \in \mathbb{N}^\ast \\ l_1 + \cdots + l_k=n}}
\Omega_{l_1}[a] \dots \Omega_{l_k}[a].
\label{eq:Magnus-Chen}
\end{equation}
The first examples are:
\begin{eqnarray*}
2!(Ra)^{[2]} &=& \Omega_1^2[a] + 2\Omega_2[a],
\\
3!(Ra)^{[3]} &=& \Omega_1^3[a] + 3\bigl(\Omega_1^2[a]\Omega_2[a] +
\Omega_2[a]\Omega_1^2[a]\bigr) + 6\Omega_3[a].
\end{eqnarray*}
that might be compared with~\eqref{eq:Lam-First}; of course
$Ra=\Omega_1[a] =C^R_1$ in the occasion. Now, both sets of
equations~\eqref{eq:Chen-Magnus} and~\eqref{eq:Magnus-Chen} simply
describe how to link Magnus' expansion to the Dyson--Chen expansional.
They purely follow from the Rota--Baxter relation as well as the CBHD
formula. Therefore they are valid for \textit{any} weight-zero
Rota--Baxter operator~$R$.

By inverting the Rota--Baxter map, we solve moreover the zero-weight
CBHD~recursion:
\begin{align*}
\chi^0_n(a) &= \sum_{k=1}^{n}\frac{(-1)^{k+1}}{k}
\sum_{\substack{l_1,l_2, \dots,l_k \in \mathbb{N}^\ast \\ l_1 +
l_2 + \cdots + l_k=n}} \bigl(a(Ra)^{[l_1-1]} (Ra)^{[l_2]} \cdots
(Ra)^{[l_k]} +
\\
&+ (Ra)^{[l_1]} a(Ra)^{[l_2-1]} \cdots (Ra)^{[l_k]} + \cdots +
(Ra)^{[l_1]} (Ra)^{[l_2]} \cdots a(Ra)^{[l_k-1]}\bigr).
\end{align*}

Next we are set to give an alternative formula
to~\eqref{eq:StrichartzMagnus}, keeping left-to-right bracketing. This
is better explained by way of example. Bring in Heaviside's step
function,
\allowdisplaybreaks{
\begin{eqnarray*}
\Th_{1,2}(t_1,t_2) := \Th(t_1-t_2) := \begin{cases} 1,\mbox{
if } t_1-t_2>0, \\ 0,\;{\rm otherwise}; \end{cases}
\end{eqnarray*}}
with its help, iterated Riemann integrals can be rewritten
\begin{equation}
R\bigl(aR(b)\bigr)(t) = \int_0^t
a(t_1)\int_0^{t_1}b(t_2)\,dt_2\,dt_1 =
\int_0^t\int_0^t\Th(t_1-t_2)a(t_1)b(t_2)\,dt_2\,dt_1.
\label{eq:Strichartz}
\end{equation}
More generally,
$$
\Th_{i,j}(t_1,t_2,\ldots,t_n) := \Th(t_i-t_j),
\sepword{for $1\le i,j\le n$,}
$$
and we can write
$$
\T e^{\int_0^t a(t_1)dt_1} = 1 + \int_0^t a(t_1)\,dt_1 +
\sum_{n=2}^\infty\int_0^t\cdots\int_0^t\Th_{1,2}\cdots\Th_{n-1,n}\,
a(t_1)\cdots a(t_n)\,dt_n \dots dt_1.
$$
For instance, for the third term of the Magnus series,
applying~\eqref{eq:Chen-Magnus},
$$
\Omega_3[a] =
\int_0^t\int_0^t\int_0^t\bigl(\Th_{1,2}\Th_{2,3} -
\thalf\Th_{1,2} - \thalf\Th_{2,3} + \tthird
\bigr)a(t_1)a(t_2)a(t_3) \,dt_3\,dt_2\,dt_1.
$$
Now, we know ---if only from theorem~(\ref{thm:Wanderlust})--- this is
a Lie element, so we can apply at once the Dynkin operator to rewrite
it with nested commutators:
$$
\Omega_3[a] =
\int_0^t\int_0^t\int_0^t\bigl(\Th_{1,2}\Th_{2,3} -
\thalf\Th_{1,2} - \thalf\Th_{2,3} + \tthird \bigr)[[a(t_1),
a(t_2)], a(t_3)]\,dt_3\,dt_2\,dt_1.
$$
We see now that the last term actually does not contribute to the
integral. With very little work, just using
$\Th_{1,2}+\Th_{2,1} =\Th_{2,3}+\Th_{3,2}=1$, one
recovers~\eqref{eq:StrichartzMagnusThree}. An explicit formula for
all terms along these lines, fully equivalent to, but simpler to
work with, than Strichartz's, is easily
obtained~\cite{MexicoOtraVez}; we do not bother to write it. We
must avow, however, that we do not see a way to write terms like
the third one in the integral above as a combination of iterations
and products of the~$R$ operators; thus we must conclude that
formulae like~\eqref{eq:StrichartzMagnus}
and~\eqref{eq:StrichartzMagnusThree} are only valid for the
Riemann integral.

\smallskip

For general zero-weight Rota--Baxter operators we may fall back
on~\eqref{eq:Chen-Magnus}. Magnus himself did not use any property of
the map~$R$ beyond integration-by-parts, and only presented the
expansion in a form equivalent to~\eqref{eq:pre-Magnus}. Of course,
even using purely the weight-zero Rota--Baxter relation, there are
many equivalent ways of writing the same. For instance, simply by
Proposition~\ref{prop:LieRB}, one finds that the term at third order
in Magnus' expansion~\eqref{eq:chapucilla} is rewritten
\begin{equation*}
\frac13 R\Bigl( \big[R\bigl([R(a), a]\bigr), a\big] \Bigr)(t) -
\frac{1}{12} \big[R\bigl([R(a), a]\bigr), R(a)\big](t).
\end{equation*}
It is worthwhile to mention that Iserles and Norsett use binary rooted
trees to achieve a better understanding of Magnus'
expansion~\cite{IsNo,IserlesAMS}.

\smallskip

We have long taken the algebraic tack. But what about convergence of
the Magnus series? Note that Dyson--Chen series converge absolutely
for all~$t$ if~$a$ is bounded, and this is why they are preferred in
quantum field theory; however this good property is not transmitted in
general to Magnus series via~\eqref{eq:Chen-Magnus}, as there is an
infinite resummation involved. Excellent bounds at small~$t$ have been
found recently~\cite{Dontmoan} for matrix systems. Strichartz
linearizes arbitrary initial-value problems, for which we cannot
expect convergence in general in the smooth category; but he does not
fail to observe that Magnus' expansion has especially good properties
for Lie--Scheffers systems~\cite[Section~3]{inflationbuster}. This is
because the closing of the involved vector fields to a
finite-dimensional Lie algebra sharply improves the estimates.
Furthermore, for those systems Magnus' exponential can be interpreted
as the exponential map of Lie~theory.

\section{Towards solving the $\theta$-weight recursions}

Let us now come back to Proposition~\ref{prop:factorization} and take
the first steps in going from the Dyson--Chen series to the
$\theta$-weight CBHD recursion. This looks somewhat hard; but recall
that Lam found, in the context of the Riemann integral, another way to
relate the terms in the Dyson--Chen series to those in the Magnus
expansion ---consult\cite{Lam,OteoRos}. In fact, Lam's findings are
true in a much more general sense, i.e., for general weight
Rota--Baxter algebras, as we will indicate here. The attentive reader
will remember the weight-zero pre-Lie product~\eqref{eq:preLie}, that
allows for the following way of writing the weight zero CBHD
recursion~$\chi^0(a)$, see~\eqref{eq:pre-Magnus}:
\begin{align}
\chi^0(a) &= a + \frac 12 a \cdot_R a + \Bigl(\frac14\bigl(a \cdot_R
(a \cdot_R a)\bigr) + \frac{1}{12} \bigl((a \cdot_R a) \cdot_R
a\bigr)\Bigr)
\nonumber \\
& +\frac{1}{24}R\Bigl(a \cdot_R \bigl( (a \cdot_R a)\cdot_R a\bigr) +
\bigl(a \cdot_R (a \cdot_R a)\bigr) \cdot_R a + (a \cdot_R a) \cdot_R
(a \cdot_R a)\Bigr)
\nonumber \\
& +\frac{1}{8} R\Bigl(a \cdot_R \bigl(a \cdot_R(a \cdot_R
a)\bigr)\Bigr) + \cdots.
\label{eq:pre-Lie-Magnus}
\end{align}
This contains in germ the main idea. Remember~\eqref{eq:LamsC} in
terms of the (double and) pre-Lie Rota--Baxter product. Lam made an
exponential Ansatz
$$
\sum_{n \geq 0} R(a)^{[n]} = \exp\Bigl(\sum_{m>0} K_m(a)\Bigr)
$$
and derived the following formulae for the~$K_i$'s in terms of
$C^R_1(a),\dots,C^R_i(a)$:
\begin{eqnarray*}
K_1(a) &=& C^R_1(a), \quad K_2(a) = \frac{1}{2}C^R_2(a), \quad K_3(a)
= \frac{1}{3}C^R_3(a) + \frac{1}{12}[C^R_2(a),C^R_1(a)],
\\
K_4(a) &=& \frac{1}{4}C^R_4(a) + \frac{1}{12}[C^R_3(a),C^R_1(a)],\dots
\end{eqnarray*}
The weight-$\theta$ Rota--Baxter relation enters at the level of
identity~\eqref{eq:RBlevel}, hence implying the particular form of
the~$K_i$'s. This naturally demands a comparison with the CBHD
recursion, respectively the generalized Spitzer identity.

\begin{thm}
Let $(A,R)$ be an associative Rota--Baxter algebra of
weight~$\theta$. Then for
$K_i=K_i(C^R_1(a),\dots,C_i^R(a);\theta)$ we have
$$
\sum_{i>0}K_it^i = -R\bigl(\chi^{\theta}\bigl(\theta^{-1}\log(1_A
- \theta at)\bigr)\bigr).
$$
\label{prop:solvingCHIviaRB}
\end{thm}

Hence, finding a formula for the $K_i$'s gives a solution, in the
sense of a closed expression, to the CBHD
recursion~$\chi^{\theta}$, which follows from the Rota--Baxter
relation. A full proof of this statement lies beyond the scope of
this work and will be provided elsewhere~\cite{KuruJosePatras}. In
the context of Hopf and Rota--Baxter theory, it is the
generalization of the shuffle relation to the quasi-shuffle (or
mixed-shuffle) identity~\cite{EbGu} underlying the algebraic
structure encoded in the Rota--Baxter relation of nonzero weight,
which generalizes the integration-by-parts
rule~\eqref{eq:integ-by-parts} corresponding to the shuffle
relation.

Of course, when~$\theta=0$ Lam's~$K_i$'s are just the
Magnus~$\Om_i$'s. These are expressed as sums of commutators of
\textit{right-to-left} bracketed integrals, when~$R$ is the Riemann
operator. This turns out to be the most efficient method for the
expansion, as well. For instance, the expression of~$K_5$ contains
just six terms, whereas $\Om_5$ is written usually with~22
terms~\cite{OteoRos}.

We close this section with a simple but striking observation flowing
>from the last theorem. Defining $u(at):=\theta^{-1}\log(1_A-\theta
at)$, we recover $-\chi^{\theta}(u(at))$ from~$\chi^0(at)$, that
is, from~\eqref{eq:pre-Lie-Magnus}, simply by using the
weight-$\theta$ pre-Lie product~\ref{eq:preLie}). A full proof of
this statement will be given elsewhere. But we show this here up
to third order. Using $\theta^{-1}\log(1_A-\theta at)=-\sum_{n>0}
\frac{{\theta}^{n-1}} {n}(at)^n$, we find for
$$
-\chi^{\theta}(u(at)) = at - \sum_{n>0}\chi_n^{\theta}(u(a))t^{n+1}
$$
the following
\allowdisplaybreaks{
\begin{eqnarray*}
\chi_{(1)}^{\theta}(u(a)) &=& \frac{1}{2}\theta a^2 -
\frac{1}{2}[R(a),a],
\\
\chi_{(2)}^{\theta}(u(a)) &=& \frac{1}{3}\theta^2 a^3 - \frac{1}{4}
\theta \bigl([R(a^2), a] + [R(a), a^2] \bigr) + \frac{1}{4}
\bigl[a,R\big([a,R(a)]\big)\bigr]
\\
&& + \frac{1}{12}\Bigl( \bigl[[a, R(a)], R(a)\bigr] - \theta
\bigl[a,[a,R(a)]\bigr] \Bigr).
\end{eqnarray*}}
Let us go back to~\eqref{eq:pre-Lie-Magnus} and use the pre-Lie
product $a\cdot_R b:=[a,R(b)]+\theta ba$ of~\eqref{eq:preLie}. We
obtain at second order
\allowdisplaybreaks{
\begin{eqnarray*}
\frac 12 a \cdot_R a &=& \frac 12 [a, R(a)] + \frac 12 \theta a^2.
\end{eqnarray*}}
At third order we calculate:
\allowdisplaybreaks{
\begin{align*}
&\frac14 \bigl(a \cdot_R (a \cdot_R a)\bigr) + \frac{1}{12} \bigl((a
\cdot_R a) \cdot_R a\bigr)
\\
& = + \frac{1}{4} [a, R(a \cdot_R a)] + \frac{1}{4} \theta (a \cdot_R
a)a + \frac{1}{12} [(a \cdot_R a), R(a)] + \frac{1}{12}\theta a(a
\cdot_R a)
\\
&= \frac{1}{4} \bigl[a,R\big([a, R(a)]\big)\bigr] + \frac{1}{4}\theta
[a, R(a^2)] + \frac{1}{4} \theta [a, R(a)]a + \frac{1}{4}\theta^2 a^3
\\
&\quad\ + \frac{1}{12} \bigl[[a ,R(a)], R(a)\bigr] + \frac{1}{12}
\theta [a^2, R(a)] + \frac{1}{12} a[a, R(a)] +
\frac{1}{12}\theta^2 a^3
\\
&= \frac{1}{3}\theta^2 a^3 + \frac{1}{12} \bigl[[a, R(a)], R(a)\bigr]
+ \frac{1}{4} \bigl[a, R\big([a, R(a)]\big)\bigr] + \frac{1}{4}\theta
\bigl([a, R(a^2)] + [a^2 ,R(a)]\bigr)
\\
&\quad\ + \frac{1}{4} \theta [a, R(a)]a - \frac{1}{6} \theta [a^2,
R(a)] + \frac{1}{12} \theta a[a, R(a)]
\\
&= \frac{1}{3}\theta^2 a^3 + \frac{1}{4}\theta \bigl([a, R(a^2)] +
[a^2, R(a)]\bigr) + \frac{1}{4} \bigl[a, R\big([a,R(a)]\big)\bigr]
\\
&\quad\ + \frac{1}{12} \bigl[[a, R(a)], R(a)\bigr] - \frac{1}{12}
\theta \bigl[a, [a, R(a)]\bigr].
\end{align*}}
Earlier in Section~10 we have seen how the Magnus expansion naturally
follows from the CBHD recursion in the limit $\theta\downarrow0$. In
turn we see here the advantage of reformulating Magnus' expansion in
terms of the Rota--Baxter pre-Lie product of weight~$\theta$ yielding
the CBHD~recursion.

\section{Conclusion and outlook}

Our purpose in this paper was twofold. Starting from the
innocent-looking dynamical system~\eqref{eq:gato-encerrado} ---of
classical Lie--Scheffers type when $G$ is an ordinary Lie group--- we
sought to reformulate it in Hopf algebraic terms, thus being led to
generalized derivation and integration (Rota--Baxter) operators.
Whereby we show that two of the three main ordinary strategies to
attack non-autonomous linear differential equations (linked
respectively to the names of Magnus and Dyson--Chen) still make sense
in the broader context. In particular, the noncommutative version of
the Bohnenblust--Spitzer identity has been found, and we blaze a trail
to solve the nonlinear recursion introduced earlier by one of us in
relation with the noncommutative Spitzer formula.

There is no doubt that the product integral method to
attack~\eqref{eq:gato-encerrado}, often linked with the name of
Fer~\cite{Fer}, is also susceptible to our kind of algebraic
reinterpretation and generalization. However, with a heavy heart, we
leave this for a later occasion: the present paper is already long
enough.

Needless to say, the programmatic purpose of this work was to
propagandize the Hopf algebra approach to differential equations. The
lure of presenting classical subjects under a new light explains why
we spent much space on a smooth transition from the standard to a
Hopf-flavoured view of dynamical systems; and indeed this article
became a powerful spur to revisit the traditional proof of the
Lie--Scheffers theorem, and plug its gaps~\cite{TresMosqueterosbis}.
On the other hand, many of our findings and procedures will surely not
raise an eyebrow of people working in sophisticated methods for
control theory ---on which we confess no expertise. There is, at any
rate, plenty left to do. Avenues open for possible research include:

\begin{itemize}
\item{} The Cari\~nena--Ramos' approach to Lie--Scheffers systems,
based on connections, should be recast in the noncommutative mould, in
the light of~\cite{Bigotes2} and~\cite{TheSecondComing}.

\item{} To relate and compare the action algebroid approach to
group \& Lie algebra actions with the Hopf algebra approach.

\item{} Investigation of the product integral method.

\item{} Further exploration of the theory of Rota--Baxter
operators as natural inverses to skewderivations; that is,
developing Rota's proposal of an algebraic theory of integration.

\item{} Definitive clarification of the noncommutative Spitzer
formula and the noncommutative Bohnenblust--Spitzer identity in the
light of Lam's expansion.

\item{} The bridge to control theory and chronological products,
via Loday's dendriform algebras in particular, should be enlarged
and strengthened. In this respect, Lie--Butcher
theory~\cite{LastWillBeFirst,LastWillBeSecond} shows great
promise.
\end{itemize}

\newpage

\appendix

\section{Pr\'ecis on group actions}

\begin{defn}
A (left) action of a Lie group~$G$ on a manifold~$M$ is a
homomorphism~$\Phi$ of~$G$ into~$\Diff M$. For $x\in M$, and $g\in G$
we denote
$$
\Phi_g := \Phi(g) \sepword{and} \Phi(g,x) := \Phi_g x.
$$
A right action is just an antihomomorphism of~$G$ into~$\Diff M$. The
orbits of~$\Phi$ are the subsets of~$M$ of the form~$\Phi(G,x)$ for a
fixed $x\in M$; they are homogeneous manifolds, on which the action is
transitive. We will call $\Phi_x$ the map from $G$ to~$M$ defined
by~$g\mapsto \Phi(g,x)$. Recall that a flow is an action of~$\R$
on~$M$. When $\Phi$ with the indicated properties is given, we say $M$
is a $G$-manifold. A Lie group action is proper if given any pair
$K,L$ of compacts subsets of~$M$, the set $\set{g\in G: gK\cap L \ne
\emptyset}$ is compact. The stabilizer or isotropy subgroups are then
compact. Proper actions, in particular compact group actions of
general Lie groups, have good properties: for instance the orbits of a
proper action are closed submanifolds of~$M$~\cite{MichorNotes}. An
action is faithful (or effective, or essential) when the map
$g\to\Phi_g$ is injective; if the kernel of this map is discrete, we
say the action is almost faithful.
\end{defn}

A good reference for Lie group actions
is~\cite[Chapter~4]{OldRedBook}. As for the examples, any Lie group
$G$ acts on itself by left and right translations $L_g,R_g:G\to G$
respectively given for each $g\in G$ by
$$
g'\mapsto gg',\qquad g'\mapsto g'g.
$$
The inverse diffeomorphisms are $L_g^{-1}=L_{g^{-1}}$
and~$R_g^{-1}=R_{g^{-1}}$. This action is free and transitive. Also
$G$ acts on itself by conjugation:
$$
g'\mapsto gg'g^{-1} =: \Ad(g)g'.
$$
This action is neither free nor transitive; it is almost faithful iff
the centre of~$G$ is discrete.

\begin{defn}
Suppose $G$ acts both on~$N$ by~$\Phi^N$ and on~$M$ by~$\Phi^M$. A
smooth map $f:N\to M$ between these manifolds is \textit{equivariant}
(with respect to the actions) if $f\circ\Phi^N_g= \Phi^M_g\circ f$ for
each $g\in G$. The maps $\Phi_x:G\to M$, where $\Phi$ is a left
(right) action are equivariant for all~$x$, with respect to the left
(right) action of $G$ on itself and~$\Phi$:
\begin{equation}
\Phi_x \circ L_g = \Phi_g \circ \Phi_x \sepword{or}
\Phi_x \circ R_g = \Phi_g \circ \Phi_x,
\label{eq:begs-the-question}
\end{equation}
as the case may be.
\end{defn}

If $G$ acts on~$M$, then $G$ also acts on~$TM$ by
$$
(g,v_x) \mapsto (\Phi_g x,T_x\Phi_gv_x) =: \Phi^T(g,v_x), \sepword{for
$v_x\in T_xM$.}
$$
When $\Phi$ is described in local coordinates, say by
$$
\Phi_i(g,x) = h_i(g,x^1,\ldots,x^n), \sepword{then} T_x\Phi_gv_x =
\sum_{j=1}^n\frac{\pd h_i}{\pd x^j} (g,x^1,\ldots,x^n)v_x^j.
$$
Clearly the map $v_x\mapsto\Phi^T(g,v_x)$ from $T_xM$ into
$T_{\Phi(g,x)}M$ is linear and the canonical projection $\tau_M:TM\to
M$ is equivariant with respect to these actions:
$\tau_M\bigl(\Phi^T(g,v_x)\bigr)=\Phi\bigl(g,\tau_M(v_x)\bigr)$. We
then say that $\Phi$ is \textit{equilinear}~\cite{DieudonneIV}. For
vector fields, then, there is the action:
\begin{equation}
(g,X) \mapsto T\Phi_g\circ X\circ\Phi_{g^{-1}}.
\label{eq:insidious-thing}
\end{equation}

Corresponding to group translations we have then equilinear left and
right actions of $G$ on~$TG$; as well as actions on~$\X(G)$. In view
of~\eqref{eq:insidious-thing}, a vector field $X^L$ on~$G$ is left
invariant if for all $g\in G$, $X^L\circ L_g=TL_g\circ X^L$; this
means that $X^L$ is $L_g$-related to itself for all $g\in G$.
Therefore the left invariant vector fields constitute a Lie
subalgebra~$\X^L(G)=:\g_L$ of~$\X(G)$. Replacing $L_g$ by $R_g$ we
obtain right invariant vector fields~$X^R\in\X^R(G)$ and a Lie
subalgebra~$\X^R(G)=:\g_R$. In particular, $X^L,X^R$ are determined by
their values in the neutral element:
\begin{equation}
X^L(g) = T_1L_gX^L(L_g^{-1}g) = T_1L_gX^L(1_G); \sepword{similarly}
X^R(g) = T_1R_gX^R(1_G);
\label{eq:morir-al-palo}
\end{equation}
for typographical simplicity we write~$T_1$ instead of~$T_{1_G}$. The
dimension of~$\X^L(G)$ or of~$\X^R(G)$ is thus that of the group. We
denote by~$X^L_\xi,X^R_\xi$ the left invariant, respectively right
invariant, vector field associated to $\xi\in T_1G$. The (complete)
flow of~$X^L_\xi$ is $(t,g)\mapsto g\exp(tX^L_\xi)$ and the flow
of~$X^R_\xi$ is $(t,g)\mapsto\exp(tX^R_\xi)g$.

We remark that $\g_L$ is the commutant of~$\g_R$ in~$\X(G)$, and
vice versa. For instance, thinking of the affine group of
orientation-preserving transformations of the line as a neighbourhood
of~$(1,0)$ with the multiplication rule:
$$
(x^1,x^2)\.(y^1,y^2) = (x^1y^1, x^1y^2 + x^2),
$$
then a basis for left (respectively right) invariant vector fields is
$$
(X^L_1, X^L_2) := (x^1\del_1, x^1\del_2); \sepword{respectively}
(X^R_1, X^R_2) := (x^1\del_1 + x^2\del_2, \del_2).
$$
With our Lie bracket, by the way: $[X^R_1, X^R_2]=X^R_2$. It is an
easy exercise to check that if $a_1(x^1,x^2)\del_1+ a_2(x^1,x^2)
\del_2$ commutes with $X^L_1,X^L_2$, then it is a linear combination
of~$X^R_1,X^R_2$ with scalar coefficients.

Consider the tangent map $T\imath:TG\to TG$ lifting the inversion
diffeomorphism $\imath:g\mapsto g^{-1}$ on the base; it carries
left invariant vector fields into right invariant ones. The vector
fields $T\imath\circ X^L_\xi$ and~$-X^R_\xi\circ\imath$
along~$\imath$ coincide, that is, $X^L_\xi$ is
$\imath$-projectable on $-X^R_\xi$. This simply because
$\bigl(g^{-1}\exp(tX^L_\xi)\bigr)^{-1}=\exp(-tX^R_\xi)g$.
Therefore $[X^L_\xi, X^L_\eta]$ projects into $[X^R_\eta,
X^R_\xi]$.

Now, $TG$ is itself a group, with product $T\mu$ lifted from the
product~$\mu:G\x G\to G$. The short exact sequence (where $T_1G$ is
the additive group of this tangent linear space)
$$
0 \to T_1G \to TG \to G \to 1
$$
splits, which means $TG\sim T_1G\rtimes G$, with $T_1G$ embedded
in~$TG$ as a normal subgroup. In particular $TG$ is a trivial vector
bundle. We have in~$TG$:
\begin{equation}
gv_{g'} = TL_gv_{g'}; \qquad v_{g'}g = TR_gv_{g'}.
\label{eq:maestros-sutiles}
\end{equation}
Clearly, the action of~$G$ on~$T_1G$ is just~$\Ad_{1_G}^T$. Henceforth
we write~$\Ad$ for this adjoint action of~$G$ on~$T_1G$. A Lie bracket
can now be defined directly on~$T_1G$ by $[\xi,\eta]:=
\ad(\xi)\eta:=T_1\Ad(\xi)\eta$. One could also transfer to~$T_1G$ the
Lie algebra structure from~$\X^L(G)$ or~$\X^R(G)$, say $[\xi,\eta]:=
[X^R_\xi,X^R_\eta](1)$. That these and other natural definitions
amount to the same is standard fare~\cite[Appendix~III]{Cirilo}. The
space $T_1G$ with any of these equivalent structures is what people
call the tangent (Lie) algebra~$\g$ of~$G$.

\begin{defn}
A left (right) Lie algebra (infinitesimal) action~$\la$ on~$M$ is a
Lie algebra homomorphism (antihomomorphism) $\g\ni\xi\mapsto\la_\xi
\in\X(M)$; we say $M$ is a $\g$-manifold. The action is said
\textit{transitive} at~$x$ when the $\la_\xi(x)$ span $T_xM$. It is
furthermore \textit{primitive} when the stabilizer $\g_x$ is a maximal
subalgebra; these concepts are analogous to the case of Lie group
actions. When the action is transitive at all points of~$M$, we say
infinitesimally transitive. Given $\la:\g\to\X(M)$, if $\lambda_\g$ is
made up of complete vector fields (in particular when~$M$ is compact,
guaranteeing completeness of all vector fields) and~$G$ is the simply
connected Lie group with Lie algebra~$\g$, then there is a unique
$\Phi:G\to\Diff M$ such that $T_1\Phi=\la$. This lifting to a group
action always exists locally. We remark as well that our choice of
sign for the bracket of vector fields insures that the derivative
$T_1\Phi$ of a left action is a left action.
\end{defn}

For the infinitesimal description of actions, the following notion is
essential.

\begin{defn}
Let $\Phi$ denotes an action of $G$ on $M$. For $\xi\in\g$, the map
$(t,x)\mapsto\Phi(\exp t\xi,x)$ is a flow on~$M$. The
\textit{fundamental vector field} or \textit{infinitesimal generator}
$\xi^\Phi_M$ of~$\Phi$ corresponding to~$\xi$ is the vector field
\begin{equation}
\xi^\Phi_M(x) := \ddto{t}\Phi(\exp t\xi,x)=T_1\Phi_x(\xi)\,.
\label{eq:acabaramos}
\end{equation}
The superscript~$\Phi$ is omitted in the notation when the action is
clear in the context. The image of~$\g$ under $T_1\Phi_x$ is the tangent
bundle $T(G\.x)$ of the $\Phi$-orbit. The corresponding differential
operator is given by
\begin{equation*}
\xi^\Phi_M f(x) := \ddto{t}f\bigl(\Phi(\exp t\xi,x)\bigr).
\end{equation*}
The anchor map $\xi\mapsto\xi^\Phi_M$ from the tangent algebra~$\g$
to~$\X(M)$ constitutes a Lie--Rinehart algebra; the corresponding Lie
algebroid will be transitive when the action of~$\g$ on~$M$ is
infinitesimally transitive.
\end{defn}

For example, when $\Phi$ is $L_g:G\to G$, we know that the
corresponding flow is $(t,g')\mapsto R_{g'}\exp t\xi$. Therefore
\begin{equation}
\xi_G(g') = T_1R_{g'}\xi = X^R_\xi(g'),
\label{eq:for-later-use}
\end{equation}
the \textit{right} invariant vector field associated to~$\xi$. By the
same token $\xi^R_G(g)=X^L_\xi(g)$.

If $M$ is a $G$-manifold, the flow of~$\xi_M$ is given by $\Phi_{\exp
t\xi}$. Indeed,
\begin{align*}
&\frac{d}{dt}\Phi(\exp t\xi,x) = \frac{d}{ds}\biggr|_{s=0}
\Phi(\exp(s+t)\xi,x)
\\
&= \frac{d}{ds}\biggr|_{s=0}\Phi(\exp s\xi,x)\circ\Phi(\exp t\xi,x)
= \xi_M\circ \Phi_{\exp t\xi}(x).
\end{align*}
As a consequence $\xi_M$ is complete. The reader will have little
difficulty in verifying the following

\begin{prop}
\label{pr:equiflows}
Let $N,M$ be $G$-manifolds with respective actions~$\Phi^N,\Phi^M$,
and $f:N\to M$ a smooth map equivariant with respect to these
actions; then $\xi_N\sim_f\xi_M$, that is $Tf\circ\xi_N=
\xi_M\circ f$. More precisely, $\xi_N\sim_f\xi_M$ iff the flows verify
$$
f\circ\Phi^N_{\exp(t\xi)}= \Phi^M_{\exp(t\xi)}\circ f.
$$
\end{prop}

\begin{prop}
\label{pr:wow}
For every $\xi,\eta\in\g$ we have
$$
[\xi_M, \eta_M] = [\xi, \eta]_M.
$$
In other words: $\xi\mapsto \xi_M$ is a left Lie algebra action.
\end{prop}

\begin{proof}
A simple calculation gives
$$
(\Ad_g\xi)_M = T\Phi_{g^{-1}}\xi_M.
$$
We obtain the result immediately by differentiation. Our
unconventional choice of sign for the Lie bracket of vector fields
avoids the obnoxious minus signs of the usual treatments.
\end{proof}

\smallskip

The action $\Phi$ of~$G$ on~$M$ lifts naturally to representations
of~$G$ on the various linear spaces associated with~$M$ ---for
instance to representations on spaces of sections of vector
bundles~\cite{ACO} or on morphisms of vector bundles. We will limit
ourselves to some simple cases, needed in the main text.
For~$f\in\F(M)$, we consider $(g\. f)(x) :=
f\bigl(\Phi(g^{-1},x)\bigr)$; then for~$T\in\F'(M)$ and
for~$D\in\D(M)$:
$$
\<g\. T, f> := \<T, g^{-1}\. f>, \sepword{respectively} (g\. D)f :=
g\. D(g^{-1}\. f).
$$
Invariant functions, distributions and differential operators are
defined in the obvious way.

\section{Differential equations on homogeneous spaces}

The problem of solving non-autonomous differential equations on
\textit{homogeneous spaces} of Lie groups is intimately linked to
Lie--Scheffers theory: given an arbitrary Lie group~$G$ and an action
of it on a manifold~$M$, for most purposes one can restrict oneself to
the orbits of the action, that is, the points of~$M/G$; these are
(immersed) submanifolds of~$M$ of the form~$G/G_x$,
with~$G_x:=\set{g\in G:\Phi(g,x)=x}$ the stabilizer of a point~$x$ of
the orbit, a closed Lie subgroup of~$G$. From this perspective,
Lie--Scheffers systems are precisely those that can be rewritten in
the form
\begin{equation}
\dot x(t)= \la_{\xi(t)}\bigl(x(t)\bigr),
\label{eq:yetLS}
\end{equation}
where $A:\R\to\g$ is a curve on the Lie algebra~$\g$ and~$\la$ denotes
an infinitesimal action. If $\la=T\Phi$ for some action~$\Phi$ of~$G$
on~$M$ and~$g(t)$ solves the initial value problem~\eqref{eq:fund}:
\begin{equation}
\dot g(t) = \xi_G(t,g(t)); \qquad g(0) = 1_G,
\label{eq:fund-bis}
\end{equation}
then the solution of~\eqref{eq:yetLS} with initial condition
$x(0)=x_0$ is given by the integrated action: $x(t)=\Phi(g(t),x_0)$.
At this point we again advise the reader to
consult~\cite{TresMosqueterosbis}.

In practice we consider transitive actions on~$M\equiv G/G_x$. Suppose
that $x_{(1)}$ is a particular solution of~\eqref{eq:yetLS} satisfying
$x(0)=x_0$. Let $g_1\in\Map(\R_t,G)$ such that $x_{(1)}(t)=
\Phi(g_1(t),x_0)$. Such curve is not unique in general; but, if $g_2$
is another one, then $g_2(t)=g_1(t)h(t)$ with $h$ in
$\Map(\R_t,G_{x_0})$. It is convenient to choose $h$ so that $g_2$ is
the fundamental solution of~\eqref{eq:yetLS}:
$$
\dot g_2(t) = T_1R_{g_2(t)}\xi(t),
$$
upon using~\eqref{eq:for-later-use} in the last equality. Then $h$ is
the fundamental solution of the Lie--Scheffers system associated to
the curve $B:\R\to\g_{x_0}$, given by~\cite{Luther}:
$$
B(t) = T_1L_{g_1(t)^{-1}}\bigl(T_1R_{g_1(t)}\xi(t) - \dot
g_1(t)\bigr).
$$
Therefore the knowledge of a particular solution of~\eqref{eq:yetLS}
that satisfies $x_{(1)}(0)=x_0$ reduces the problem of finding the
fundamental solution for~$G$ to finding the fundamental solution for
the subgroup $G_{x_0}$. Naturally if more particular solutions are
known, whose values at $0$ are $x_1,\dots,x_r$, then we can reduce the
problem to solving a Lie--Scheffers system in the subgroup
$G_{x_0}\cap\cdots\cap G_{x_r}$. When this group is discrete, one can
explicitly compute the fundamental solution for~$G$, from which the
general solution of the original Lie--Scheffers system can be derived.
This is known as the Lie reduction method.

A variant of the Lie reduction method was studied in the language of
gauge theory in~\cite{PBPepinArturo}. Without actually invoking
connections, we illustrate the approach in this last reference with
the Riccati equation~\eqref{eq:manes-de-Riccati}. The latter seeks the
integral curves of the vector field along~$\pi_2:\R_t\x M\to M$:
$$
\bar Y = \bigl(a_0(t) + a_1(t)x + a_2(t)x^2\bigr)\frac{\pd}{\pd x}.
$$
{}For vector fields $E_+=\pd/\pd x,H=x\,\pd/\pd x$ and
$E_-=x^2\,\pd/\pd x$ we observe the commutation relations
\begin{equation}
[H,E_+] = E_+; \qquad [E_+,E_-] = -2H; \qquad [H,E_-] = -E_-,
\label{eq:porca-miseria}
\end{equation}
exactly those of the matrices $E'_+ := \begin{pmatrix} 0 & 1 \\ 0 & 0
\end{pmatrix}$; $H' :=\begin{pmatrix} \thalf & 0 \\ 0 & -\thalf
\end{pmatrix}$; $E'_-:= \begin{pmatrix} 0 & 0 \\ -1 & 0
\end{pmatrix}$. Therefore $E_\pm,H$ realize the (perfect) Lie
algebra~$\sll(2;\R)$ of the group~$SL(2;\R)$. The corresponding flows
of~$\R$ are respectively
$$
x_0 \longmapsto x_0 + t; \qquad x_0 \longmapsto x_0e^t; \qquad x_0
\longmapsto \frac{x_0}{1-x_0t};
$$
the last one blows up for $x_0>0$ in finite time, indicating that
$E_-$ is not complete. This can be corrected by adding to~$\R$ the
point at infinity. More precisely, we have the well-known action of
the projective group $SL(2;\R)/Z_2$ on the projective
line~$\R\cup\infty$ ---to wit, the projectivization of the fundamental
action of $SL(2;\R)$ on~$\R^2$. Just as well, in the spirit of this
article, we can decide to regard the action as a local one, defined on
the open set of~$SL(2;\R)\x\R$ given by the pairs such that
$cx+d\ne0$.

Now, consider the group $\Map\bigl(\R_t,SL(2;\R)\bigr)$ of curves
acting on the set of Riccati equations (that is, the group of
automorphisms of the trivial principal bundle $SL(2;\R)\x\R_t\to\R_t$)
corresponding to the indicated action, expressed by:
\begin{equation*}
\Phi\bigl(A(t),x(t)\bigr) = \Phi\biggl(\begin{pmatrix} \a(t) & \b(t)
\\ \ga(t) & \dl(t) \end{pmatrix},x(t)\biggr) = \frac{\a(t)x(t) +
\b(t)}{\ga(t)x(t) + \dl(t)},
\end{equation*}
together with the other obvious cases. When $x(t)$ is a solution of
the Riccati equation~\eqref{eq:manes-de-Riccati}, then $x'(t):=
\Phi\bigl(A(t),x(t)\bigr)$ is also a solution of a Riccati equation
with coefficients
\begin{align*}
\begin{pmatrix} a'_2(t) \\ a'_1(t)  \\ a'_0(t) \end{pmatrix}
= \begin{pmatrix} \dl^2 & -\dl\ga & \ga^2 \\
-2\b\dl & \a\dl+\b\ga & -2\a\ga \\
\b^2 & -\a\b & \a^2 \end{pmatrix}
\begin{pmatrix} a_2(t) \\ a_1(t)  \\ a_0(t) \end{pmatrix}
+ \begin{pmatrix} \ga\dot\dl -\dl\dot\ga  \\
\dl\dot\a -\a\dot\dl +\b\dot\ga- \ga\dot\b  \\
\a\dot\b -\b\dot\a \end{pmatrix}.
\end{align*}
The second term on the right hand side is a 1-cocycle for the linear
action on the coefficients of the Riccati equation given by the first
term. If a particular solution $x_{(1)}(t)$
of~\eqref{eq:manes-de-Riccati} is known, the element $A_1(t) =
\begin{pmatrix} 1 & 0\\ -x_{(1)}^{-1}(t) & 1\end{pmatrix} \in
\Map\bigl(\R_t,SL(2;\R)\bigr)$, transforms the original Riccati
equation into the linear equation $dx'/dt=
\bigl(2x_{(1)}^{-1}(t)a_0(t)+a_1(t)\bigr)x'+a_0$, thereby reducing the
group $SL(2;\R)$ to the subgroup $A(1;\R)$. When a second particular
solution $x_{(2)}(t)$ of~\eqref{eq:manes-de-Riccati} is given, then
$x'=x_{(1)}x_{(2)}/(x_{(1)}-x_{(2)})$ satisfies the linear equation,
therefore we obtain the corresponding homogeneous linear equation
using the matrix $A_2 = \begin{pmatrix} 1 & -x_{(1)}x_{(2)}(x_{(1)} -
x_{(2)})^{-1}\\ 0 & 1\end{pmatrix}$. Concretely, the change of
variables
$$
x'' = \Phi(A_2,x') = \Phi(A_2A_1,x) = \frac{x_{(1)}^2(x-x_{(2)})}
{(x_{(2)}-x_{(1)})(x - x_{(1)})}
$$
leads to the homogeneous linear equation $dx''/dt=
\bigl(2x_{(1)}^{-1}(t)a_0(t)+a_1(t)\bigr)x''$. Finally, if $x_{(3)}$
is a third particular solution of~\eqref{eq:manes-de-Riccati}, then
$z=x_{(1)}^2(x_{(2)} -x_{(3)})/(x_{(2)} -x_{(1)})(x_{(1)} -x_{(3)})$
solves this linear equation, thus if $A_3 = \begin{pmatrix} z^{-1/2} &
0 \\ 0 & z^{1/2}\end{pmatrix}$, the transformation
$$
x'''=\Phi(A_3A_2A_1,x)
=\frac{(x -x_{(2)})(x_{(1)}-x_{(3)})}{(x-x_{(1)})(x_{(2)}-x_{(3)})}
$$
gives the reduced equation $dx'''/dt=0$, which is the superposition
principle~\eqref{eq:desmanes-de-LS} for the Riccati equation.

\smallskip

We are not likely to find an exact solution for~\eqref{eq:fund-bis} in
most cases. This is one reason why we concentrate on approximate
solutions in this paper. To attack~\eqref{eq:eqingr}, it is generally
a good strategy to move on to an equivalent system on the tangent
algebra of~$G$ ---a coordinate space for~$G$ which enjoys the
advantage of being a linear space. To effect properly the method of
working on the tangent algebra, one needs to ponder equivariant maps
between homogeneous spaces. We go to this in the next appendix.

\section{More on the same}

Consider again the canonical action of~$G$ on~its tangent
algebra~$\g$, and let $f:\g\to G$ be a local coordinate map. The
exponential map is an example, but of course there are slight
variants of it (see below); or we could employ, if available, the
Cayley map~\cite{KostantMichor}. A local action $B^f$ of~$G$
on~$\g$ is constructed by~$B^f_g = f^{-1}\circ L_g\circ f$. This
is a (somewhat skew) generalized version of the CBHD map, since,
if $f$ is the exponential map, then for $\eta\in\g$ we obtain:
$$
B^{\exp}(g,\eta) = \log(g\exp\eta) = \log g + \eta + \CBHD(\log
g,\eta),
$$
with the notation of Section~8. Similarly for right actions.

By definition the map $f$ is equivariant with respect to $B^f$ and
left translations. Since the maps $\Phi_x$ are also equivariant, their
composition $\Phi_x\circ f:\g\to G\to M$ is equivariant, and we have
the following commutative diagram relating the flows on~$M,G$
and~$\g$:
\[
\xymatrix{&\g \ar[r]^f & G \ar[r]^{\Phi_x} & M\\
& \g \ar[u]^{B^f_{e^{t\xi}}} \ar[r]_f
& G \ar[u]^{L_{e^{t\xi}}} \ar[r]_{\Phi_x}
& M \ar[u]^{\Phi_{e^{t\xi}}},}
\]
with the notation of Section~4 for $\exp(t\xi)$. By
Proposition~\ref{pr:equiflows}, this commutative diagram can be
extended to:
\[
\xymatrix{& T\g\approx\g\x\g \ar[r]^{\qquad Tf} & TG \ar[r]^{T\Phi_x}
& TM \\ & \g \ar[u]^{\xi_\g} \ar[r]^f & G \ar[u]^{\xi_G}
\ar[r]^{\Phi_x} & M \ar[u]^{\xi_M(x)} \\ & \g \ar[u]^{B^f_{e^{t\xi}}}
\ar[r]_f & G \ar[u]^{L_{e^{t\xi}}} \ar[r]_{\Phi_x} & M
\ar[u]^{\Phi_{e^{t\xi}}};}
\]
in particular
\begin{equation}
\xi_M(x) \circ \Phi_x \circ f = T\Phi_x \circ Tf \circ \xi_\g.
\label{eq:estas-avisado}
\end{equation}
The overarching question is now: what is the concrete description
of~$\xi_\g$? This we answer next, and we obtain a congenial reply.
Write $g=f(u)$ with $u\in\g$, to distinguish the role of the points
of~$\g$ as coordinates for~$G$. The map $T_uf:T\g\to TG$ can be
factorized into a map from~$T_u\g\approx\g$ to~$\g$, say $A^f_u$, and
the translation~$T_1R_{f(u)}$. Now, in view
of~\eqref{eq:for-later-use}, $Tf\circ\xi_\g=\xi_G\circ f$ gives
$$
T_1R_{f(u)}\circ A^f_u\circ\xi_\g(u) = T_1R_{f(u)}\,\xi;
$$
therefore
\begin{equation}
\xi_\g(u) = (A^f_u)^{-1}\xi,
\label{eq:te-lo-dije}
\end{equation}
where $A^f_u=T_{f(u)}R_{f^{-1}(u)}\circ T_u f$ is the \textit{Darboux
derivative} of~$f$, a map $T\g\to\g\x\g$ yielding the pullback via~$f$
of the right Maurer--Cartan form on~$G$ (a $\g$-valued 1-form
on~$\g$). Then one recovers the `static' version
of~\eqref{eq:madre-del-cordero} from a slightly different viewpoint.
Note the double role of~$\g$ in the construction: on the one hand, its
elements are parameters of the infinitesimal generators on~$M$; on the
other hand they serve as coordinates of the \textit{linear} space on
which we want to solve a differential equation equivalent to the one
originally given on~$M$. The general Darboux derivative for
group-valued maps on manifolds is a key ingredient in the study of
connections via transitive Lie algebroids~\cite{Mackenzie}.

In summary, a differential equation on a homogenous
$G$-manifold~$M$ ---described by infinitesimal generators of the Lie
group action along the projection $\R_t\x M\to M$---  has been
transformed to a `pulled-back' equation on the tangent algebra
of the group:
\begin{equation}
\dot u = \xi_\g(u;t),
\label{eq:mama-de-Tarzan}
\end{equation}
by means of the commutative diagram
\[
\xymatrix{& T\g \ar[r]^{T\Phi_x\circ Tf} & TM \\
& \R_t\x\g \ar[u]^{\xi_\g} \ar[r]^{\!\!\!\!\Phi_x\circ f}
& \R_t\x M \ar[u]^{\xi_M},}
\]
where $\xi_\g$ is the vector field along $\R_t\x\g\to\g$ associated to
the curve $t\mapsto\xi(t)$, explicitly given by
formula~\eqref{eq:te-lo-dije}. The equation evolving on the Lie
algebra is susceptible of attack by geometrical integration
techniques, a point made in~\cite{Engo}.

To exemplify, let us look at Riccati's
equation~\eqref{eq:manes-de-Riccati} again. Consider
$$
L:= a_0(.)E'_+ + a_1(.)H' + a_2(.)E'_- \in
\Map\bigl(\R_t,\mathfrak{sl}(2;\R)\bigr).
$$
We know that if we are able somehow to solve the equation
\begin{equation}
\frac{dg}{dt} = L(t)g(t), \sepword{with} g(t_0) = 1_{SL(2;\R)},
\label{eq:gato-encerradisimo}
\end{equation}
then~\eqref{eq:manes-de-Riccati} is entirely solved by the `Green
operator'
$$
x(t) = \Phi(g(t),x_0).
$$
where~$\Phi$ is the integrated action considered in the previous
section. To search for that solution, let us bring in a variant of the
exponential map~\cite{PesoNormando}. Using canonical coordinates of
the second kind for the element $g(t)\in G$, write:
\begin{align}
g(t) = f(u(t)) &:= \exp(u^0(t)E'_+)\exp(u^1(t)H')\exp(u^2(t)E'_-)
\nonumber \\
&=: \exp(u^0(t)L'_0)\exp(u^1(t)L'_1)\exp(u^2(t)L'_2).
\label{eq:suerte-o-verdad}
\end{align}
Therefore~$f$ denotes the defined locally bijective map
from~$\mathfrak{sl}(2;\R)$ onto~$SL(2;\R)$, with $u\equiv
(u^0,u^1,u^2)$. (Incidentally, this means that we seek the general
solution of the Riccati equation under the form
\begin{equation*}
x(t) = \frac{e^{u^1(t)}x_0}{1 - u^2(t)x_0} + u^0(t);
\end{equation*}
then, taking $x_0=\infty,0,1$, three particular solutions are
obtained, and the reader will see at once that the superposition
formula~\eqref{eq:desmanes-de-LS} follows from here.)

Replacing~$g(t)$ in~\eqref{eq:gato-encerradisimo}
by~\eqref{eq:suerte-o-verdad}, upon using the commutation relations
we obtain
\begin{align*}
\frac{dg(t)}{dt} g^{-1}(t)
&= \dot u^0 E'_+ + \dot u^1 e^{u^0 E'_+}H' e^{-u^0 E'_+}
+ \dot u^2 e^{u^0 E'_+}e^{u^1 H'}E'_ - e^{-u^1 H'} e^{-u^0 E'_+}
\\
&= \dot u^0 E'_+ + \dot u^1 \exp(u^0 \ad E'_+) H' + \dot u^2
\exp(u^0 \ad E'_+) \exp(u^1 \ad H') E'_-
\\
&= \dot u^0 E'_+ + \dot u^1 (H' - u^0 E'_+) +
\dot u^2 e^{-u^1} \exp(u^0 \ad E'_+)E'_-
\\
&= \dot u^0 E'_+ + \dot u^1 (H' - u^0 E'_+)
+ \dot u^2 e^{-u^1} (E'_- - 2u^0 H' + (u^0)^2 E'_+)
\\
&= (\dot u^0 - u^0\dot u^1 + (u^0)^2 e^{-u^1} \dot u^2) E'_+
+ (\dot u^1 - 2u^0 e^{-u^1} \dot u^2) H' + e^{-u^1} \dot u^2 E'_-
\\
&= a_0(t) E'_+ + a_1(t) H' + a_2(t) E'_-.
\end{align*}
This leads to the following differential equations for the
$u$-variables:
\begin{align}
\dot u^0 &= a_0(t) + a_1(t)u^0 + a_2(t)(u^0)^2
\nonumber \\
\dot u^1 &= a_1(t) + 2a_2(t)u^0
\label{eq:mi-patin}
\\
\dot u^2 &= a_2(t)e^{u^1},
\nonumber
\end{align}
to be solved under the initial conditions $u^0(t_0)=u^1(t_0)=
u^2(t_0)=0$. With the chosen map~$f$, the first equation of this
system is the same Riccati equation we started with. This we contrived
to make the point again that \textit{one} particular solution needs to
be known, for the general solution to be obtainable by quadratures.
The explicit form of the Darboux derivative~$A^f_u$ in our example is
$$
(u^0,u^1,u^2; v^0,v^1,v^2) \mapsto
\sum_{k=0}^2\Ad_{\prod_{i=0}^{k-1}\exp(u^iL'_i)}v^kL'_k.
$$
The inversion of this map was just performed, with
result the field corresponding to the system~\eqref{eq:mi-patin}, to
wit,
$$
\bigl(a_0(t) + a_1(t)u^0 + a_2(t)(u^0)^2\bigr)\pd_0 + \bigl(a_1(t)
+ 2a_2(t)u^0\bigr)\pd_1 + a_2(t)e^{u^1}\pd_2,
$$
which is our~\eqref{eq:mama-de-Tarzan}. Once the latter equation is
solved, the rest is obvious: as repeatedly said, one just uses the map
$\Phi\circ f$, to go back to~$M$.

The perceptive reader would ask at this point: what about transferring
the convolution algebra in Section~4 to the tangent algebra, too? We
know nowadays that for conjugation invariant distributions this can be
done~\cite{Doomuch}.

\smallskip

Further work on connections \textit{\`a la} Lie--Rinehart in the
respect of Lie--Scheffers systems is in progress~\cite{NowHector}.

\section{Fa\`a~di~Bruno Hopf algebra and the Lie--Engel theorem}

Due to its fundamental nature, the Hopf algebra we conjure next is
ubiquitous. Let $\Diff^+_0(\R)$ be the group of
orientation-preserving formal diffeomorphisms of~$\R$ (similarly
for~$\C$) leaving~0 fixed. We think of them as exponential power
series:
\begin{equation}
f(t) = \sum_{n=1}^\infty \frac{f_n}{n!}
\,t^n \sepword{with} f_1 > 0. \label{eq:carpe-diem}
\end{equation}
On $\Diff^+_0(\R)$ we consider the coordinate functions
$$
a_n(f) := f_n = f^{(n)}(0), \quad n \geq 1.
$$
Now,
$$
h(t) = \sum_{k=1}^\infty\frac{f_k}{k!}\biggl(\sum_{l=1}^\infty
\frac{g_l}{l!}\,t^l\biggr)^k,
$$
where $h$ is the composition $f\circ g$ of two such diffeomorphisms.
Therefore, from Cauchy's product formula, the $n$th coefficient $h_n
=a_n(h)$ is
$$
h_n = \sum_{k=1}^n \frac{f_k}{k!}\,\sum_{l_i\geq 1,\,l_1+\cdots+l_k=n}
\frac{n!\,g_{l_1}\cdots g_{l_k}}{l_1!\cdots l_k!}.
$$
To rewrite $h_n$ in a compact form, it is convenient to
introduce the notation
$$
\binom{n}{\la;k} := \frac{n!}{\la_1!\la_2!\dots\la_n!
(1!)^{\la_1}(2!)^{\la_2}\dots(n!)^{\la_n}}.
$$
Then, taking in consideration that the sum $l_1+\cdots+l_k=n$ can be
rewritten as
$$
\la_1 + 2\la_2 + \cdots + n\la_n = n, \sepword{where} \la_1 + \cdots +
\la_n = k
$$
if there are $\la_1$ copies of~$1,\,\la_2$ copies of~$2$, and so on,
among the $l_i$; and that the number of contributions from~$g$ of this
type is precisely the multinomial coefficient
$$
\binom{k}{\la_1 \cdots \la_n} = \frac{k!}{\la_1!\cdots\la_n!},
$$
it follows:
\begin{equation}
h_n = \sum_{k=1}^n f_k \sum_{\la\vdash n,|\la|=k}\binom{n}{\la;k}\,
g_1^{\la_1} \dots g_n^{\la_n} =: \sum_{k=1}^n f_k \,
B_{n,k}(g_1,\dots,g_{n+1-k}).
\label{eq:wonder-series}
\end{equation}
We have used notations of the theory of partitions of integers.
The $B_{n,k}$ are called the (partial, exponential) \textit{Bell
polynomials}, often defined via the expansion
$$
\exp\biggl(u\sum_{m\geq 1}g_m\frac{t^m}{m!}\biggr) = 1 + \sum_{n\geq1}
\frac{t^n}{n!}\biggl[\sum_{k=1}^n u^k B_{n,k}(g_1,\dots,g_{n+1-k})
\biggr],
$$
which is a particular case of~\eqref{eq:wonder-series}. Each
$B_{n,k}$ is a homogeneous polynomial of degree~$k$.

According to~\eqref{eq:dies-illa}, a coproduct on
$\Rr(\Diff^+_0(\R))$, which we realize as the polynomial algebra
$\R[a_1,a_2,\dots]$, is given by $\Dl a_n(g,f)=a_n(f\circ g)$. This
entails
\begin{equation}
\Dl a_n = \sum_{k=1}^n \sum_{\la\vdash n,|\la|=k}\binom{n}{\la;k}
a_1^{\la_1} a_2^{\la_2}\dots a_n^{\la_n} \ox a_k.
\label{eq:quam-minimum}
\end{equation}
The flip of~$f$ and~$g$ is done to keep the tradition of writing the
linear part on the right of the tensor product; this amounts to taking
the opposite coalgebra structure. With~\eqref{eq:quam-minimum} we have
a bialgebra structure. In a Hopf algebra grouplike elements are
invertible: $g^{-1}=Sg$. Since $a_1$ is grouplike, to have an antipode
one must either adjoin an inverse $a_1^{-1}$, or put $a_1 = 1$. The
latter is equivalent to work with the subgroup $\Diff^+_{0,1}(\R)$ of
$\Diff^+_0(\R)$, of diffeomorphisms tangent to the identity at~$0$,
that is, to consider power series such that $f_0=0$ and $f_1=1$. The
coproduct formula is accordingly simplified to:
$$
\Dl a_n = \sum_{k=1}^n\sum_{\la\vdash n,|\la|=k}\binom{n}{\la;k}
a_2^{\la_2}a_3^{\la_3}\cdots \ox a_k = \sum_{k=1}^n
B_{n,k}(1,\dots,a_{n+1-k}) \ox a_k.
$$
The resulting graded connected Hopf algebra $\F=\Rr^{\rm cop}
(\Diff^+_{0,1}(\R))$, where the superindex stands for the opposite
coalgebra structure, was baptized Fa\`a~di~Bruno Hopf algebra by Joni
and Rota~\cite{JoniR}. The degree is then given by $|a_n|=n-1$.

Several comments are in order. Formula~\eqref{eq:wonder-series}
can be directly expressed in terms of partitions of finite sets;
consult~\cite{Quaoar} or~\cite[Chapter~5]{Livingstone}. The happy
fact that the algebra of representative functions~$\Rr^{\rm
cop}(\Diff^+_{0,1}(\R))$ is graded is related to the linearity of
the product~$f\circ g$ in one of the coordinates. This also means
that~$\Diff^+_{0,1}(\R)$ is the inverse limit of
finite-dimensional matrix groups, and that it possesses a
(necessarily unique) right invariant connection with vanishing
torsion and curvature. Also, although `formal' may sound a bit
dismissive, one should remember that, in view of E.~Borel's
theorem, expression~\eqref{eq:carpe-diem} does represent a smooth
function; and that~\eqref{eq:wonder-series} can be used to show
without having recourse to complex variables that the composition
of analytic functions on appropriate domains is
analytic~\cite{Krank}.

\smallskip

Let us turn our attention to the dual of the Fa\`a~di~Bruno Hopf
algebra~$\F$. Since we are dealing with a graded connected Hopf
algebra it is natural to consider the \textit{graded dual}, that we
denote simply by~$\F'$; for which $\F''=\F$. (From the discussion in
Sections~4 and~5 we know there exist bigger duals, for instance~$\F'$
does not have grouplike elements apart from its unit~$\eta$.) Let
$a_n'$ be the linear functionals defined by $\dst{a_n'}{P}=\del P/\del
a_n(0)$, where~$P$ is a polynomial in~$\R[a_2,a_3,\dots]$. In
particular the~$a_n'$ kill non-trivial products of the~$a_q$
generators. Also, taking in consideration that the counit
$\eta(P)=P(0)$ of~$\F$ is the unit in~$\F'$
$$
\dst{\Dl a_n'}{P\ox Q} = \dst{a_n'}{m(P\ox Q)} =\dst{a_n'}{PQ} =
\pds{(PQ)}{a_n}(0) = \dst{a_n'\ox 1 + 1 \ox a_n'}{P\ox Q},
$$
Thus the $a_n'$ are primitive. Using the definition of the Bell
polynomials
$$
\dst{a_n'a_m'}{a_q} = \dst{a_n'\ox a_m'}{\Dl a_q} = \begin{cases}
\binom{m+n-1}{n} & \text{if } q = m+n-1 \\ 0 & \text{otherwise}.
\end{cases}
$$
On the other hand, note that
$$
\Dl(a_q a_r) = a_q a_r \ox 1 + 1\ox a_q a_r + a_q \ox a_r +
a_r \ox a_q + R,
$$
where $R$ is either vanishing or a sum of terms of the form $b\ox c$
with~$b$ or~$c$ a monomial in~$a_2,a_3,\dots$ of degree greater
than~$1$. Therefore
$$
\dst{a_n'a_m'}{a_q a_r}  = \dst{a_n' \ox a_m'}{\Dl(a_q a_r)}
=\begin{cases}
1 &\text{if } n = q \neq m = r \text{ or } n = r \neq m = q,
\\
2 &\text{if } m = n = q = r,
\\
0 &\text{otherwise}.
\end{cases}
$$
Similarly, since all the terms of the coproduct of three or more
$a_q$'s are the tensor product of two monomials where at least one of
them is of order greater than~$1$, it follows that
$$
\dst{a'_n a'_m}{a_{q_1}a_{q_2}a_{q_3}\cdots}=0.
$$
Collecting all this together,
$$
a'_na'_m = \binom{m-1+n}{n}a'_{n+m-1} + \bigl(1 +
\dl_{nm}\bigr)(a_na_m)'.
$$
In particular,
$$
[a'_n, a'_m] := a'_na'_m - a'_ma'_n = (m - n)\frac{(n+m-1)!}{n!m!}
\,a_{n+m-1}'.
$$
Therefore, taking $b'_n := (n+1)!a'_{n+1}$, we get the simpler looking
\begin{equation}
[b'_n, b'_m] = (m-n)\,b'_{n+m}.
\label{eq:credula postero}
\end{equation}
The Cartier--Milnor--Moore theorem implies that $\F'$ is the
enveloping algebra of the Lie algebra spanned by the~$b'_n$ with
commutators~\eqref{eq:credula postero}. Obviously $\F'$ can be
realized by the vector fields $Z_n:=x^{n+1}\,\pd/\pd x$, for $n\ge1$,
on the real line~\cite{ConnesMHopf}.

Consider the `regular representation' of~$\F$ given by
$\<a\triangleright a', b>:=\<a', ba>$ on~$\F'$. Since
$$
\<b\triangleright(a \triangleright a'), c> =\<a\triangleright a', cb>
= \<a', cba> =\<ba\triangleright a', c> \sepword{and}
\<1\triangleright a', b> := \<a', b>,
$$
we do obtain a left module algebra over~$\F$. Let now $a$ be a
primitive element of $\F$; using the Sweedler notation:
\begin{align*}
\<a\triangleright b'a', c> &= \<b'a', ca> = \<b'\ox a', \Dl(ca)> =
\<b'\ox a', \Dl c\Dl a>
\\
&= \<b'\ox a', c_{(1)}\ox c_{(2)}(a\ox1 + 1\ox a)>
\\
&= \<a\triangleright b'\ox a'+ b' \ox a\triangleright a', c_{(1)}\ox
c_{(2)}>
\\
&= \<(a\triangleright b')a'+ b'(a\triangleright a'),c>,
\end{align*}
so $a$ acts as a derivation. In particular if $a\triangleright a'
=a\triangleright b'=0$, then $a\triangleright(a'b')=0$, hence the
kernel of the map $a\triangleright\cdot$ is a Lie subalgebra of
vector fields, and we conclude that primitive elements of~$\F$
identify finite-dimensional Lie subalgebras of vector fields. Now,
the space $P(\F)$ of primitive elements of~$\F$ has
\textit{dimension~two}. Indeed, $P(\F)=\bigl(\R
1\oplus{\F'_+}^2\bigr)^\perp$, where $\F'_+:=\ker\eta$ is the
augmentation ideal of~$\F'$. By~\eqref{eq:credula postero} there
is a dual basis of~$\F'$ made of products, except for its first
two elements. Hence $\dim P(\F) = 2$. A basis of $P(\F)$ is given
by $\{a_2,a_3-\frac{3}{2}a_2^2\}$. This yields the equations
$y''=0$ and $y'y'''-3(y'')^2/2=0$, respectively solved by
dilations and by the action of~$SL(2;\R)$ we know; translations do
not show up because we made $a_1=1$.

The previous argument, together with the part of classical
one~\cite{LieAngel} ---more recently rehearsed
in~\cite[Section~XIX]{DieudonneIV} or in~\cite{Luther}--- to the
effect that infinitesimally transitive actions on the line must
correspond to Lie algebras of vector fields of dimension at most
three, shows that Riccati's is the only nonlinear Lie--Scheffers
differential equation on the real (or complex) line. Whether or not it
is simpler to think in Hopf algebraic terms seems largely a matter of
taste. We do contend that the Beatus Fa\`a~di~Bruno algebra is too
fundamental an object to ignore.

\subsection*{Acknowledgments}

We are very thankful to J.~C.~V\'arilly for helpful comments. The
first named author acknowledges partial financial support through
grants BFM-2003-02532 and DGA-Grupos Consolidados~E24/1. KEF greatly
acknowledges support from the European Post-Doctoral Institute. HF
acknowledges support from the Vicerrector\'{\i}a de Investigaci\'on of
the Universidad de Costa Rica. The first and fourth named authors are
grateful for discussions with J.~Grabowski and G.~Marmo, during a
visit to Universit\`a di Napoli Federico II in November~2005. JMGB
moreover acknowledges partial support from CICyT, Spain, through
grant~FIS2005-02309. KEF, HF and JMGB thank the Departamento de
F\'{\i}sica Te\'orica of the Universidad de Zaragoza for its warm
hospitality.

\end{document}